\documentclass{amsart}
\usepackage{amsfonts}

\usepackage{amsmath}

\newtheorem{theorem}{Theorem}
\theoremstyle{plain}
\newtheorem{acknowledgement}{Acknowledgement}

\newtheorem{corollary}{Corollary}

\newtheorem{definition}{Definition}
\newtheorem{example}{Example}

\newtheorem{lemma}{Lemma}

\newtheorem{proposition}{Proposition}
\newtheorem{remark}{Remark}

\numberwithin{equation}{section}
\input{tcilatex}

\begin{document}
\title[Stochastic $\Sigma $-convergence]{Stochastic $\Sigma $-convergence
and applications}
\author{Mamadou Sango}
\address{M. Sango, Department of Mathematics and Applied Mathematics,
University of Pretoria, Pretoria 0002, South Africa}
\email{mamadou.sango@up.ac.za}
\author{Jean Louis Woukeng}
\address{J.L. Woukeng, Department of Mathematics and Computer Science,
University of Dschang, P.O. Box 67, Dschang, Cameroon}
\curraddr{J.L. Woukeng, Department of Mathematics and Applied Mathematics,
University of Pretoria, Pretoria 0002, South Africa}
\email{jwoukeng@yahoo.fr}
\date{April, 2011}

\begin{abstract}
Motivated by the fact that in nature almost all phenomena behave randomly in
some scales and deterministically in some other scales, we build up a
framework suitable to tackle both deterministic and stochastic
homogenization problems simultaneously, and also separately. Our approach,
the \textit{stochastic }$\Sigma $\textit{-convergence}, can be seen either
as a multiscale stochastic approach since deterministic homogenization
theory can be seen as a special case of stochastic homogenization theory
(see Theorem \ref{t2.2}), or as a conjunction of the stochastic and
deterministic approaches, both taken globally, but also each separately. One
of the main applications of our results is the homogenization of a model of
rotating fluids.
\end{abstract}

\maketitle

\section{Introduction}

A wide range of scientific and engineering problems involve multiscale
phenomena. Roughly speaking, each matter is characterized by its own
geometric dimensions which are very often several order of magnitude larger.
The study and the understanding of these issues demand the development of
new mathematical tools and methods. Homogenization theory is such a tool
which now occupies a central place in contemporary mathematical research.

Deterministic problems in the periodic setting prominently featured in the
first decade of the development of the theory till the pioneering works of
Kozlov \cite{Kozlov1, Kozlov2}, Papanicolaou and Varadhan \cite{PAPANICOLAOU}
in stochastic homogenization in the late 1970s. Since then intense research
activities have been undertaken with a great wealth of results as shown by
the vast existing literature to date, see e.g., \cite{BBL, Bourgeat1,
Bourgeat2, 10, 14, 15, 28, Pankov1, Nils1, Nils2, Wright1, Wright2, Wright3,
Zhikov2, Zhikov3, Zhikov4, Zhikov5}. It is worth noting the interesting work
on stochastic homogenization in the framework of viscosity solutions by
several prominent mathematicians \cite{CAFFARELLI1, CAFFARELLI2, VARADHAN1,
LIONS-SOUGANIDIS1, LIONS-SOUGANIDIS2}.

In order to deal with deterministic homogenization theory beyond the
periodic setting, Nguetseng \cite{26}, following Zhikov and Krivenko \cite%
{Zhikov4}, introduced the concept of homogenization algebras. This theory
relies heavily on ergodic theory (but not the ergodicity!) because in
applications, the assumption of ergodicity of the homogenization algebra
considered is fundamental.\emph{\ It is important to note that there was a
gap between the periodic homogenization theory and the stochastic
homogenization theory, gap which was filled by Nguetseng's deterministic
homogenization theory. However as we will see in the present work, this
recent deterministic theory can be viewed as a special case of a generalized
version of the stochastic homogenization theory of Bourgeat et al. \cite{10}
which we construct.} Indeed, Theorem \ref{t2.2} (see Section 2) allows to
build on the spectrum of an algebra with mean value, a dynamical system
whose invariant measure is precisely the measure related to the mean value
defined on the algebra. As a result of the above-mentioned theorem, we get a
generalization of all the results presented in \cite{26, CMP}, those in \cite%
{26, CMP} being the special case corresponding to ergodic algebras; see
Section 4.

The two theories mentioned above have the specificity to be used to solve
either stochastic homogenization problems only (for the first one) or
deterministic homogenization problems only (for the second one).
Unfortunately, as we know, in nature, very few phenomena behave, either just
randomly or deterministically; most of these phenomena behave randomly in
some scales, and deterministically in other scales.

Motivated by this vision of the physical nature, we rely on these two
theories and hence on their associated convergence methods (the stochastic
two scale convergence in the mean \cite{10} and the $\Sigma $-convergence %
\cite{26, CMA}) to propose a general method of solving coupled -
deterministic and stochastic - homogenization problems. Our method, the 
\textit{stochastic }$\Sigma $\textit{-convergence}, combines the macroscopic
and microscopic [random and deterministic] scales, and has therefore the
advantage of taking both the simplicity and the efficiency of the
macroscopic models, as well as the accuracy of the coupled
random-deterministic microscopic models. Moreover our multiscale approach is
motivated by the fact that the usual monoscale approach has proven to be
inadequate because of prohibitively large number of variables involved in
each physical problem. One can also give at least two reasons quite natural.
Firstly, a scale can not be at the same time deterministic and random.
Secondly, the application of our results to natural phenomena; see Sections
5 and 6. To be more precise, our method permits henceforth to treat
deterministic homogenization problems without resorting to the ergodicity
assumption on one hand, and on the other hand allow viewing the stochastic
two-scale convergence in the mean \cite{10} in a more general angle as
generalizing the $\Sigma $-convergence \cite{26, CMP}.

We hope that the theory developed in the present paper will find
applications in the emerging field of homogenization of stochastic partial
differential equations undertaken in the papers \cite{BENSOUSSAN},\ \cite%
{WANG} in the periodic case and in \cite{SANGO} in the case of non
periodically perforated domains.

The paper is organized as follows. In Section 2 we give some preliminary
results related to the theory of dynamical systems on abstract probability
spaces. We also define and give some fundamental properties of generalized
Besicovitch spaces. Section 3 is devoted to the study of the concept of
stochastic $\Sigma $-convergence. We prove therein some compactness results.
In Sections 4, 5 and 6, we give several applications of the earlier results.
We begin in Section 4 by showing how the results of Section 3 apply and how
they generalize the existing results; this is illustrated by the study of a
rather simple linear operator in divergence form. We then compare our
results with the already existing ones. In Section 5 we study the
homogenization problem for the well-known nonlinear Reynolds equation. One
important achievement of our results is obtained in Section 6 where we solve
the coupled stochastic-deterministic homogenization problem related to the
following Stokes equation: 
\begin{equation*}
\begin{array}{l}
-\overset{N}{\underset{i,j=1}{\sum }}\frac{\partial }{\partial x_{i}}\left(
a_{ij}(x,T(x/\varepsilon _{1})\omega ,x/\varepsilon _{2})\frac{\partial 
\mathbf{u}_{\varepsilon }}{\partial x_{j}}\right) +\mathbf{h}^{\varepsilon
}\times \mathbf{u}_{\varepsilon }+\func{grad}p_{\varepsilon }=\mathbf{f}%
\text{\ in }Q \\ 
\;\;\;\;\;\;\;\;\;\;\;\;\;\;\;\;\;\;\;\;\;\;\;\;\;\;\;\;\;\;\;\;\;\func{div}%
\mathbf{u}_{\varepsilon }=0\text{\ in }Q \\ 
\;\;\;\;\;\;\;\;\;\;\;\;\;\;\;\;\;\;\;\;\;\;\;\;\;\;\;\;\;\;\;\;\;\mathbf{u}%
_{\varepsilon }=0\text{\ on }\partial Q.%
\end{array}%
\end{equation*}%
We get the following homogenization result which is, to our knowledge, new.

\begin{theorem}
Assume

\begin{enumerate}
\item $a_{ij}(x,\omega ,\cdot )\in A$ for all $(x,\omega )\in Q\times \Omega 
$, $1\leq i,j\leq N$, and $\mathbf{h}\in L^{\infty }(\Omega ;A)^{N}$.
\end{enumerate}

For each $0<\varepsilon <1$ and for a.e. $\omega \in \Omega $ let $\mathbf{u}%
_{\varepsilon }(\cdot ,\omega )=(u_{\varepsilon }^{k}(\cdot ,\omega ))\in 
\mathbb{H}_{0}^{1}(Q)$ be the (unique) solution of the above Stokes
equation. Then as $\varepsilon \rightarrow 0$, 
\begin{equation*}
\mathbf{u}_{\varepsilon }\rightarrow \mathbf{u}_{0}\text{ stoch. in }%
L^{2}(Q\times \Omega )^{N}\text{-weak}
\end{equation*}%
and 
\begin{equation*}
\frac{\partial u_{\varepsilon }^{k}}{\partial x_{j}}\rightarrow \frac{%
\partial u_{0}^{k}}{\partial x_{j}}+\overline{D}_{j,\omega }u_{1}^{k}+\frac{%
\overline{\partial }u_{2}^{k}}{\partial y_{j}}\text{ stoch. in }%
L^{2}(Q\times \Omega )\text{-weak }\Sigma \text{ }(1\leq j,k\leq N)
\end{equation*}%
where $\mathbf{u}=(\mathbf{u}_{0},\mathbf{u}_{1},\mathbf{u}_{2})$ is the
unique solution to the following variational problem: 
\begin{equation*}
a(\mathbf{u},\mathbf{v})+\iint_{Q\times \Omega }(\widetilde{\mathbf{h}}%
\times \mathbf{u}_{0})\cdot \mathbf{v}_{0}dxd\mu =\left\langle \mathbf{f},%
\mathbf{v}_{0}\right\rangle \text{ for all }\mathbf{v}=(\mathbf{v}_{0},%
\mathbf{v}_{1},\mathbf{v}_{2})\in \mathbb{F}_{0}^{1}
\end{equation*}%
with: 
\begin{eqnarray*}
a(\mathbf{u},\mathbf{v}) &=&\sum_{i,j,k=1}^{N}\iint_{Q\times \Omega \times
\Delta (A)}\widehat{a}_{ij}(x,\omega ,s)\left( \frac{\partial u_{0}^{k}}{%
\partial x_{j}}+\overline{D}_{j,\omega }u_{1}^{k}+\partial _{j}\widehat{%
u_{2}^{k}}\right) \\
&&\;\;\;\;\ \ \ \ \ \ \ \ \times \left( \frac{\partial v_{0}^{k}}{\partial
x_{i}}+\overline{D}_{i,\omega }v_{1}^{k}+\partial _{i}\widehat{v_{2}^{k}}%
\right) dxd\mu d\beta ;
\end{eqnarray*}%
\begin{equation*}
\widetilde{\mathbf{h}}(\omega )=\int_{\Delta (A)}\widehat{\mathbf{h}}(\omega
,s)d\beta ;\ \ \ \ \ \ \ \ \ \ \ \ \ \ \ \ \ \ \ \ \ \ \ \ \ \ \ \ \ \ \ \ \ 
\end{equation*}%
\begin{equation*}
\left\langle \mathbf{f},\mathbf{v}_{0}\right\rangle =\int_{\Omega }\left( 
\mathbf{f}(\cdot ,\omega ),\mathbf{v}_{0}(\cdot ,\omega )\right)
_{H^{-1}(Q)^{N},H_{0}^{1}(Q)^{N}}d\mu
\end{equation*}%
and%
\begin{equation*}
\partial _{j}\widehat{u_{2}^{k}}=\mathcal{G}_{1}(\overline{\partial }%
u_{2}^{k}/\partial y_{j})\text{ (and a same definition for }\partial _{i}%
\widehat{v_{2}^{k}}\text{).}
\end{equation*}
\end{theorem}

Unless otherwise specified, vector spaces throughout are assumed to be
complex vector spaces, and scalar functions are assumed to be complex
valued. We shall always assume that the numerical spaces $\mathbb{R}^{m}$
and their open sets are each equipped with the Lebesgue measure.

\section{Preliminaries on dynamical systems and generalized Besicovitch
spaces}

\subsection{Stochastic vector calculus}

We begin \ by recalling the definition of the notion of a dynamical system.
Let $(\Omega ,\mathcal{M},\mu )$ denote a probability space. An $N$%
-dimensional dynamical system on $\Omega $ is a family of invertible
mappings $T\left( x\right) :\Omega \rightarrow \Omega $, $x\in \mathbb{R}%
^{N} $, such that the following conditions hold:

\begin{itemize}
\item[(i)] (\textit{Group property}) $T\left( 0\right) =id_{\Omega }$ and $%
T\left( x+y\right) =T(x)\circ T(y)$ for all $x,y\in \mathbb{R}^{N}$;

\item[(ii)] (\textit{Invariance}) The mappings $T\left( x\right) :\Omega
\rightarrow \Omega $ are measurable and $\mu $-measure preserving, i.e., $%
\mu \left( T\left( x\right) F\right) =\mu \left( F\right) $ for each $x\in 
\mathbb{R}^{N}$ and every $F\in \mathcal{M}$;

\item[(iii)] (\textit{Measurability}) For each $F\in \mathcal{M}$, the set $%
\left\{ \left( x,\omega \right) \in \mathbb{R}^{N}\times \Omega :T\left(
x\right) \omega \in F\right\} $ is measurable with respect to the product $%
\sigma $-algebra $\mathcal{L}\otimes \mathcal{M}$, where $\mathcal{L}$ is
the $\sigma $-algebra of Lebesgue measurable sets.
\end{itemize}

We recall that in (i) above, the symbol $\circ $ denotes the usual
composition of mappings, and in (iii), $\mathcal{L}\otimes \mathcal{M}$ is
the $\sigma $-algebra generated by the family $\{L\times M:L\in \mathcal{L}$
and $M\in \mathcal{M}\}$, $L\times M$ being the Cartesian product of the
sets $L$ and $M$.

If $\Omega $ is a compact topological space, by a continuous $N$-dimensional
dynamical system on $\Omega $ is meant any family of mappings $T(x):\Omega
\rightarrow \Omega $, $x\in \mathbb{R}^{N}$, satisfying the above group
property (i) and the following condition: The mapping $(x,\omega )\mapsto
T(x)\omega $ is continuous from $\mathbb{R}^{N}\times \Omega $ to $\Omega $.

Let $1\leq p\leq \infty $. An $N$-dimensional dynamical system $T\left(
x\right) :\Omega \rightarrow \Omega $ induces a $N$-parameter group of
isometries $U\left( x\right) :L^{p}(\Omega )\rightarrow L^{p}(\Omega )$
defined by 
\begin{equation*}
\left( U\left( x\right) f\right) (\omega )=f\left( T(x)\omega \right) \text{%
,\ \ }f\in L^{p}(\Omega )
\end{equation*}%
which is strongly continuous, i.e., $U\left( x\right) f\rightarrow f$\ in $%
L^{p}(\Omega )$ as $x\rightarrow 0$; see \cite[p. 223]{20} or\ \cite[p. 131]%
{28}. We denote by $D_{i,p}$ ($1\leq i\leq N$) the generator of $U(x)$ along
the $i$th coordinate direction, and by $\mathcal{D}_{i,p}$ its domain. Thus,
for $f\in L^{p}(\Omega )$, $f$ is in $\mathcal{D}_{i,p}$ if and only if the
limit $D_{i,p}f$ defined by 
\begin{equation*}
D_{i,p}f(\omega )=\lim_{\tau \rightarrow 0}\frac{f\left( T(\tau e_{i})\omega
\right) -f(\omega )}{\tau }
\end{equation*}%
exists strongly in $L^{p}(\Omega )$, where $e_{i}$ denotes the vector $%
\left( \delta _{ij}\right) _{1\leq j\leq N}$, $\delta _{ij}$ being the
Kronecker $\delta $. One can naturally define higher order derivatives by
setting $D_{p}^{\alpha }=D_{1,p}^{\alpha _{1}}\cdot \cdot \cdot
D_{N,p}^{\alpha _{N}}$ for $\alpha =(\alpha _{1},...,\alpha _{N})\in \mathbb{%
N}^{N}$, where $D_{i,p}^{\alpha _{i}}=D_{i,p}\circ \cdots \circ D_{i,p}$, $%
\alpha _{i}$-times.

Now we need to define the stochastic analog of the smooth functions on $%
\mathbb{R}^{N}$. To this end, we set $\mathcal{D}_{p}(\Omega )=\cap
_{i=1}^{N}\mathcal{D}_{i,p}$ and define 
\begin{equation*}
\mathcal{D}_{p}^{\infty }(\Omega )=\left\{ f\in L^{p}\left( \Omega \right)
:D_{p}^{\alpha }f\in \mathcal{D}_{p}(\Omega )\text{ for all }\alpha \in 
\mathbb{N}^{N}\right\} .
\end{equation*}%
It is a fact that each element of $\mathcal{D}_{\infty }^{\infty }(\Omega )$
possesses stochastic derivatives of any order that are bounded. So as in %
\cite{4} we denote it by the suggestive symbol $\mathcal{C}^{\infty }(\Omega
)$, and also as in \cite{4} it can be shown that $\mathcal{C}^{\infty
}(\Omega )$ is dense in $L^{p}(\Omega )$, $1\leq p<\infty $. At this level,
one can naturally define the concept of \textit{stochastic distribution}: by
a stochastic distribution on $\Omega $ is meant any continuous linear
mapping from $\mathcal{C}^{\infty }(\Omega )$ to the complex field $\mathbb{C%
}$. We recall that $\mathcal{C}^{\infty }(\Omega )$ is endowed with its
natural topology defined by the family of seminorms $N_{n}(f)=\sup_{\left%
\vert \alpha \right\vert \leq n}\sup_{\omega \in \Omega }\left\vert
D_{\infty }^{\alpha }f(\omega )\right\vert $ (where $\left\vert \alpha
\right\vert =\alpha _{1}+...+\alpha _{N}$ for $\alpha =(\alpha
_{1},...,\alpha _{N})\in \mathbb{N}^{N}$). We denote the space of stochastic
distributions by $\left( \mathcal{C}^{\infty }(\Omega )\right) ^{\prime }$.
One can also define the stochastic weak derivative of $f\in \left( \mathcal{C%
}^{\infty }(\Omega )\right) ^{\prime }$ as follows: For any $\alpha \in 
\mathbb{N}^{N}$, $D^{\alpha }f$ stands for the stochastic distribution
defined by 
\begin{equation*}
\left( D^{\alpha }f\right) (\phi )=(-1)^{\left\vert \alpha \right\vert
}f\left( D_{\infty }^{\alpha }\phi \right) \;\;\forall \phi \in \mathcal{C}%
^{\infty }(\Omega ).
\end{equation*}%
As $\mathcal{C}^{\infty }(\Omega )$ is dense in $L^{p}(\Omega )$ ($1\leq
p<\infty $), it is immediate that $L^{p}(\Omega )\subset \left( \mathcal{C}%
^{\infty }(\Omega )\right) ^{\prime }$ so that one may define the stochastic
weak derivative of any $f\in L^{p}(\Omega )$, and it verifies the following
functional equation: 
\begin{equation*}
\left( D^{\alpha }f\right) (\phi )=(-1)^{\left\vert \alpha \right\vert
}\int_{\Omega }fD_{\infty }^{\alpha }\phi d\mu \text{\ for all }\phi \in 
\mathcal{C}^{\infty }(\Omega )\text{.}
\end{equation*}%
In particular, for $f\in \mathcal{D}_{i,p}$ we have $-\int_{\Omega
}fD_{i,\infty }\phi d\mu =\int_{\Omega }\phi D_{i,p}fd\mu $ for all $\phi
\in \mathcal{C}^{\infty }(\Omega )$ so that we may identify $D_{i,p}f$ with $%
D^{\alpha _{i}}f$, where $\alpha _{i}=\left( \delta _{ij}\right) _{1\leq
j\leq N}$. Conversely, if $f\in L^{p}\left( \Omega \right) $ is such that
there exists $f_{i}\in L^{p}\left( \Omega \right) $ with $\left( D^{\alpha
_{i}}f\right) (\phi )=-\int_{\Omega }f_{i}\phi d\mu $ for all $\phi \in 
\mathcal{C}^{\infty }(\Omega )$, then $f\in \mathcal{D}_{i,p}$ and $%
D_{i,p}f=f_{i}$. Therefore, endowing $\mathcal{D}_{p}(\Omega )$ with the
natural graph norm 
\begin{equation*}
\left\Vert f\right\Vert _{\mathcal{D}_{p}(\Omega )}^{p}=\left\Vert
f\right\Vert _{L^{p}\left( \Omega \right) }^{p}+\sum_{i=1}^{N}\left\Vert
D_{i,p}f\right\Vert _{L^{p}\left( \Omega \right) }^{p}\;\;(f\in \mathcal{D}%
_{p}(\Omega ))
\end{equation*}%
we obtain a Banach space representing the stochastic generalization of the
Sobolev spaces $W^{1,p}\left( \mathbb{R}^{N}\right) $, and so, we denote it
by $W^{1,p}(\Omega )$.

Now, returning to the general setting of dynamical systems, we recall that a
function $f\in L^{p}\left( \Omega \right) $ is said to be \textit{invariant
for }$T$ (relative to $\mu $) if for any $x\in \mathbb{R}^{N}$, $f\circ
T\left( x\right) =f$ $\mu $-a.e. on $\Omega $. We denote by $%
I_{nv}^{p}\left( \Omega \right) $ the set of functions in $L^{p}\left(
\Omega \right) $ that are invariant for $T$. The set $I_{nv}^{p}\left(
\Omega \right) $ is a closed vector subspace of $L^{p}\left( \Omega \right) $%
. The dynamical system $T$ is said to be ergodic if every $T$-invariant
function $f\in I_{nv}^{p}\left( \Omega \right) $ is constant. We have the
following very useful properties for functions in $L^{1}\left( \Omega
\right) $.

\begin{itemize}
\item[(P1)] For $f\in \mathcal{D}_{1}^{\infty }\left( \Omega \right) $, and
for $\mu $-a.e. $\omega \in \Omega $, the function $x\mapsto f(T(x)\omega )$
is in $\mathcal{C}^{\infty }(\mathbb{R}^{N})$ and further $D_{x}^{\alpha
}f(T(x)\omega )=\left( D_{1}^{\alpha }f\right) (T(x)\omega )$ for any $%
\alpha \in \mathbb{N}^{N}$.

\item[(P2)] For $f\in L^{1}\left( \Omega \right) $, we have $f\in
I_{nv}^{1}\left( \Omega \right) $ if and only if $D_{i,1}f=0$ for each $%
1\leq i\leq N$.
\end{itemize}

Let $1<p<\infty $. Thanks to (P2) above, one can easily check that, for $%
f\in L^{p}\left( \Omega \right) $, $f$ is in $I_{nv}^{p}\left( \Omega
\right) $ if and only if $D_{i,p}f=0$ for all $1\leq i\leq N$, since $%
D_{i,p} $ is the restriction to $L^{p}(\Omega )$ of $D_{i,1}$. So if we
endow $\mathcal{C}^{\infty }(\Omega )$ with the seminorm 
\begin{equation}
\left\| u\right\| _{\#,p}^{p}=\sum_{i=1}^{N}\left\| D_{i,p}u\right\|
_{L^{p}\left( \Omega \right) }^{p}\ \ \ \ \ (u\in \mathcal{C}^{\infty
}(\Omega ))  \label{2.0}
\end{equation}%
we obtain a locally convex space which is\ generally non separated and non
complete. We denote by $\mathcal{W}^{1,p}(\Omega )$ the separated completion
of $\mathcal{C}^{\infty }(\Omega )$ with respect to the seminorm $\left\|
\cdot \right\| _{\#,p}$, and we denote by $I_{p}$ the canonical mapping of $%
\mathcal{C}^{\infty }(\Omega )$ into its separated completion $\mathcal{W}%
^{1,p}(\Omega )$. It is to be noted that $\mathcal{W}^{1,p}(\Omega )$ is
also the separated completion of $\mathcal{C}^{\infty }(\Omega
)/(I_{nv}^{p}(\Omega )\cap \mathcal{C}^{\infty }(\Omega ))$ with respect to
the same seminorm since for $u\in \mathcal{C}^{\infty }(\Omega )$ we have $%
\left\| u\right\| _{\#,p}=0$ if and only if $u\in I_{nv}^{p}(\Omega )$, that
is $u\in I_{nv}^{p}(\Omega )\cap \mathcal{C}^{\infty }(\Omega )$. The
following property is obtained through the theory of completion of uniform
spaces; see, e.g., \cite[Chap. II, Sect. 3, no 7]{9}.

\medskip

The gradient operator $D_{\omega ,p}=(D_{1,p},...,D_{N,p}):\mathcal{C}%
^{\infty }(\Omega )\rightarrow L^{p}\left( \Omega \right) ^{N}$ extends by
continuity to a unique mapping $\overline{D}_{\omega ,p}=(\overline{D}%
_{1,p},...,\overline{D}_{N,p}):\mathcal{W}^{1,p}(\Omega )\rightarrow
L^{p}\left( \Omega \right) ^{N}$ with the properties 
\begin{equation*}
D_{i,p}=\overline{D}_{i,p}\circ I_{p}\ \ \ \ \ \ \ \ \ 
\end{equation*}%
and 
\begin{equation*}
\left\| u\right\| _{\mathcal{W}^{1,p}(\Omega )}\equiv \left\| u\right\|
_{\#,p}=\left( \sum_{i=1}^{N}\left\| \overline{D}_{i,p}u\right\|
_{L^{p}\left( \Omega \right) }^{p}\right) ^{1/p}\text{\ for }u\in \mathcal{W}%
^{1,p}(\Omega ).
\end{equation*}%
Moreover, the mapping $\overline{D}_{\omega ,p}$ is an isometric embedding
of $\mathcal{W}^{1,p}(\Omega )$ into a closed subspace of $L^{p}\left(
\Omega \right) ^{N}$, so that the Banach space $\mathcal{W}^{1,p}(\Omega )$
is reflexive. By duality we define the operator div$_{\omega ,p^{\prime
}}:L^{p^{\prime }}\left( \Omega \right) ^{N}\rightarrow \left( \mathcal{W}%
^{1,p}(\Omega )\right) ^{\prime }$ ($p^{\prime }=p/(p-1)$) by 
\begin{equation*}
\left\langle \text{div}_{\omega ,p^{\prime }}u,w\right\rangle =-\left\langle
u,\overline{D}_{\omega ,p}w\right\rangle \text{\ for all }w\in \mathcal{W}%
^{1,p}\left( \Omega \right) \text{ and }u=(u_{i})\in L^{p^{\prime }}(\Omega
)^{N}\text{,}
\end{equation*}%
where $\left\langle u,\overline{D}_{\omega ,p}w\right\rangle
=\sum_{i=1}^{N}\int_{\Omega }u_{i}\overline{D}_{i,p}wd\mu $. The operator div%
$_{\omega ,p^{\prime }}$ just defined extends the natural divergence
operator defined in $\mathcal{C}^{\infty }(\Omega )$ since for all $f\in 
\mathcal{C}^{\infty }(\Omega )$ we have $D_{i,p}f=\overline{D}%
_{i,p}(I_{p}(f))$.

The following result will be of great interest in the next sections.

\begin{proposition}
\label{p2.1}Let $\mathbf{v}\in L^{p}\left( \Omega \right) ^{N}$ satisfying 
\begin{equation*}
\int_{\Omega }\mathbf{v}\cdot \mathbf{g}d\mu =0\text{ for all }\mathbf{g}\in 
\mathcal{V}_{\func{div}}=\{\mathbf{f}\in \mathcal{C}^{\infty }(\Omega )^{N}:%
\text{\emph{div}}_{\omega ,p^{\prime }}\mathbf{f}=0\}.
\end{equation*}%
Then there exists $u\in \mathcal{W}^{1,p}(\Omega )$ such that $\mathbf{v}=%
\overline{D}_{\omega ,p}u$.
\end{proposition}

\begin{proof}
We need to check the following: (1) div$_{\omega ,p^{\prime }}$ is closed;
(2) $($div$_{\omega ,p^{\prime }})^{\ast }=-\overline{D}_{\omega ,p}$ where $%
($div$_{\omega ,p^{\prime }})^{\ast }$ is the adjoint operator of div$%
_{\omega ,p^{\prime }}$; (3) Ran($\overline{D}_{\omega ,p}$) is closed in $%
L^{p}\left( \Omega \right) ^{N}$ and finally, (4) $\mathbf{v}$ is orthogonal
to the kernel of div$_{\omega ,p^{\prime }}$. Indeed (1)-(3) will yield Ran$(%
\overline{D}_{\omega ,p})=(\ker ($div$_{\omega ,p^{\prime }}))^{\perp }$ by
a well-known result (see, e.g., \cite[Chap. 13, p. 352, Thm 13.8]{30}) where 
$(\ker ($div$_{\omega ,p^{\prime }}))^{\perp }$ denote the orthogonal
complement of $\ker ($div$_{\omega ,p^{\prime }})$, and finally the
proposition will follow at once from (4). So let us check them.

(1) is trivial, (2) is a mere consequence of the definition of div$_{\omega
,p^{\prime }}$. As for (3), if $\mathbf{v}_{n}=\overline{D}_{\omega
,p}u_{n}\in \,$Ran$(\overline{D}_{\omega ,p})$ is such that $\mathbf{v}%
_{n}\rightarrow \mathbf{v}$ in $L^{p}(\Omega )^{N}$, then $(u_{n})_{n}$ is a
Cauchy sequence in $\mathcal{W}^{1,p}(\Omega )$ and so, converges in $%
\mathcal{W}^{1,p}(\Omega )$ towards some $u\in \mathcal{W}^{1,p}(\Omega )$,
that is, $\overline{D}_{\omega ,p}u_{n}\rightarrow \overline{D}_{\omega ,p}u$
in $L^{p}\left( \Omega \right) ^{N}$, hence $\mathbf{v}=\overline{D}_{\omega
,p}u$. Finally for (4) it suffices to show that $\mathcal{V}_{\func{div}}$
is dense in $\ker ($div$_{\omega ,p^{\prime }})$. To see this, let $\mathbf{g%
}\in \ker ($div$_{\omega ,p^{\prime }})$; arguing as in the proof of %
\cite[Lemma 2.3 (b)]{10} there exists a sequence $(\mathbf{g}%
_{n})_{n}\subset \mathcal{V}_{\func{div}}$ such that $\mathbf{g}%
_{n}\rightarrow \mathbf{g}$ in $L^{p}\left( \Omega \right) ^{N}$. The proof
is complete.
\end{proof}

We end this subsection with some definitions. Let $f$ be a measurable
function in $\Omega $; for a fixed $\omega \in \Omega $ the function $%
x\mapsto f(T(x)\omega )$, $x\in \mathbb{R}^{N}$, is called a realization of $%
f$ and the mapping $(x,\omega )\mapsto f(T(x)\omega )$ is called a
stationary process. The process is said to be stationary ergodic if the
dynamical system $T$ is ergodic. We will also use the notation div$_{\omega
} $ instead of div$_{\omega ,p^{\prime }}$, accordingly.

In the forthcoming sections we will adopt the following notation: $\overline{%
D}_{\omega }$ will stand for $\overline{D}_{\omega ,p}$, and, $\overline{D}%
_{i,p}$ (resp. $D_{i,p}$) will be denoted by $\overline{D}_{i,\omega }$
(resp. $D_{i,\omega }$) if there is no danger of confusion.

\subsection{Homogenization supralgebras}

We use a new concept of homogenization algebras. This concept has just been
defined in a more recent paper \cite{CMP}. It is more general than those
defined in the papers \cite{26, Zhikov4} because we do not need the algebra
to be separable (as in \cite{26}), or to consist of functions that are
uniformly continuous (as in \cite{Zhikov4}). Before we go any further, we
need to give some preliminaries. Let $\mathcal{H}=(H_{\varepsilon
})_{\varepsilon >0}$ be the action of $\mathbb{R}_{+}^{\ast }$ (the
multiplicative group of positive real numbers) on the numerical space $%
\mathbb{R}^{N}$ defined as follows: 
\begin{equation}
H_{\varepsilon }(x)=\frac{x}{\varepsilon _{1}}\;\;(x\in \mathbb{R}^{N})
\label{2.1}
\end{equation}%
where $\varepsilon _{1}$ is a positive function of $\varepsilon $ tending to
zero with $\varepsilon $. For given $\varepsilon >0$, let 
\begin{equation*}
u^{\varepsilon }(x)=u(H_{\varepsilon }(x))\;\;(x\in \mathbb{R}%
^{N}).\;\;\;\;\;\;
\end{equation*}%
For $u\in L_{\text{loc}}^{1}(\mathbb{R}_{y}^{N})$ (as usual, $\mathbb{R}%
_{y}^{N}$ denotes the numerical space $\mathbb{R}^{N}$ of variables $%
y=(y_{1},...,y_{N})$), $u^{\varepsilon }$ lies in $L_{\text{loc}}^{1}(%
\mathbb{R}_{x}^{N})$. More generally, if $u$ lies in $L_{\text{loc}}^{p}(%
\mathbb{R}^{N})$ (resp. $L^{p}(\mathbb{R}^{N})$), $1\leq p<+\infty $, then
so also does $u^{\varepsilon }$.

A function $u\in \mathcal{B}(\mathbb{R}_{y}^{N})$ (the $\mathcal{C}$%
*-algebra of bounded continuous complex functions on $\mathbb{R}_{y}^{N}$)
is said to have a mean value for $\mathcal{H}$, if there exists a complex
number $M(u)$ such that $u^{\varepsilon }\rightarrow M(u)$ in $L^{\infty }(%
\mathbb{R}_{x}^{N})$-weak $\ast $ as $\varepsilon \rightarrow 0$. The
complex number $M(u)$ is called the mean value of $u$ (for $\mathcal{H}$).
It is evident that this defines a mapping $M$ which is a positive linear
form (on the space of functions $u\in \mathcal{B}(\mathbb{R}_{y}^{N})$ with
mean value) attaining the value $1$ on the constant function $1$ and
verifying the inequality $\left\vert M(u)\right\vert \leq \left\Vert
u\right\Vert _{\infty }\equiv \sup_{y\in \mathbb{R}^{N}}\left\vert
u(y)\right\vert $\ for all such $u$. The mapping $M$ is called the \textit{%
mean value on} $\mathbb{R}^{N}$ \textit{for} $\mathcal{H}$. It is also a
fact, as the characteristic function of all relatively compact set in $%
\mathbb{R}^{N}$ lies in $L^{1}(\mathbb{R}^{N})$, that 
\begin{equation}
M(u)=\lim_{r\rightarrow +\infty }\frac{1}{\left\vert B_{r}\right\vert }%
\int_{B_{r}}u(y)dy\;\;\;\;\;\;\;\;\;\;\;\;\;\;  \label{2.3}
\end{equation}%
where $B_{r}$ stands for the bounded open ball in $\mathbb{R}^{N}$ with
radius $r$, and $\left\vert B_{r}\right\vert $ denotes its Lebesgue measure.
Expression (\ref{2.3}) also works for $u\in L_{\text{loc}}^{1}(\mathbb{R}%
^{N})$ provided that the above limit makes sense. In connection with the
dynamical systems, we have the following Birkhoff ergodic theorem (see \cite%
{16}).

\begin{theorem}[\textit{Birkhoff ergodic theorem}]
\label{t2.0}Let $T$ be a dynamical system acting on the probability space $%
(\Omega ,\mathcal{M},\mu )$. Let $f\in L^{p}(\Omega )$, $p\geq 1$. Then for
almost all $\omega \in \Omega $ the realization $x\mapsto f(T(x)\omega )$
possesses a mean value in the sense of \emph{(\ref{2.3})}. Furthermore, the
mean value $M(f(T(\cdot )\omega ))$ is invariant and 
\begin{equation*}
\int_{\Omega }f(\omega )d\mu =\int_{\Omega }M(f(T(\cdot )\omega ))d\mu .
\end{equation*}%
Moreover if the dynamical system $T$ is ergodic, then 
\begin{equation*}
M(f(T(\cdot )\omega ))=\int_{\Omega }fd\mu \text{ \ for }\mu \text{-a.e. }%
\omega \in \Omega .
\end{equation*}
\end{theorem}

\begin{definition}
\label{d2.1}\emph{By a }homogenization supralgebra\emph{\ (or }$H$\emph{%
-supralgebra, in short) on }$\mathbb{R}^{N}$\emph{\ for }$\mathcal{H}$\emph{%
\ we mean any closed subalgebra of }$\mathcal{B}(\mathbb{R}^{N})$\emph{\
which contains the constants, is closed under complex conjugation and whose
elements possess a mean value for }$\mathcal{H}$\emph{.}
\end{definition}

\begin{remark}
\label{r2.1}\emph{From the above definition we see that the concept of }$H$%
\emph{-supralgebra is more general than those of }$H$\emph{-algebra \cite{29}
and of algebra with mean value \cite{Zhikov4}. In fact any separable }$H$%
\emph{-supralgebra is an }$H$\emph{-algebra while any algebra with mean
value is an }$H$\emph{-supralgebra as any uniformly continuous function is
continuous.}
\end{remark}

Let $A$ be an $H$-supralgebra on $\mathbb{R}^{N}$ (for $\mathcal{H}$). It is
known that $A$ (endowed with the sup norm topology) is a commutative $%
\mathcal{C}$*-algebra with identity. We denote by $\Delta (A)$ the spectrum
of $A$ and by $\mathcal{G}$ the Gelfand transformation on $A$. We recall
that $\Delta (A)$ (a subset of the topological dual $A^{\prime }$ of $A$) is
the set of all nonzero multiplicative linear functionals on $A$, and $%
\mathcal{G}$ is the mapping from $A$ to $\mathcal{C}(\Delta (A))$ such that $%
\mathcal{G}(u)(s)=\left\langle s,u\right\rangle $ ($s\in \Delta (A)$), where 
$\left\langle ,\right\rangle $ denotes the duality pairing between $%
A^{\prime }$ and $A$. We endow $\Delta (A)$ with the relative weak$\ast $
topology on $A^{\prime }$. Then using the well-known theorem of Stone (see
e.g., either \cite{21} or more precisely \cite[Theorem IV.6.18, p. 274]{16})
one can easily show that the spectrum $\Delta (A)$ is a compact topological
space, and the Gelfand transformation $\mathcal{G}$ is an isometric
isomorphism identifying $A$ with $\mathcal{C}(\Delta (A))$ (the continuous
functions on $\Delta (A)$) as $\mathcal{C}$*-algebras. Next, since each
element of $A$ possesses a mean value, this yields a map $u\mapsto M(u)$
(denoted by $M$ and called the mean value) which is a nonnegative continuous
linear functional on $A$ with $M(1)=1$, and so provides us with a linear
nonnegative functional $\psi \mapsto M_{1}(\psi )=M(\mathcal{G}^{-1}(\psi ))$
defined on $\mathcal{C}(\Delta (A))=\mathcal{G}(A)$, which is clearly
bounded. Therefore, by the Riesz-Markov theorem, $M_{1}(\psi )$ is
representable by integration with respect to some Radon measure $\beta $ (of
total mass $1$) in $\Delta (A)$, called the $M$\textit{-measure} for $A$ %
\cite{26}. It is evident that we have 
\begin{equation}
M(u)=\int_{\Delta (A)}\mathcal{G}(u)d\beta \text{\ for }u\in A\text{.}
\label{3.5'}
\end{equation}

The spectrum of a Banach algebra is a very abstract concept. We give in the
following result a characterization of the spectrum of some particular
Banach algebras.

\begin{proposition}
\label{p2.0}Let $A$ be an $H$-supralgebra. Assume $A$ separates the points
of $\mathbb{R}^{N}$. Then $\Delta (A)$ is the Stone-\v{C}ech
compactification of $\mathbb{R}^{N}$.
\end{proposition}

\begin{proof}
For each $y\in \mathbb{R}^{N}$ define an element $\phi _{y}\in \Delta (A)$
by $\phi _{y}(u)=u(y)$, $u\in A$. Then the mapping $\phi :y\mapsto \phi _{y}$
from $\mathbb{R}^{N}$ into $\Delta (A)$ is continuous and has dense range.
In fact since the topology in $\Delta (A)$ is the weak$\ast $ one and
further the mappings $y\mapsto \phi _{y}(u)=u(y)$, $u\in A$, are continuous
on $\mathbb{R}^{N}$, it follows that $\phi $ is continuous. Now supposing
that $\phi (\mathbb{R}^{N})$ is not dense in $\Delta (A)$ we derive the
existence of a nonempty open subset $U$ of $\Delta (A)$ such that $U\cap
\phi (\mathbb{R}^{N})=\emptyset $. Then by Urysohn's lemma there exists $%
v\in \mathcal{C}(\Delta (A))$ with $v\neq 0$ and $\left. v\right\vert
_{\Delta (A))\backslash U}=0$. By the Gelfand representation theorem, $v=%
\mathcal{G}(u)$ for some $u\in A$. But then 
\begin{equation*}
u(y)=\phi _{y}(u)=\mathcal{G}(u)(\phi _{y})=v(\phi _{y})=0
\end{equation*}%
for all $y\in \mathbb{R}^{N}$, contradicting $u\neq 0$. Thus $\phi (\mathbb{R%
}^{N})$ is dense in $\Delta (A)$.

Next, every $f$ in $A$ (viewed as element of $\mathcal{B}(\mathbb{R}^{N})$)
extends continuously to $\Delta (A)$ in the sense that there exists $%
\widehat{f}\in \mathcal{C}(\Delta (A))$ such that $\widehat{f}(\phi
(y))=f(y) $ for all $y\in \mathbb{R}^{N}$ (just take $\widehat{f}=\mathcal{G}%
(f)$). Finally assume that $A$ separates the points of $\mathbb{R}^{N}$.
Then the mapping $\phi :\mathbb{R}^{N}\rightarrow \phi (\mathbb{R}^{N})$ is
a homeomorphism. In fact, we only need to prove that $\phi $ is injective.
For that, let $y,z\in \mathbb{R}^{N}$ with $y\neq z$; since $A$ separates
the points of $\mathbb{R}^{N}$, there exists a function $u\in A$ such that $%
u(y)\neq u(z)$, hence $\phi _{y}\neq \phi _{z}$, and our claim is justified.
We therefore conclude that the couple $(\Delta (A),\phi )$ is the Stone-\v{C}%
ech compactification of $\mathbb{R}^{N}$.
\end{proof}

The following result is classically known.

\begin{proposition}
\label{p2.5}\emph{(1)} Assume $A=\mathcal{C}_{\text{\emph{per}}}(Y)$\ is the
algebra of $Y$-periodic continuous functions on $\mathbb{R}_{y}^{N}$\ ($Y=(-%
\frac{1}{2},\frac{1}{2})^{N}$). Then its spectrum is the $N$-dimensional
torus $\mathbb{T}^{N}=\mathbb{R}^{N}/\mathbb{Z}^{N}$.\emph{\ (2)} Assume $%
A=AP(\mathbb{R}_{y}^{N})$\ is the algebra of all almost periodic continuous
functions on $\mathbb{R}_{y}^{N}$ defined as the vector space consisting of
all functions defined on $\mathbb{R}_{y}^{N}$ that are uniformly
approximated by finite linear combinations of the functions in the set $%
\{\exp (2i\pi k\cdot y):k\in \mathbb{R}^{N}\}$. Then its spectrum $\Delta
(AP(\mathbb{R}_{y}^{N}))$\ is a compact topological group homeomorphic to
the Bohr compactification of $\mathbb{R}^{N}$.
\end{proposition}

Next, the partial derivative of index $i$ ($1\leq i\leq N$) on $\Delta (A)$
is defined to be the mapping $\partial _{i}=\mathcal{G}\circ \partial
/\partial y_{i}\circ \mathcal{G}^{-1}$ (usual composition) of $\mathcal{D}%
^{1}(\Delta (A))=\{\varphi \in \mathcal{C}(\Delta (A)):\mathcal{G}%
^{-1}(\varphi )\in A^{1}\}$ into $\mathcal{C}(\Delta (A)),$ where $%
A^{1}=\{\psi \in \mathcal{C}^{1}(\mathbb{R}^{N}):$ $\psi ,\partial \psi
/\partial y_{i}\in A$ ($1\leq i\leq N$)$\}$. Higher order derivatives are
defined analogously, and one also defines the space $A^{m}$ (integers $m\geq
1$) to be the space of all $\psi \in \mathcal{C}^{m}(\mathbb{R}_{y}^{N})$
such that $D_{y}^{\alpha }\psi =\frac{\partial ^{\left\vert \alpha
\right\vert }\psi }{\partial y_{1}^{\alpha _{1}}\cdot \cdot \cdot \partial
y_{N}^{\alpha _{N}}}\in A$ for every $\alpha =(\alpha _{1},...,\alpha
_{N})\in \mathbb{N}^{N}$ with $\left\vert \alpha \right\vert \leq m$, and we
set $A^{\infty }=\cap _{m\geq 1}A^{m}$. At the present time, let $\mathcal{D}%
(\Delta (A))=\{\varphi \in \mathcal{C}(\Delta (A)):$ $\mathcal{G}%
^{-1}(\varphi )\in A^{\infty }\}.$ Endowed with a suitable locally convex
topology, $A^{\infty }$ (resp. $\mathcal{D}(\Delta (A))$) is a Fr\'{e}chet
space and further, $\mathcal{G}$ viewed as defined on $A^{\infty }$ is a
topological isomorphism of $A^{\infty }$ onto $\mathcal{D}(\Delta (A))$. It
is worth recalling that $A^{\infty }$ is the deterministic analog of the
space $\mathcal{C}^{\infty }(\Omega )$ defined in Subsection 2.1.

Analogously to the space $\mathcal{D}^{\prime }(\mathbb{R}^{N})$, we now
define the space of distributions on $\Delta (A)$ to be the space of all
continuous linear form on $\mathcal{D}(\Delta (A))$. We denote it by $%
\mathcal{D}^{\prime }(\Delta (A))$ and we endow it with the strong dual
topology. It is an easy exercise to see that if $A^{\infty }$ is dense in $A$
(this is the case when, e.g., $A$ is translation invariant and moreover each
element of $A$ is uniformly continuous; see \cite[Proposition 2.3]{40} for
the justification. We will also see at the end of this subsection that this
density result is a fact when dealing with such kind of $H$-supralgebras
since one may connect their spectrums to a dynamical system, and then
recover the said density result by just using the results of Subsection 2.1)
then $L^{p}(\Delta (A))$ ($1\leq p\leq \infty $) is a subspace of $\mathcal{D%
}^{\prime }(\Delta (A))$ (with continuous embedding), so that one may define
the Sobolev spaces on $\Delta (A)$ as follows. 
\begin{equation*}
W^{1,p}(\Delta (A))=\{u\in L^{p}(\Delta (A)):\text{ }\partial _{i}u\in
L^{p}(\Delta (A))\text{ (}1\leq i\leq N\text{)}\}\;(1\leq p<\infty )
\end{equation*}%
where the derivative $\partial _{i}u$ is taken in the distribution sense on $%
\Delta (A)$. We equip $W^{1,p}(\Delta (A))$ with the norm

\begin{equation*}
\begin{array}{l}
||u||_{W^{1,p}(\Delta (A))}=\left[ ||u||_{L^{p}(\Delta
(A))}^{p}+\sum_{i=1}^{N}||\partial _{i}u||_{L^{p}(\Delta (A))}^{p}\right] ^{%
\frac{1}{p}}\text{ \thinspace }\left( u\in W^{1,p}(\Delta (A))\right) , \\ 
1\leq p<\infty ,%
\end{array}%
\end{equation*}%
which makes it a Banach space. To that space are attached some other spaces
such as $W^{1,p}(\Delta (A))/\mathbb{C}=\{u\in W^{1,p}(\Delta
(A)):\int_{\Delta (A)}ud\beta =0\}$ and its separated completion $%
W_{\#}^{1,p}(\Delta (A))$; we refer to \cite{26} for a documented
presentation of these spaces.

As we have said a while ago, we end this subsection with an important result
connecting the dynamical systems to the spectrum of some $H$-supralgebras.

\begin{theorem}
\label{t2.2}Let $A$ be an $H$-supralgebra on $\mathbb{R}^{N}$. Suppose $A$
is translation invariant and that each of its elements is uniformly
continuous (thus $A$ is an algebra with mean value). Then the translations $%
T(y):\mathbb{R}^{N}\rightarrow \mathbb{R}^{N}$, $T(y)x=x+y$, extend to a
group of homeomorphisms $T(y):\Delta (A)\rightarrow \Delta (A)$, $y\in 
\mathbb{R}^{N}$, which forms a continuous $N$-dimensional dynamical system
on $\Delta (A)$ whose invariant probability measure is precisely the $M$%
-measure $\beta $ for $A$.
\end{theorem}

\begin{proof}
As $A$ is translation invariant, each translation $T(y)$ induces an
isometric isomorphism still denoted by $T(y)$, from $A$ onto $A$, defined by 
$T(y)u=u(\cdot +y)$ for $u\in A$. So define $\widetilde{T}(y):\mathcal{C}%
(\Delta (A))\rightarrow \mathcal{C}(\Delta (A))$ by 
\begin{equation*}
\widetilde{T}(y)\mathcal{G}(u)=\mathcal{G}(T(y)u)\;\;(u\in A)
\end{equation*}%
where $\mathcal{G}$ denotes the Gelfand transformation on $A$. Then $%
\widetilde{T}(y)$ is an isometric isomorphism of $\mathcal{C}(\Delta (A))$
onto itself; this is easily seen by the fact that $\mathcal{G}$ is an
isometric isomorphism of $A$ onto $\mathcal{C}(\Delta (A))$. Therefore, by
the classical Banach-Stone theorem there exists a unique homeomorphism $%
\overline{T}(y)$ of $\Delta (A)$ onto itself. The family thus constructed is
in fact a continuous $N$-dimensional dynamical system. Indeed the group
property easily comes from the equality $\mathcal{G}(T(y)u)(s)=\mathcal{G}%
(u)(\overline{T}(y)s)$ ($y\in \mathbb{R}^{N}$, $s\in \Delta (A)$, $u\in A$).
As far as the continuity property is concerned, let $(y_{n})_{n}$ be a
sequence in $\mathbb{R}^{N}$ and $(s_{d})_{d}$ be a net in $\Delta (A)$ such
that $y_{n}\rightarrow y$ in $\mathbb{R}^{N}$ and $s_{d}\rightarrow s$ in $%
\Delta (A)$. Then the uniform continuity of $u\in A$ leads to $%
T(y_{n})u\rightarrow T(y)u$ in $\mathcal{B}(\mathbb{R}^{N})$, and the
continuity of $\mathcal{G}$ gives $\mathcal{G}(T(y_{n})u)\rightarrow 
\mathcal{G}(T(y)u)$, the last convergence result being uniform in $\mathcal{C%
}(\Delta (A))$. Hence $\mathcal{G}(T(y_{n})u)(s_{d})\rightarrow \mathcal{G}%
(T(y)u)(s)$, which is equivalent to $\mathcal{G}(u)(\overline{T}%
(y_{n})s_{d})\rightarrow \mathcal{G}(u)(\overline{T}(y)s)$. As $\mathcal{C}%
(\Delta (A))$ separates the points of $\Delta (A)$, this yields $\overline{T}%
(y_{n})s_{d}\rightarrow \overline{T}(y)s$ in $\Delta (A)$, which implies
that the mapping $(y,s)\mapsto \overline{T}(y)s$, from $\mathbb{R}^{N}\times
\Delta (A)$ to $\Delta (A)$, is continuous. It remains to check that $\beta $
is the invariant measure for $\overline{T}$. But this easily comes from the
invariance under translations' property of the mean value and of the
integral representation (\ref{3.5'}). We keep using the notation $T(y)$ for $%
\overline{T}(y)$, and the proof is complete.
\end{proof}

With the above result, one may directly consider deterministic
homogenization theory in algebras with mean value as a particular case of
stochastic homogenization theory. That is why in the sequel, our results in
these particular $H$-supralgebras could be viewed as particular ones of
reiterated stochastic homogenization theory. However they are no less
important because so far, although widely used, the results stated in
Section 3 have never been proven before.

\subsection{The generalized Besicovitch spaces}

We can define the generalized Besicovitch spaces associated to a $H$%
-supralgebra. The notations are those of the preceding subsection. Let $A$
be a $H$-supralgebra on $\mathbb{R}^{N}$. Let $1\leq p<\infty $. If $u\in A$%
, then $\left| u\right| ^{p}\in A$ with $\mathcal{G}(\left| u\right|
^{p})=\left| \mathcal{G}(u)\right| ^{p}$. Hence the limit $%
\lim_{r\rightarrow +\infty }\frac{1}{\left| B_{r}\right| }\int_{B_{r}}\left|
u(y)\right| ^{p}dy$ exists and we have 
\begin{equation*}
\lim_{r\rightarrow +\infty }\frac{1}{\left| B_{r}\right| }\int_{B_{r}}\left|
u(y)\right| ^{p}dy=M(\left| u\right| ^{p})=\int_{\Delta (A)}\left| \mathcal{G%
}(u)\right| ^{p}d\beta .
\end{equation*}%
Hence, for $u\in A$, put 
\begin{equation*}
\left\| u\right\| _{p}=\left( M(\left| u\right| ^{p})\right)
^{1/p}.\;\;\;\;\;\;\;\;\;\;\;\;\;\;\;\;\;\;\;\;
\end{equation*}%
This defines a seminorm on $A$ with which $A$ is in general not separated
and not complete. First we denote by $B_{A}^{p}$ the closure of $A$ with
respect to $\left\| \cdot \right\| _{p}$. Then It is known that $B_{A}^{p}$
is a complete seminormed vector space verifying $B_{A}^{q}\subset B_{A}^{p}$
for $1\leq p\leq q<\infty $. From this last property one may naturally
define the space $B_{A}^{\infty }$ as follows: 
\begin{equation*}
B_{A}^{\infty }=\{f\in \cap _{1\leq p<\infty }B_{A}^{p}:\sup_{1\leq p<\infty
}\left\| f\right\| _{p}<\infty \}\text{.}\;\;\;\;\;\;\;\;\;
\end{equation*}%
We endow $B_{A}^{\infty }$ with the seminorm $\left[ f\right] _{\infty
}=\sup_{1\leq p<\infty }\left\| f\right\| _{p}$, which makes it a complete
seminormed space. We recall that the spaces $B_{A}^{p}$ ($1\leq p\leq \infty 
$) are not in general Fr\'{e}chet spaces since they are not separated in
general. The following properties are worth noticing \cite{CMP, SNW}:

\begin{itemize}
\item[(\textbf{1)}] The Gelfand transformation $\mathcal{G}:A\rightarrow 
\mathcal{C}(\Delta (A))$ extends by continuity to a unique continuous linear
mapping, still denoted by $\mathcal{G}$, of $B_{A}^{p}$ into $L^{p}(\Delta
(A))$. Furthermore if $u\in B_{A}^{p}\cap L^{\infty }(\mathbb{R}_{y}^{N})$
then $\mathcal{G}(u)\in L^{\infty }(\Delta (A))$ and $\left\Vert \mathcal{G}%
(u)\right\Vert _{L^{\infty }(\Delta (A))}\leq \left\Vert u\right\Vert
_{L^{\infty }(\mathbb{R}_{y}^{N})}$.

\item[(\textbf{2)}] The mean value $M$ viewed as defined on $A$, extends by
continuity to a positive continuous linear form (still denoted by $M$) on $%
B_{A}^{p}$ satisfying $M(u)=\int_{\Delta (A)}\mathcal{G}(u)d\beta $ ($u\in
B_{A}^{p}$). Furthermore, $M(\tau _{a}u)=M(u)$ for each $u\in B_{A}^{p}$ and
all $a\in \mathbb{R}^{N}$, where $\tau _{a}u(y)=u(y-a)$ for almost all $y\in 
\mathbb{R}^{N}$.

\item[(\textbf{3)}] Let $1\leq p,q,r<\infty $ be such that $\frac{1}{p}+%
\frac{1}{q}=\frac{1}{r}\leq 1$. The usual multiplication $A\times
A\rightarrow A$; $(u,v)\mapsto uv$, extends by continuity to a bilinear form 
$B_{A}^{p}\times B_{A}^{q}\rightarrow B_{A}^{r}$ with 
\begin{equation*}
\left\| uv\right\| _{r}\leq \left\| u\right\| _{p}\left\| v\right\| _{q}%
\text{\ for }(u,v)\in B_{A}^{p}\times B_{A}^{q}.
\end{equation*}
\end{itemize}

The following result will be of great interest in the work.

\begin{proposition}
\label{p2.3}Let $A$ be a $H$-supralgebra on $\mathbb{R}^{N}$. Assume each
element of $A$ is uniformly continuous and moreover $A$ is translation
invariant (i.e. $\tau _{a}u=u(\cdot +a)\in A$ for all $u\in A$ and all $a\in 
\mathbb{R}^{N}$). Then $A^{\infty }$ is dense in $B_{A}^{p}$.
\end{proposition}

\begin{proof}
Since $A$ is an algebra with mean value, the result follows from %
\cite[Proposition 2.4]{SNW}.
\end{proof}

Now, let $u\in B_{A}^{p}$ ($1\leq p<\infty $); then $\left| u\right| ^{p}\in
B_{A}^{1}$ (this is easily seen) and so, by part (\textbf{2)} above one has $%
M(\left| u\right| ^{p})=\int_{\Delta (A)}\left| \mathcal{G}(u)\right|
^{p}d\beta =\left\| \mathcal{G}(u)\right\| _{L^{p}(\Delta (A))}^{p}$. Thus
for $u\in B_{A}^{p}$ we have $\left\| u\right\| _{p}=\left( M(\left|
u\right| ^{p})\right) ^{1/p}$, and $\left\| u\right\| _{p}=0$ if and only if 
$\mathcal{G}(u)=0$. Unfortunately, the mapping $\mathcal{G}$ (defined on $%
B_{A}^{p}$) is not in general injective. So let $\mathcal{N}=Ker\mathcal{G}$
(the kernel of $\mathcal{G}$) and let 
\begin{equation*}
\mathcal{B}_{A}^{p}=B_{A}^{p}/\mathcal{N}.\;\;\;\;\;\;\;\;\;
\end{equation*}%
Endowed with the norm 
\begin{equation*}
\left\| u+\mathcal{N}\right\| _{\mathcal{B}_{A}^{p}}=\left\| u\right\|
_{p}\;\;(u\in B_{A}^{p}),\;\;
\end{equation*}%
$\mathcal{B}_{A}^{p}$ is a Banach space with the following property.

\begin{theorem}[\protect\cite{CMP}]
\label{t2.3}The mapping $\mathcal{G}:B_{A}^{p}\rightarrow L^{p}(\Delta (A))$
induces an isometric isomorphism $\mathcal{G}_{1}$ of $\mathcal{B}_{A}^{p}$
onto $L^{p}(\Delta (A))$.
\end{theorem}

As a first consequence of the preceding theorem one can define the mean
value of $u+\mathcal{N}$ (for each $u\in B_{A}^{p}$) as follows: 
\begin{equation}
M_{1}(u+\mathcal{N})=M(u)\text{, so that }M_{1}(u+\mathcal{N}%
)=\lim_{r\rightarrow +\infty }\frac{1}{\left| B_{r}\right| }%
\int_{B_{r}}u(y)dy.  \label{0.1}
\end{equation}

One crucial result that can be derived from the preceding theorem is the
following

\begin{corollary}
\label{c2.1}The following hold true:

\begin{itemize}
\item[(i)] The spaces $\mathcal{B}_{A}^{p}$ are reflexive for $1<p<\infty $;

\item[(ii)] The topological dual of the space $\mathcal{B}_{A}^{p}$ $(1\leq
p<\infty )$ is the space $\mathcal{B}_{A}^{p^{\prime }}$ $(p^{\prime
}=p/(p-1))$, the duality being given by 
\begin{equation*}
\begin{array}{l}
\left\langle u+\mathcal{N},v+\mathcal{N}\right\rangle _{\mathcal{B}%
_{A}^{p^{\prime }},\mathcal{B}_{A}^{p}}=M(uv)=\int_{\Delta (A)}\mathcal{G}%
_{1}(u+\mathcal{N})\mathcal{G}_{1}(v+\mathcal{N})d\beta \\ 
\text{for }u\in B_{A}^{p^{\prime }}\text{ and }v\in B_{A}^{p}\text{.}%
\end{array}%
\end{equation*}
\end{itemize}
\end{corollary}

This result is easily proven by using the properties of the $L^{p}$-spaces
and the above isometric isomorphism.

\begin{remark}
\label{r2.1'}\emph{The space }$\mathcal{B}_{A}^{p}$\emph{\ is the separated
completion of }$B_{A}^{p}$\emph{\ and the canonical mapping of }$B_{A}^{p}$%
\emph{\ into }$\mathcal{B}_{A}^{p}$\emph{\ is just the canonical surjection
of }$B_{A}^{p}$\emph{\ onto }$\mathcal{B}_{A}^{p}$\emph{; see once more %
\cite[Chap. II, Sect. 3, no 7]{9} for the theory of completion.}
\end{remark}

Another definition which will be of great interest in the forthcoming
sections is

\begin{definition}
\label{d2.2}\emph{\ An $H$-supralgebra }$A$ \emph{on }$\mathbb{R}^{N}$\emph{%
\ is }ergodic\emph{\ if for every }$u\in B_{A}^{1}$\emph{\ such that }$%
\left\| u-u(\cdot +a)\right\| _{1}=0$\emph{\ for every }$a\in \mathbb{R}^{N}$%
\emph{\ we have }$\left\| u-M(u)\right\| _{1}=0$\emph{.}
\end{definition}

The above definition is equivalent to say that any $B_{A}^{1}$-translation
invariant function is $B_{A}^{1}$-constant, that is, the dynamical system $T$
defined on $\mathbb{R}^{N}$\ by $T(y)x=x+y$ is ergodic in the sense of
Subsection 2.1.

An equivalent property stated by Casado Diaz and Gayte is given in the
following proposition.

\begin{proposition}[\protect\cite{12}]
\label{p2.4}An $H$-supralgebra $A$\ on $\mathbb{R}^{N}$\ is ergodic\ if and
only if 
\begin{equation}
\lim_{r\rightarrow +\infty }\left\| \frac{1}{\left| B_{r}\right| }%
\int_{B_{r}}u(\cdot +y)dy-M(u)\right\| _{p}=0\text{\ for all }u\in B_{A}^{p}%
\text{,\ }1\leq p<\infty \text{.}  \label{Eq.1}
\end{equation}
\end{proposition}

The following result provides us with a few examples of ergodic $H$%
-supralgebras (see next section for its application).

\begin{lemma}[\protect\cite{CMP}]
\label{l2.0}Let $A$ be an $H$-supralgebra on $\mathbb{R}^{N}$ with the
following property: For any $u\in A$, 
\begin{equation}
\lim_{r\rightarrow +\infty }\frac{1}{\left| B_{r}\right| }%
\int_{B_{r}}u(x+y)dx=M(u)\text{\ uniformly with respect to }y.  \label{Eq.2}
\end{equation}%
Then $A$ is ergodic.
\end{lemma}

In order to simplify the presentation of the paper we will from now on, use
the same letter $u$ (if there is no danger of confusion) to denote the
equivalence class of an element $u\in B_{A}^{p}$. The symbol $\varrho $ will
denote the canonical mapping of $B_{A}^{p}$ onto $\mathcal{B}%
_{A}^{p}=B_{A}^{p}/\mathcal{N}$.

Our goal here is to define another space attached to $\mathcal{B}_{A}^{p}$.
Let $u\in \mathcal{D}^{\prime }(\Delta (A))$, and let $\alpha \in \mathbb{N}%
^{N}$. We know that $\partial ^{\alpha }u\in \mathcal{D}^{\prime }(\Delta
(A))$ exists and is defined by 
\begin{equation}
\left\langle \partial ^{\alpha }u,\varphi \right\rangle =(-1)^{\left| \alpha
\right| }\left\langle u,\partial ^{\alpha }\varphi \right\rangle \text{\ for
any }\varphi \in \mathcal{D}(\Delta (A))\text{.}  \label{2.7'}
\end{equation}%
This leads to the following definition.

\begin{definition}
\label{d2.3}\emph{The }formal derivative of order $\alpha \in \mathbb{N}^{N}$
\emph{is defined to be the formal operator on }$\mathcal{B}_{A}^{p}$\emph{\
given by }%
\begin{equation}
\overline{D}_{y}^{\alpha }=\mathcal{G}_{1}^{-1}\circ \partial ^{\alpha
}\circ \mathcal{G}_{1}  \label{2.2'}
\end{equation}%
\emph{where }$\partial ^{\alpha }$\emph{\ is defined above. In particular,
for }$\alpha =(\delta _{ij})_{1\leq j\leq N}$\emph{\ with }$1\leq i\leq N$%
\emph{, }$\overline{D}_{y}^{\alpha }$\emph{\ is denoted by }$\overline{%
\partial }/\partial y_{i}$\emph{\ and is called the }formal derivative of
index\emph{\ }$i$\emph{.}
\end{definition}

\begin{remark}
\label{r2.2}\emph{For }$u\in B_{A}^{1,p}$\emph{\ (that is the space of }$%
u\in B_{A}^{p}$\emph{\ such that }$D_{y}u\in (B_{A}^{p})^{N}$\emph{) we have 
} 
\begin{equation*}
\mathcal{G}_{1}\left( \varrho \left( \frac{\partial u}{\partial y_{i}}%
\right) \right) =\mathcal{G}\left( \frac{\partial u}{\partial y_{i}}\right)
=\partial _{i}\mathcal{G}\left( u\right) =\partial _{i}\mathcal{G}_{1}\left(
\varrho (u)\right) =(\text{\emph{by definition}})\,\mathcal{G}_{1}\left( 
\frac{\overline{\partial }}{\partial y_{i}}(\varrho (u))\right) ,
\end{equation*}%
\emph{hence } 
\begin{equation*}
\varrho \left( \frac{\partial u}{\partial y_{i}}\right) =\frac{\overline{%
\partial }}{\partial y_{i}}(\varrho
(u)),\;\;\;\;\;\;\;\;\;\;\;\;\;\;\;\;\;\;\;\;\;\;
\end{equation*}%
\emph{or equivalently, } 
\begin{equation}
\varrho \circ \frac{\partial }{\partial y_{i}}=\frac{\overline{\partial }}{%
\partial y_{i}}\circ \varrho \text{\ \emph{on }}B_{A}^{1,p}.  \label{2.3'}
\end{equation}
\end{remark}

\bigskip We return for a while to the framework of the preceding subsection
and assume that the hypotheses of Theorem \ref{t2.2}\ are satisfied. Let $%
\{T(y):y\in \mathbb{R}^{N}\}$ be the dynamical system constructed in Theorem %
\ref{t2.2}. We know by the results of Subsection 2.1 that $T(y)$ induces a $%
N $-parameter group of isometries $U(y):L^{p}(\Delta (A))\rightarrow
L^{p}(\Delta (A))$. By the properties of $\mathcal{G}_{1}$, this also yields
a $N$-parameter group of isometries $\mathcal{G}_{1}^{-1}\circ U(y)\circ 
\mathcal{G}_{1}:\mathcal{B}_{A}^{p}\rightarrow \mathcal{B}_{A}^{p}$. We
denote by $D_{i,p}$ the generators of $\mathcal{G}_{1}^{-1}\circ U(y)\circ 
\mathcal{G}_{1}$. Now, let $u\in A^{1}$; we have $\partial _{i}\mathcal{G}%
(u)=\mathcal{G}(\frac{\partial u}{\partial y_{i}})=\mathcal{G}_{1}(\varrho (%
\frac{\partial u}{\partial y_{i}}))$, so that $\varrho (\frac{\partial u}{%
\partial y_{i}})=\frac{\overline{\partial }}{\partial y_{i}}(\varrho (u))$
by the preceding remark. But since $\frac{\partial u}{\partial y_{i}}$ is
the derivative along the direction $e_{i}=(\delta _{ij})_{1\leq i\leq N}$ of
the dynamical system induced by the translations in $\mathbb{R}^{N}$, it is
immediate that 
\begin{equation*}
D_{i,p}(\varrho (u))=\frac{\overline{\partial }}{\partial y_{i}}(\varrho
(u)),
\end{equation*}%
so that 
\begin{equation}
D_{i,p}=\frac{\overline{\partial }}{\partial y_{i}}.\;\;\;\;\;\;\;\;\;
\label{2.12}
\end{equation}%
The above equality is crucial in the process of viewing homogenization in
algebras with mean value as a special case of stochastic homogenization.
Indeed in the case when $\Omega =\Delta (A)$, it allows to just replace $%
\mathcal{C}^{\infty }(\Omega )$ by the space $\mathcal{G}_{1}(\varrho
(A^{\infty }))=\mathcal{G}(A^{\infty })=\mathcal{D}(\Delta (A))$ which plays
exactly the same role since firstly, it is dense in $L^{p}(\Delta (A))$ for
all $1\leq p<\infty $ and secondly, for all $u\in \mathcal{D}(\Delta
(A))\equiv \mathcal{C}^{\infty }(\Delta (A))$ we have $u\in L^{\infty
}(\Delta (A))$ and $\partial ^{\alpha }u\in L^{\infty }(\Delta (A))$ for all 
$\alpha \in \mathbb{N}^{N}$. This remark will be particularly used in
Section 6 when dealing with the homogenization of some Stokes' type
equations.

Now, set (for $1\leq p<\infty $) 
\begin{equation*}
\mathcal{B}_{A}^{1,p}=\left\{ u\in \mathcal{B}_{A}^{p}:\frac{\overline{%
\partial }u}{\partial y_{i}}\in \mathcal{B}_{A}^{p},\;\text{for }1\leq i\leq
N\right\} .\;\;\;\;\;\;\;\;\;\;
\end{equation*}%
We endow $\mathcal{B}_{A}^{1,p}$ with the norm 
\begin{equation*}
\left\Vert u\right\Vert _{\mathcal{B}_{A}^{1,p}}=\left[ \left\Vert
u\right\Vert _{p}^{p}+\sum_{i=1}^{N}\left\Vert \frac{\overline{\partial }u}{%
\partial y_{i}}\right\Vert _{p}^{p}\right] ^{1/p}\text{\ \ }(u\in \mathcal{B}%
_{A}^{1,p})
\end{equation*}%
which makes it a Banach space with the property that the restriction of $%
\mathcal{G}_{1}$ to $\mathcal{B}_{A}^{1,p}$ is an isometric isomorphism from 
$\mathcal{B}_{A}^{1,p}$ onto $W^{1,p}(\Delta (A))$. However we will be
mostly concerned with the subspace $\mathcal{B}_{A}^{1,p}/\mathbb{C}$ of $%
\mathcal{B}_{A}^{1,p}$ consisting of functions $u\in \mathcal{B}_{A}^{1,p}$
with $M_{1}(u)\equiv M(u)=0$. Equipped with the seminorm 
\begin{equation*}
\left\Vert u\right\Vert _{\mathcal{B}_{A}^{1,p}/\mathbb{C}}=\left\Vert 
\overline{D}_{y}u\right\Vert _{p}:=\left[ \sum_{i=1}^{N}\left\Vert \frac{%
\overline{\partial }u}{\partial y_{i}}\right\Vert _{p}^{p}\right]
^{1/p}\;\;(u\in \mathcal{B}_{A}^{1,p}/\mathbb{C})
\end{equation*}%
where $\overline{D}_{y}=(\overline{\partial }/\partial y_{i})_{1\leq i\leq
N} $, $\mathcal{B}_{A}^{1,p}/\mathbb{C}$ is a locally convex topological
space which is in general not separated and not complete. We denote by $%
\mathcal{B}_{\#A}^{1,p}$ the separated completion of $\mathcal{B}_{A}^{1,p}/%
\mathbb{C}$ with respect to $\left\Vert \cdot \right\Vert _{\mathcal{B}%
_{A}^{1,p}/\mathbb{C}}$, and by $J_{1}$ the canonical mapping of $\mathcal{B}%
_{A}^{1,p}/\mathbb{C}$ into $\mathcal{B}_{\#A}^{1,p}$. By the theory of
completion of the uniform spaces \cite[Chap. II, Sect. 3, no 7]{9} it is a
fact that the mapping $\overline{\partial }/\partial y_{i}:\mathcal{B}%
_{A}^{1,p}/\mathbb{C}\rightarrow \mathcal{B}_{A}^{p}$ extends by continuity
to a unique continuous linear mapping still denoted by $\overline{\partial }%
/\partial y_{i}:\mathcal{B}_{\#A}^{1,p}\rightarrow \mathcal{B}_{A}^{p}$ and
satisfying 
\begin{equation}
\frac{\overline{\partial }}{\partial y_{i}}\circ J_{1}=\frac{\overline{%
\partial }}{\partial y_{i}}\text{\ and }\left\Vert u\right\Vert _{\mathcal{B}%
_{\#A}^{1,p}}=\left\Vert \overline{D}_{y}u\right\Vert _{p}\;\;(u\in \mathcal{%
B}_{\#A}^{1,p})  \label{2.4'}
\end{equation}%
where $\overline{D}_{y}=(\overline{\partial }/\partial y_{i})_{1\leq i\leq
N} $. Since $\mathcal{G}_{1}$ is an isometric isomorphism of $\mathcal{B}%
_{A}^{1,p}$ onto $W^{1,p}(\Delta (A))$ we have by (\ref{2.2'}) that the
restriction of $\mathcal{G}_{1}$ to $\mathcal{B}_{A}^{1,p}/\mathbb{C}$ sends
isometrically and isomorphically $\mathcal{B}_{A}^{1,p}/\mathbb{C}$ onto $%
W^{1,p}(\Delta (A))/\mathbb{C}$. So by \cite[Chap. II, Sect. 3, no 7]{9}
there exists a unique isometric isomorphism $\overline{\mathcal{G}}_{1}:%
\mathcal{B}_{\#A}^{1,p}\rightarrow W_{\#}^{1,p}(\Delta (A))$ such that 
\begin{equation}
\overline{\mathcal{G}}_{1}\circ J_{1}=J\circ \mathcal{G}_{1}\;\;\;\;\;\;\;\;%
\;\;\;\;\;\;\;\;\;\;\;\;\;  \label{2.5'}
\end{equation}%
and 
\begin{equation}
\partial _{i}\circ \overline{\mathcal{G}}_{1}=\mathcal{G}_{1}\circ \frac{%
\overline{\partial }}{\partial y_{i}}\;\;(1\leq i\leq N).\;\;\;\;\;\;\;\;\;
\label{2.6'}
\end{equation}%
We recall that $J$ is the canonical mapping of $W^{1,p}(\Delta (A))/\mathbb{C%
}$ into its separated completion $W_{\#}^{1,p}(\Delta (A))$ while $J_{1}$ is
the canonical mapping of $\mathcal{B}_{A}^{1,p}/\mathbb{C}$ into $\mathcal{B}%
_{\#A}^{1,p}$. Furthermore, as $J_{1}(\mathcal{B}_{A}^{1,p}/\mathbb{C})$ is
dense in $\mathcal{B}_{\#A}^{1,p}$ (this is classical), it follows that if $%
A^{\infty }$ is dense in $A$ (this is the case when $A$ is an algebra with
mean value), then $(J_{1}\circ \varrho )(A^{\infty }/\mathbb{C})$ is dense
in $\mathcal{B}_{\#A}^{1,p}$, where $A^{\infty }/\mathbb{C}=\{u\in A^{\infty
}:M(u)=0\}$.

\section{The stochastic $\Sigma $-convergence}

In this section we define the concept of stochastic $\Sigma $-convergence
which is the generalization of both two-scale convergence in the mean (of
Bourgeat et al. \cite{10}) and $\Sigma $-convergence (of Nguetseng \cite{26}%
). In all that follows, $Q$ is an open subset of $\mathbb{\mathbb{R}}^{N}$
and $A$ is an $H$-supralgebra on $\mathbb{\mathbb{R}}_{y}^{N}$. We use the
letter $\mathcal{G}$ to denote the Gelfand transformation on $A$. Points in $%
\Delta (A)$ are denoted by $s$. We still denote by $M$ the mean value on $%
\mathbb{\mathbb{R}}^{N}$ for the action $\mathcal{H}$ (see Section 2). The
compact space $\Delta (A)$ is equipped with the $M$-measure $\beta $ for $A$%
. Next, let $(\Omega ,\mathcal{M},\mu )$ denote a probability space and let $%
\{T(y):y\in \mathbb{\mathbb{R}}^{N}\}$ denote a $N$-dimensional dynamical
system acting on the probability space $(\Omega ,\mathcal{M},\mu )$. Points
in $\Omega $ are denoted by $\omega $. Finally, let $\varepsilon _{1}$ and $%
\varepsilon _{2}$ be two well separated functions of $\varepsilon $ tending
towards zero with $\varepsilon $, that is, $0<\varepsilon _{1},\varepsilon
_{2},\varepsilon _{2}/\varepsilon _{1}\rightarrow 0$ as $\varepsilon
\rightarrow 0$, and such that the functions $x\mapsto x/\varepsilon _{1}$
and $x\mapsto x/\varepsilon _{2}$ define two actions of $\mathbb{R}%
_{+}^{\ast }$ on $\mathbb{\mathbb{R}}^{N}$.

\begin{definition}
\label{d3.1}\emph{A bounded sequence }$(u_{\varepsilon })_{\varepsilon >0}$%
\emph{\ in }$L^{p}(Q\times \Omega )$\emph{\ (}$1\leq p<\infty $\emph{) is
said to }weakly stochastically $\Sigma $-converge\emph{\ in }$L^{p}(Q\times
\Omega )$\emph{\ to some }$u_{0}\in L^{p}(Q\times \Omega ;\mathcal{B}%
_{A}^{p})$\emph{\ if as }$\varepsilon \rightarrow 0$\emph{, we have } 
\begin{equation}
\int_{Q\times \Omega }u_{\varepsilon }(x,\omega )f\left( x,T\left( \frac{x}{%
\varepsilon _{1}}\right) \omega ,\frac{x}{\varepsilon _{2}}\right) dxd\mu
\rightarrow \iint_{Q\times \Omega \times \Delta (A)}\widehat{u}_{0}(x,\omega
,s)\widehat{f}(x,\omega ,s)dxd\mu d\beta \;\;\;\;\;\;\;\;  \label{3.1}
\end{equation}%
\emph{for every }$f\in \mathcal{C}_{0}^{\infty }(Q)\otimes \mathcal{C}%
^{\infty }(\Omega )\otimes A$\emph{, where }$\widehat{u}_{0}=\mathcal{G}%
_{1}\circ u_{0}$\emph{\ and }$\widehat{f}=\mathcal{G}\circ f=\mathcal{G}%
_{1}\circ (\varrho \circ f)$\emph{. We express this by writing} $%
u_{\varepsilon }\rightarrow u_{0}$\ \textit{stoch. in }$L^{p}(Q\times \Omega
)$\textit{-weak} $\Sigma $.
\end{definition}

We recall that $\mathcal{C}_{0}^{\infty }(Q)\otimes \mathcal{C}^{\infty
}(\Omega )\otimes A$ is the space of functions of the form 
\begin{equation*}
f(x,\omega ,y)=\sum_{\text{finite}}\varphi _{i}(x)\psi _{i}(\omega )g_{i}(y)%
\text{,\ \ }(x,\omega ,y)\in Q\times \Omega \times \mathbb{R}^{N}\text{,}
\end{equation*}%
with $\varphi _{i}\in \mathcal{C}_{0}^{\infty }(Q)$, $\psi _{i}\in \mathcal{C%
}^{\infty }(\Omega )$ and $g_{i}\in A$. Such functions are dense in $%
\mathcal{C}_{0}^{\infty }(Q)\otimes L^{p^{\prime }}(\Omega )\otimes A$ ($%
p^{\prime }=p/(p-1)$ for $1<p<\infty $, since $\mathcal{C}^{\infty }(\Omega
) $ is dense in $L^{p^{\prime }}(\Omega )$) and hence in $\mathcal{K}%
(Q;L^{p^{\prime }}(\Omega ))\otimes A$ ($\mathcal{K}(Q;L^{p^{\prime
}}(\Omega ))$ being the space of continuous functions of $Q$ into $%
L^{p^{\prime }}(\Omega )$ with compact support containing in $Q$; see e.g., %
\cite[ Proposition 5]{Bou1} for the denseness result). As $\mathcal{K}%
(Q;L^{p^{\prime }}(\Omega ))$ is dense in $L^{p^{\prime }}(Q;L^{p^{\prime
}}(\Omega ))=L^{p^{\prime }}(Q\times \Omega )$ and $L^{p^{\prime }}(Q\times
\Omega )\otimes A$ is dense in $L^{p^{\prime }}(Q\times \Omega ;A)$, the
uniqueness of the stochastic $\Sigma $-limit is ensured.

Before continuing our study, we need to make a comparison between the weak
stochastic $\Sigma $-convergence and other existing convergence methods
closed to it. For that, we must first state these convergence schemes:

\begin{itemize}
\item[(1)] A sequence $(u_{\varepsilon })_{\varepsilon >0}\subset L^{p}(Q)$%
\emph{\ }($1\leq p<\infty $) is said to weakly $\Sigma $-converge\ in $%
L^{p}(Q)$\ to some $v_{0}\in L^{p}(Q;\mathcal{B}_{A}^{p})$\ if as $E\ni
\varepsilon \rightarrow 0$, we\emph{\ }have\emph{\ } 
\begin{equation}
\int_{Q}u_{\varepsilon }(x)f\left( x,\frac{x}{\varepsilon _{2}}\right)
dx\rightarrow \iint_{Q\times \Delta (A)}\widehat{v}_{0}(x,s)\widehat{f}%
(x,s)dxd\beta  \label{3.1'}
\end{equation}%
for every $f\in L^{p^{\prime }}(Q;A)$\ ($1/p^{\prime }=1-1/p$), where\emph{\ 
}$\widehat{v}_{0}=\mathcal{G}_{1}\circ v_{0}$\ and\emph{\ }$\widehat{f}=%
\mathcal{G}_{1}\circ (\varrho \circ f)=\mathcal{G}\circ f$.

\item[(2)] A sequence $(u_{\varepsilon })_{\varepsilon >0}\subset
L^{p}(Q\times \Omega )$\emph{\ }($1\leq p<\infty $) is said to\emph{\ }%
stochastically two-scale converge\ in the mean\ to some\emph{\ }$v_{0}\in
L^{p}(Q\times \Omega )$\ if as\emph{\ }$\varepsilon \rightarrow 0$, we have 
\begin{equation}
\int_{Q\times \Omega }u_{\varepsilon }(x,\omega )f\left( x,T\left( \frac{x}{%
\varepsilon _{1}}\right) \omega \right) dxd\mu \rightarrow \iint_{Q\times
\Omega }v_{0}(x,\omega )f(x,\omega )dxd\mu  \label{3.2'}
\end{equation}%
for all admissible functions (in the sense of \cite[Section 3]{10}) $f\in
L^{p^{\prime }}\left( Q\times \Omega \right) $. We denote it by $%
u_{\varepsilon }\rightarrow u_{0}$\textit{\ stoch. in }$L^{p}\left( Q\times
\Omega \right) $\textit{-weak.}
\end{itemize}

\begin{remark}
\label{r3.0}\emph{The weak stochastic }$\Sigma $\emph{-convergence method
generalizes the above two convergence methods. Indeed, it is very important
to note that both the above definitions (\ref{3.1'}) and (\ref{3.2'}) imply
the boundedness of the sequence }$u_{\varepsilon }$\emph{\ either in }$%
L^{p}(Q)$\emph{\ or in }$L^{p}(Q\times \Omega )$\emph{, accordingly. With
this in mind, we see that if in (\ref{3.1}) we take }$f\in \mathcal{C}%
_{0}^{\infty }(Q)\otimes \mathcal{C}^{\infty }(\Omega )$\emph{, that is }$f$%
\emph{\ is constant with respect to }$y\in \mathbb{R}^{N}$\emph{, and next
using the density of the latter space in }$L^{p^{\prime }}(Q\times \Omega )$%
\emph{, then (\ref{3.1}) reads as (\ref{3.2'}) with }$v_{0}(x,\omega
)=\int_{\Delta (A)}\widehat{u}_{0}(x,\omega ,s)d\beta $\emph{\ by choosing
in }$L^{p^{\prime }}(Q\times \Omega )$\emph{\ admissible functions. If
besides we take in (\ref{3.1}) }$f\in \mathcal{C}_{0}^{\infty }(Q)\otimes A$%
\emph{, that is }$f$\emph{\ not depending upon the random variable }$\omega $%
\emph{\ and further if we choose }$u_{\varepsilon }$\emph{\ not depending on 
}$\omega $\emph{, then using the density of }$\mathcal{C}_{0}^{\infty
}(Q)\otimes A$\emph{\ in }$L^{p^{\prime }}(Q;A)$\emph{\ we readily get (\ref%
{3.1'}) with }$\widehat{v}_{0}(x,s)=\int_{\Omega }\widehat{u}_{0}(x,\omega
,s)d\mu $\emph{.}
\end{remark}

The following result is easily verified; its proof is left to the reader.

\begin{proposition}
\label{p2}Let $(u_{\varepsilon })_{\varepsilon >0}$ be a sequence in $%
L^{p}\left( Q\times \Omega \right) $. If $u_{\varepsilon }\rightarrow u_{0}$
stoch. in $L^{p}\left( Q\times \Omega \right) $-weak $\Sigma $, then $%
(u_{\varepsilon })_{\varepsilon >0}$ stochastically two-scale converges in
the mean towards $v_{0}(x,\omega )=\int_{\Delta (A)}\widehat{u}_{0}(x,\omega
,s)d\beta $ and 
\begin{equation*}
\int_{\Omega }u_{\varepsilon }\left( \cdot ,\omega \right) \psi (\omega
)d\mu \rightarrow \iint_{\Omega \times \Delta (A)}\widehat{u}_{0}\left(
\cdot ,\omega ,s\right) \psi (\omega )d\mu d\beta \text{ in }L^{1}(Q)\text{%
-weak }\forall \psi \in I_{nv}^{p^{\prime }}\left( \Omega \right) .
\end{equation*}
\end{proposition}

The next results provide us with a few examples of sequences that weakly
stochastically $\Sigma $-converge.

\begin{proposition}
\label{p3.3}Let $f\in \mathcal{K}(Q;\mathcal{C}^{\infty }(\Omega ;A))$.
Then, as $\varepsilon \rightarrow 0$, 
\begin{equation}
\int_{Q\times \Omega }f\left( x,T\left( \frac{x}{\varepsilon _{1}}\right)
\omega ,\frac{x}{\varepsilon _{2}}\right) dxd\mu \rightarrow \iint_{Q\times
\Omega \times \Delta (A)}\widehat{f}\left( x,\omega ,s\right) dxd\mu d\beta .
\label{Eqn1}
\end{equation}
\end{proposition}

\begin{proof}
Since $\mathcal{C}_{0}^{\infty }(Q)\otimes \mathcal{C}^{\infty }(\Omega
)\otimes A$ is dense in $\mathcal{K}(Q;\mathcal{C}^{\infty }(\Omega ;A))$ we
first check (\ref{Eqn1}) for $f$ in $\mathcal{C}_{0}^{\infty }(Q)\otimes 
\mathcal{C}^{\infty }(\Omega )\otimes A$. However, it is sufficient to do it
for $f$ under the form $f(x,\omega ,y)=\varphi (x)\psi (\omega )g(y)$ with $%
\varphi \in \mathcal{C}_{0}^{\infty }(Q)$, $\psi \in \mathcal{C}^{\infty
}(\Omega )$ and $g\in A$. But for such a $f$ we have 
\begin{eqnarray*}
\int_{Q\times \Omega }f\left( x,T\left( \frac{x}{\varepsilon _{1}}\right)
\omega ,\frac{x}{\varepsilon _{2}}\right) dxd\mu &=&\int_{Q}\left(
\int_{\Omega }\psi \left( T\left( \frac{x}{\varepsilon _{1}}\right) \omega
\right) d\mu \right) \varphi (x)g\left( \frac{x}{\varepsilon _{2}}\right) dx
\\
&=&\int_{Q}\left( \int_{\Omega }\psi (\omega )d\mu \right) \varphi
(x)g\left( \frac{x}{\varepsilon _{2}}\right) dx \\
&=&\left( \int_{\Omega }\psi (\omega )d\mu \right) \int_{Q}\varphi
(x)g\left( \frac{x}{\varepsilon _{2}}\right) dx
\end{eqnarray*}%
where the second equality above is due to the Fubini's theorem and to the
fact that the measure $\mu $ is invariant under the maps $T(y)$. But, as $%
\varepsilon \rightarrow 0$, we have the following well-known convergence
result: 
\begin{equation*}
\int_{Q}\varphi (x)g\left( \frac{x}{\varepsilon _{2}}\right) dx\rightarrow
\iint_{Q\times \Delta (A)}\varphi (x)\widehat{g}(s)dxd\beta \text{ as }%
\varepsilon \rightarrow 0.
\end{equation*}%
Hence the sequence 
\begin{equation*}
\int_{Q\times \Omega }f\left( x,T\left( \frac{x}{\varepsilon _{1}}\right)
\omega ,\frac{x}{\varepsilon _{2}}\right) dxd\mu \rightarrow \iint_{Q\times
\Omega \times \Delta (A)}\widehat{f}(x,\omega ,s)dxd\mu d\beta .
\end{equation*}

Now, let $f\in \mathcal{K}(Q;\mathcal{C}^{\infty }(\Omega ;A))$ and let $%
\eta >0$ be arbitrarily fixed. Let $K\subset Q$ be a compact set such that
supp$f\subset K$. By a density argument we choose $\phi $ in $\mathcal{C}%
_{0}^{\infty }(Q)\otimes \mathcal{C}^{\infty }(\Omega )\otimes A$ with supp$%
\phi \subset K$, such that $\left\Vert f-\phi \right\Vert _{\infty }\leq
\eta /(3\left\vert K\right\vert )$, $\left\vert K\right\vert $ denoting the
Lebesgue volume of $K$. By the decomposition 
\begin{equation*}
\begin{array}{l}
\int_{Q\times \Omega }f^{\varepsilon }dxd\mu -\iint_{Q\times \Omega \times
\Delta (A)}\widehat{f}dxd\mu d\beta =\int_{Q\times \Omega }(f^{\varepsilon
}-\phi ^{\varepsilon })dxd\mu \\ 
\ \ \ \ \ +\int_{Q\times \Omega }\phi ^{\varepsilon }dxd\mu -\iint_{Q\times
\Omega \times \Delta (A)}\widehat{\phi }dxd\mu d\beta +\iint_{Q\times \Omega
\times \Delta (A)}(\widehat{\phi }-\widehat{f})dxd\mu d\beta ,%
\end{array}%
\end{equation*}%
it follows readily that there exists $\varepsilon _{0}>0$ such that 
\begin{equation*}
\left\vert \int_{Q\times \Omega }f^{\varepsilon }dxd\mu -\iint_{Q\times
\Omega \times \Delta (A)}\widehat{f}dxd\mu d\beta \right\vert \leq \eta 
\text{\ for }0<\varepsilon \leq \varepsilon _{0}\text{.}
\end{equation*}%
This completes the proof.
\end{proof}

As a result, we have the following corollaries.

\begin{corollary}
\label{c3.1}Let $u\in \mathcal{K}(Q;\mathcal{C}^{\infty }(\Omega ;A))$ and
let $1\leq p<\infty $. Then, as $\varepsilon \rightarrow 0$,

\begin{itemize}
\item[(i)] $u^{\varepsilon }\rightarrow \varrho (u)$ stoch. in $%
L^{p}(Q\times \Omega )$-weak $\Sigma $, where $\varrho $ denote the
canonical mapping of $B_{A}^{p}$ into $\mathcal{B}_{A}^{p}$, and the
function $\varrho (u)$ is defined by $\varrho (u)(x,\omega )=\varrho
(u(x,\omega ))$ for a.e. $(x,\omega )\in Q\times \Omega $;

\item[(ii)] $\left\| u^{\varepsilon }\right\| _{L^{p}(Q\times \Omega
)}\rightarrow \left\| \widehat{u}\right\| _{L^{p}(Q\times \Omega \times
\Delta (A))}$.
\end{itemize}
\end{corollary}

\begin{proof}
(i) For each $f\in \mathcal{C}_{0}^{\infty }(Q)\otimes \mathcal{C}^{\infty
}(\Omega )\otimes A$ we have $uf\in \mathcal{K}(Q;\mathcal{C}^{\infty
}(\Omega ;A))$, hence part (i) follows readily by Proposition \ref{p3.3}.
For (ii), since $\mathcal{K}(Q;\mathcal{C}^{\infty }(\Omega ;A))$ is a
Banach algebra, it is easily shown that, for $1\leq p<\infty $ we have $%
\left\vert u\right\vert ^{p}\in \mathcal{K}(Q;\mathcal{C}^{\infty }(\Omega
;A))$ whenever $u\in \mathcal{K}(Q;\mathcal{C}^{\infty }(\Omega ;A))$, so
that once again by Proposition \ref{p3.3} we have, as $\varepsilon
\rightarrow 0$, 
\begin{equation*}
\int_{Q\times \Omega }\left\vert u^{\varepsilon }\right\vert ^{p}dxd\mu
\rightarrow \iint_{Q\times \Omega \times \Delta (A)}\left\vert \widehat{u}%
\right\vert ^{p}dxd\mu d\beta \text{.}
\end{equation*}%
The proof is complete.
\end{proof}

\begin{corollary}[\textit{Lower-semicontinuity property}]
\label{c3.2}Let $(u_{\varepsilon })_{\varepsilon >0}$ be a sequence in $%
L^{p}(Q\times \Omega )$ $(1\leq p<\infty )$ such that $u_{\varepsilon
}\rightarrow u_{0}$ stoch. in $L^{p}(Q\times \Omega )$-weak $\Sigma $ as $%
\varepsilon \rightarrow 0$, where $u_{0}\in L^{p}(Q\times \Omega ;\mathcal{B}%
_{A}^{p})$. Then 
\begin{equation}
\left\Vert u_{0}\right\Vert _{L^{p}(Q\times \Omega ;\mathcal{B}%
_{A}^{p})}\leq \underset{\varepsilon \rightarrow 0}{\lim \inf }\left\Vert
u_{\varepsilon }\right\Vert _{L^{p}(Q\times \Omega )}.\;\;\;\;  \label{3.8}
\end{equation}
\end{corollary}

\begin{proof}
Let $f\in \mathcal{C}_{0}^{\infty }(Q)\otimes \mathcal{C}^{\infty }(\Omega
)\otimes A$. We have 
\begin{equation}
\left\vert \int_{Q\times \Omega }u_{\varepsilon }f^{\varepsilon }dxd\mu
\right\vert \leq \left\Vert u_{\varepsilon }\right\Vert _{L^{p}(Q\times
\Omega )}\left\Vert f^{\varepsilon }\right\Vert _{L^{p^{\prime }}(Q\times
\Omega )}.  \label{3.81}
\end{equation}%
Then taking $\lim \inf_{\varepsilon \rightarrow 0}$ of both sides of (\ref%
{3.81}) and using the equality 
\begin{equation*}
\lim_{\varepsilon \rightarrow 0}\left\Vert f^{\varepsilon }\right\Vert
_{L^{p^{\prime }}(Q\times \Omega )}=\left\Vert \widehat{f}\right\Vert
_{L^{p^{\prime }}(Q\times \Omega \times \Delta (A))}\text{ (see part (ii) of
Corollary \ref{c3.1} above)}
\end{equation*}%
one arrives at 
\begin{equation}
\left\vert \iint_{Q\times \Omega \times \Delta (A)}\widehat{u}_{0}\widehat{f}%
dxd\mu d\beta \right\vert \leq \left\Vert \widehat{f}\right\Vert
_{L^{p^{\prime }}(Q\times \Omega \times \Delta (A))}\underset{\varepsilon
\rightarrow 0}{\lim \inf }\left\Vert u_{\varepsilon }\right\Vert
_{L^{p}(Q\times \Omega )}.  \label{3.82}
\end{equation}%
The space $\mathcal{G}(\mathcal{C}_{0}^{\infty }(Q)\otimes \mathcal{C}%
^{\infty }(\Omega )\otimes A)=\mathcal{C}_{0}^{\infty }(Q)\otimes \mathcal{C}%
^{\infty }(\Omega )\otimes \mathcal{C}(\Delta (A))$ being dense in $%
L^{p^{\prime }}(Q\times \Omega \times \Delta (A))$, (\ref{3.82}) still holds
with $v\in L^{p^{\prime }}(Q\times \Omega \times \Delta (A))$ instead of $%
\widehat{f}$. Consequently 
\begin{eqnarray*}
\left\Vert \widehat{u}_{0}\right\Vert _{L^{p}(Q\times \Omega \times \Delta
(A))} &=&\sup_{\left\Vert v\right\Vert _{L^{p^{\prime }}(Q\times \Omega
\times \Delta (A))}\leq 1}\left\vert \iint_{Q\times \Omega \times \Delta (A)}%
\widehat{u}_{0}vdxd\mu d\beta \right\vert \\
&\leq &\underset{\varepsilon \rightarrow 0}{\lim \inf }\left\Vert
u_{\varepsilon }\right\Vert _{L^{p}(Q\times \Omega )}\text{.}
\end{eqnarray*}%
The lemma follows.
\end{proof}

Throughout the paper the letter $E$ will denote any ordinary sequence $%
E=(\varepsilon _{n})$ (integers $n\geq 0$) with $0<\varepsilon _{n}\leq 1$
and $\varepsilon _{n}\rightarrow 0$ as $n\rightarrow \infty $. Such a
sequence will be termed a \textit{fundamental sequence}.

The usefulness of the next result will be brought to light in the sequel.
Prior to that we need one further definition.

\begin{definition}
\label{d3.1'}\emph{A function }$u\in L^{1}(Q\times \Omega ;B_{A}^{1})$\emph{%
\ is said to be }admissible\emph{\ if the trace function }$(x,\omega
)\mapsto u(x,T(x/\varepsilon _{1})\omega ,x/\varepsilon _{2})$\emph{\
(denoted by }$u^{\varepsilon }$\emph{), from }$Q\times \Omega $\emph{\ to }$%
\mathbb{C}$\emph{, is well-defined as an element of }$L^{1}(Q\times \Omega )$%
\emph{\ and satisfies the following convergence result:} 
\begin{equation}
\int_{Q\times \Omega }\left\vert u^{\varepsilon }\right\vert dxd\mu
\rightarrow \diint_{Q\times \Omega \times \Delta (A)}\left\vert \widehat{u}%
\right\vert dxd\mu d\beta \text{\ \emph{as} }\varepsilon \rightarrow 0\text{.%
}  \label{adm}
\end{equation}
\end{definition}

One can verify that any function in each of the following spaces is
admissible: $\mathcal{K}(Q;L^{p}(\Omega ;A))$ (the space of continuous
functions $f:\mathbb{R}^{N}\rightarrow L^{p}(\Omega ;A)$ with compact
support contained in $Q$, $1\leq p\leq \infty $), $\mathcal{C}(\overline{Q}%
;L^{\infty }(\Omega ;A))$ (for any bounded domain $Q$ in $\mathbb{R}^{N}$).

\begin{proposition}
\label{p3.4}Let $(u_{\varepsilon })_{\varepsilon \in E}\subset L^{p}(Q\times
\Omega )$ $(1<p<\infty )$ be a sequence which is weakly stochastically $%
\Sigma $-convergent in $L^{p}(Q\times \Omega )$ to some $u_{0}\in
L^{p}(Q\times \Omega ;\mathcal{B}_{A}^{p})$. Then as $E\ni \varepsilon
\rightarrow 0$ we have \emph{(\ref{3.1})} (in Definition \emph{\ref{d3.1}})
for any admissible function $f\in \mathcal{K}(Q;L^{p^{\prime }}(\Omega
;B_{A}^{p^{\prime },\infty }))$ where $B_{A}^{p^{\prime },\infty
}=B_{A}^{p^{\prime }}\cap L^{\infty }(\mathbb{R}^{N})$.
\end{proposition}

\begin{proof}
The space $\mathcal{K}(Q)\otimes \mathcal{C}^{\infty }(\Omega )\otimes
B_{A}^{p^{\prime },\infty }$ is dense in $\mathcal{K}(Q;L^{p^{\prime
}}(\Omega ;B_{A}^{p^{\prime },\infty }))$. Indeed $\mathcal{C}^{\infty
}(\Omega )\otimes B_{A}^{p^{\prime },\infty }$ is dense in $L^{p^{\prime
}}(\Omega ;B_{A}^{p^{\prime },\infty })$, so that by \cite[p. 46]{Bou1}, our
claim is justified. With this in mind, we firstly check (\ref{3.1}) for $%
f\in \mathcal{K}(Q)\otimes \mathcal{C}^{\infty }(\Omega )\otimes
B_{A}^{p^{\prime },\infty }$. It suffices to verify this for $f$ under the
form 
\begin{eqnarray*}
f(x,\omega ,y) &=&\varphi (x)\psi (\omega )v(y)\text{\ \ }(x\in Q,\omega \in
\Omega ,y\in \mathbb{R}^{N})\text{\ with} \\
\varphi &\in &\mathcal{K}(Q),\psi \in \mathcal{C}^{\infty }(\Omega )\text{
and }v\in B_{A}^{p^{\prime },\infty }\text{.}
\end{eqnarray*}%
Let $f$ be as above. Let $\delta >0$ be freely fixed, and let $w\in A$ be
such that $\left\Vert v-w\right\Vert _{p^{\prime }}\leq \delta $ (where we
have used here the density of $A$ in $B_{A}^{p^{\prime }}$). Set 
\begin{equation*}
g(x,\omega ,y)=\varphi (x)\psi (\omega )w(y)\text{\ \ }(x\in Q,\omega \in
\Omega ,y\in \mathbb{R}^{N}),
\end{equation*}%
which gives a function $g\in \mathcal{K}(Q)\otimes \mathcal{C}^{\infty
}(\Omega )\otimes A$. We have 
\begin{eqnarray*}
&&\int_{Q\times \Omega }u_{\varepsilon }f^{\varepsilon }dxd\mu
-\iint_{Q\times \Omega \times \Delta (A)}\widehat{u}_{0}\widehat{f}dxd\mu
d\beta \\
&=&\int_{Q\times \Omega }u_{\varepsilon }\varphi (x)\psi (T(x/\varepsilon
_{1})\omega )[v^{\varepsilon }(x/\varepsilon _{2})-w(x/\varepsilon
_{2})]dxd\mu \\
&&+\int_{Q\times \Omega }u_{\varepsilon }g^{\varepsilon }dxd\mu
-\iint_{Q\times \Omega \times \Delta (A)}\widehat{u}_{0}\widehat{g}dxd\mu
d\beta \\
&&+\iint_{Q\times \Omega \times \Delta (A)}\widehat{u}_{0}\varphi \psi (%
\widehat{w}-\widehat{v})dxd\mu d\beta \\
&=&(I)+(II)+(III).
\end{eqnarray*}%
As far as $(I)$ is concerned, we have 
\begin{equation*}
\left\vert (I)\right\vert \leq \left\Vert \varphi \right\Vert _{\infty
}\left\Vert \psi \right\Vert _{\infty }\left\Vert u_{\varepsilon
}\right\Vert _{L^{p}(Q\times \Omega )}\left( \int_{K}\left\vert
v^{\varepsilon }-w^{\varepsilon }\right\vert ^{p^{\prime }}dx\right)
^{1/p^{\prime }}
\end{equation*}%
where $K$ is a compact subset of $\mathbb{R}^{N}$ containing the support of $%
\varphi $. But $v$ and $w$ possess mean value, so that, as $\varepsilon
\rightarrow 0$, 
\begin{equation*}
\int_{K}\left\vert v^{\varepsilon }-w^{\varepsilon }\right\vert ^{p^{\prime
}}dx\rightarrow M(\left\vert v-w\right\vert ^{p^{\prime }})\left\vert
K\right\vert \text{ since }\left\vert v-w\right\vert ^{p^{\prime }}\in
B_{A}^{1}\text{,}
\end{equation*}%
$\left\vert K\right\vert $ denoting the Lebesgue measure of $K$. In view of
the equality $\left\Vert v-w\right\Vert _{p^{\prime }}=[M(\left\vert
v-w\right\vert ^{p^{\prime }})]^{1/p^{\prime }}$, we have $\lim_{E\ni
\varepsilon \rightarrow 0}\left\vert (I)\right\vert \leq c\delta $ where $c$
is a positive constant independent of $\delta $. For $(III)$, we have 
\begin{eqnarray*}
&&\left\vert \iint_{Q\times \Omega \times \Delta (A)}\widehat{u}_{0}\varphi
\psi (\widehat{w}-\widehat{v})dxd\mu d\beta \right\vert \\
&\leq &\left\Vert \widehat{u}_{0}\right\Vert _{L^{p}(Q\times \Omega \times
\Delta (A))}\left\Vert \varphi \right\Vert _{\infty }\left\Vert \psi
\right\Vert _{\infty }\left\Vert \widehat{w}-\widehat{v}\right\Vert
_{L^{p^{\prime }}(\Delta (A))} \\
&=&c\left\Vert v-w\right\Vert _{p^{\prime }}\leq c\delta
\end{eqnarray*}%
where $c=\left\Vert \widehat{u}_{0}\right\Vert _{L^{p}(Q\times \Omega \times
\Delta (A))}\left\Vert \varphi \right\Vert _{\infty }\left\Vert \psi
\right\Vert _{\infty }$. Next, since 
\begin{equation*}
\int_{Q\times \Omega }u_{\varepsilon }g^{\varepsilon }dxd\mu \rightarrow
\iint_{Q\times \Omega \times \Delta (A)}\widehat{u}_{0}\widehat{g}dxd\mu
d\beta
\end{equation*}%
it follows that 
\begin{equation*}
\lim_{E\ni \varepsilon \rightarrow 0}\left\vert \int_{Q\times \Omega
}u_{\varepsilon }\varphi (x)\psi (T(x/\varepsilon _{1})\omega
)w(x/\varepsilon _{2})dxd\mu -\iint_{Q\times \Omega \times \Delta (A)}%
\widehat{u}_{0}\varphi \psi \widehat{w}dxd\mu d\beta \right\vert \leq c\delta
\end{equation*}%
where $c>0$ is independent of $\delta $, hence (\ref{3.1}) follows with the
above taken $f$, since $\delta $ is arbitrary. In view of the density of $%
\mathcal{K}(Q)\otimes \mathcal{C}^{\infty }(\Omega )\otimes B_{A}^{p^{\prime
},\infty }$ in $\mathcal{K}(Q;L^{p^{\prime }}(\Omega ;B_{A}^{p^{\prime
},\infty }))$ the result follows by repeating the same way of proceeding as
done above.
\end{proof}

The next result is a mere consequence of the preceding result. Its easy
proof is left to the reader.

\begin{corollary}
\label{c3.3}Let $u\in \mathcal{K}(Q;L^{\infty }(\Omega ;B_{A}^{p,\infty }))$ 
$(1<p<\infty )$ be an admissible function in the sense of Definition \emph{%
\ref{d3.1'}}. Then the sequence $(u^{\varepsilon })_{\varepsilon >0}$ is
weakly stochastically $\Sigma $-convergent in $L^{p}(Q\times \Omega )$ to $%
\varrho (u)$.
\end{corollary}

The following result is the point of departure of all the compactness
results involved in this paper.

\begin{theorem}
\label{t3.1}Any bounded sequence $(u_{\varepsilon })_{\varepsilon \in E}$ in 
$L^{p}(Q\times \Omega )$ (where $E$ is a fundamental sequence and $%
1<p<\infty $) admits a subsequence which is weakly stochastically $\Sigma $%
-convergent in $L^{p}(Q\times \Omega )$.
\end{theorem}

Its proof relies on the following result whose proof can be found in \cite%
{CMP}.

\begin{proposition}
\label{p3.1}Let $X$ be a subspace (not necessarily closed) of a reflexive
Banach space $Y$ and let $f_{n}:X\rightarrow \mathbb{C}$ be a sequence of
linear functionals (not necessarily continuous). Assume there exists a
constant $c>0$ such that 
\begin{equation}
\underset{n}{\lim \sup }\left\vert f_{n}(x)\right\vert \leq c\left\Vert
x\right\Vert \;\;\text{for all }x\in X.  \label{3.2}
\end{equation}%
where $\left\Vert \cdot \right\Vert $ denotes the norm in $Y$. Then there
exist a subsequence $(f_{n_{k}})_{k}$ of $(f_{n})$ and a functional $f\in
Y^{\prime }$ such that $\lim_{k}f_{n_{k}}(x)=f(x)$ for all $x\in X$.
\end{proposition}

\begin{proof}[\textit{Proof of Theorem} \ref{t3.1}]
Let $Y=L^{p^{\prime }}(Q\times \Omega \times \Delta (A))$, $X=\mathcal{C}%
_{0}^{\infty }(Q)\otimes \mathcal{C}^{\infty }(\Omega )\otimes \mathcal{C}%
(\Delta (A))$. Let us define the mapping $L_{\varepsilon }$ by 
\begin{equation*}
L_{\varepsilon }(\widehat{f})=\int_{Q\times \Omega }u_{\varepsilon
}f^{\varepsilon }dxd\mu \text{\ \ (}\widehat{f}\in \mathcal{C}_{0}^{\infty
}(Q)\otimes \mathcal{C}^{\infty }(\Omega )\otimes \mathcal{C}(\Delta (A))=%
\mathcal{G}(\mathcal{C}_{0}^{\infty }(Q)\otimes \mathcal{C}^{\infty }(\Omega
)\otimes A)\text{).}
\end{equation*}%
where $f^{\varepsilon }(x,\omega )=f(x,T(x/\varepsilon _{1})\omega
,x/\varepsilon _{2})$ for $(x,\omega )\in Q\times \Omega $. Then 
\begin{equation*}
\underset{\varepsilon }{\lim \sup }\left\vert L_{\varepsilon }(\widehat{f}%
)\right\vert \leq c\left\Vert \widehat{f}\right\Vert _{L^{p^{\prime
}}(Q\times \Omega \times \Delta (A))}\text{\ for all }\widehat{f}\in X.
\end{equation*}%
Indeed one has the inequality $\left\vert L_{\varepsilon }(f)\right\vert
\leq c\left\Vert f^{\varepsilon }\right\Vert _{L^{p^{\prime }}(Q\times
\Omega )}$ and thus, as $\varepsilon \rightarrow 0$, $\left\Vert
f^{\varepsilon }\right\Vert _{L^{p^{\prime }}(Q\times \Omega )}\rightarrow
\left\Vert \widehat{f}\right\Vert _{L^{p^{\prime }}(Q\times \Omega \times
\Delta (A))}$ (see Corollary \ref{c3.1}). We therefore apply Proposition \ref%
{p3.1} with the above notation to get the existence of a subsequence $%
E^{\prime }$ of $E$ and of a unique $v_{0}\in L^{p}(Q\times \Omega \times
\Delta (A))$ such that 
\begin{equation*}
\int_{Q\times \Omega }u_{\varepsilon }f^{\varepsilon }dxd\mu \rightarrow
\iint_{Q\times \Omega \times \Delta (A)}v_{0}(x,\omega ,s)\widehat{f}%
(x,\omega ,s)dxd\mu d\beta
\end{equation*}%
\ for all $\widehat{f}\in X$. But $v_{0}=\mathcal{G}_{1}\circ u_{0}$ where $%
u_{0}\in L^{p}(Q\times \Omega ;\mathcal{B}_{A}^{p})$, and so the result
follows.
\end{proof}

In order to deal with the convergence of a product of sequences we need to
define the concept of strong stochastic $\Sigma $-convergence.

\begin{definition}
\label{d3.2}\emph{A sequence }$(u_{\varepsilon })_{\varepsilon >0}\subset
L^{p}(Q\times \Omega )$\emph{\ }$(1\leq p<\infty )$\emph{\ is said to }%
strongly stochastically $\Sigma $-converge\emph{\ in }$L^{p}(Q\times \Omega
) $\emph{\ to some }$u_{0}\in L^{p}(Q\times \Omega ;\mathcal{B}_{A}^{p})$%
\emph{\ if it is weakly stochastically }$\Sigma $\emph{-convergent and
further satisfies the following condition: }%
\begin{equation}
\left\| u_{\varepsilon }\right\| _{L^{p}(Q\times \Omega )}\rightarrow
\left\| \widehat{u}_{0}\right\| _{L^{p}(Q\times \Omega \times \Delta (A))}.
\label{3.12}
\end{equation}%
\emph{We denote this by }$u_{\varepsilon }\rightarrow u_{0}$\emph{\ stoch.
in }$L^{p}(Q\times \Omega )$\emph{-strong }$\Sigma $\emph{.}
\end{definition}

\begin{remark}
\label{r3.1}\emph{(1) By the above definition, the uniqueness of the limit
of such a sequence is ensured. (2) By the Corollary \ref{c3.1} it is
immediate that for any }$u\in \mathcal{K}(Q;\mathcal{C}^{\infty }(\Omega
;A)) $\emph{, the sequence }$(u^{\varepsilon })_{\varepsilon >0}$\emph{\ is
strongly stochastically }$\Sigma $\emph{-convergent to }$\varrho (u)$\emph{.}
\end{remark}

The next result will be of capital interest in the homogenization process.

\begin{theorem}
\label{t3.2}Let $1<p,q<\infty $ and $r\geq 1$ be such that $1/r=1/p+1/q\leq
1 $. Assume $(u_{\varepsilon })_{\varepsilon \in E}\subset L^{q}(Q\times
\Omega )$ is weakly stochastically $\Sigma $-convergent in $L^{q}(Q\times
\Omega )$ to some $u_{0}\in L^{q}(Q\times \Omega ;\mathcal{B}_{A}^{q})$, and 
$(v_{\varepsilon })_{\varepsilon \in E}\subset L^{p}(Q\times \Omega )$ is
strongly stochastically $\Sigma $-convergent in $L^{p}(Q\times \Omega )$ to
some $v_{0}\in L^{p}(Q\times \Omega ;\mathcal{B}_{A}^{p})$. Then the
sequence $(u_{\varepsilon }v_{\varepsilon })_{\varepsilon \in E}$ is weakly
stochastically $\Sigma $-convergent in $L^{r}(Q\times \Omega )$ to $%
u_{0}v_{0}$.
\end{theorem}

\begin{proof}
We assume without lost of generality that our sequences are real values.
This assumption is fully motivated by the fact that in general, almost only
linear problems are of complex coefficients, and so in that case, the
linearity permits to work with real coefficients. This being so, we will
deeply make use of the following simple inequalities proved in \cite{Zhikov1}%
: 
\begin{equation}
\begin{array}{l}
0\leq \left\vert a+tb\right\vert ^{p}-\left\vert a\right\vert
^{p}-pt\left\vert a\right\vert ^{p-2}ab\leq c\left\vert t\right\vert
^{1+s}(\left\vert a\right\vert ^{p}+\left\vert b\right\vert ^{p}) \\ 
\text{for each }\left\vert t\right\vert \leq 1\text{ and for every }a,b\in 
\mathbb{R}\text{, where } \\ 
s=\min (p-1,1)>0\text{ and }c>0\text{ is independent of }a,b\text{.}%
\end{array}
\label{3.13}
\end{equation}%
We proceed in two steps.\medskip

\noindent \textbf{Step 1.} Set $p^{\prime }=p/(p-1)$, and let us first show
that the sequence $z_{\varepsilon }=\left\vert v_{\varepsilon }\right\vert
^{p-2}v_{\varepsilon }$ is weakly stochastically $\Sigma $-convergent to $%
\left\vert v_{0}\right\vert ^{p-2}v_{0}$ in $L^{p^{\prime }}(Q\times \Omega
) $. To this end, let $z\in L^{p^{\prime }}(Q\times \Omega ;\mathcal{B}%
_{A}^{p^{\prime }})$ denote the weak stochastic $\Sigma $-limit of $%
(z_{\varepsilon })_{\varepsilon \in E}$ in $L^{p^{\prime }}(Q\times \Omega )$
(up to a subsequence if necessary; in fact it is easily seen that this
sequence is bounded in $L^{p^{\prime }}(Q\times \Omega )$). Let $\varphi \in 
\mathcal{C}_{0}^{\infty }(Q)\otimes \mathcal{C}^{\infty }(\Omega )\otimes A$
with $\left\Vert \varphi \right\Vert _{L^{p}(Q\times \Omega ;A)}\leq 1$. We
have by the second inequality in (\ref{3.13}) that 
\begin{eqnarray*}
\int_{Q\times \Omega }\left\vert v_{\varepsilon }+t\varphi ^{\varepsilon
}\right\vert ^{p}dxd\mu &\leq &\int_{Q\times \Omega }\left\vert
v_{\varepsilon }\right\vert ^{p}dxd\mu +pt\int_{Q\times \Omega
}z_{\varepsilon }\varphi ^{\varepsilon }dxd\mu \\
&&+c_{1}\left\vert t\right\vert ^{1+s}
\end{eqnarray*}%
for $\left\vert t\right\vert \leq 1$, $c_{1}$ being a positive constant
independent of $\varepsilon $ (since the sequence $(v_{\varepsilon
})_{\varepsilon \in E}$ is bounded in $L^{p}(Q\times \Omega )$). Taking the $%
\lim \inf_{E\ni \varepsilon \rightarrow 0}$ in the above inequality we get,
by virtue of (\ref{3.12}) (in Definition \ref{d3.2}) and the lower
semicontinuity property (\ref{3.8}) (in Corollary \ref{c3.2}) that 
\begin{eqnarray*}
\iint_{Q\times \Omega \times \Delta (A)}\left\vert \widehat{v}_{0}+t\widehat{%
\varphi }\right\vert ^{p}dxd\mu d\beta &\leq &\iint_{Q\times \Omega \times
\Delta (A)}\left\vert \widehat{v}_{0}\right\vert ^{p}dxd\mu d\beta \\
&&+pt\iint_{Q\times \Omega \times \Delta (A)}\widehat{z}\widehat{\varphi }%
dxd\mu d\beta +c_{1}\left\vert t\right\vert ^{1+s}.
\end{eqnarray*}%
On the other hand, the first inequality in (\ref{3.13}) yields 
\begin{eqnarray*}
\iint_{Q\times \Omega \times \Delta (A)}\left\vert \widehat{v}_{0}+t\widehat{%
\varphi }\right\vert ^{p}dxd\mu d\beta &\geq &\iint_{Q\times \Omega \times
\Delta (A)}\left\vert \widehat{v}_{0}\right\vert ^{p}dxd\mu d\beta \\
&&+pt\iint_{Q\times \Omega \times \Delta (A)}\left\vert \widehat{v}%
_{0}\right\vert ^{p-2}\widehat{v}_{0}dxd\mu d\beta ,
\end{eqnarray*}%
hence 
\begin{equation*}
pt\iint_{Q\times \Omega \times \Delta (A)}\left\vert \widehat{v}%
_{0}\right\vert ^{p-2}\widehat{v}_{0}dxd\mu d\beta \leq pt\iint_{Q\times
\Omega \times \Delta (A)}\widehat{z}\widehat{\varphi }dxd\mu d\beta
+c_{1}\left\vert t\right\vert ^{1+s}.
\end{equation*}%
Now, taking in the above inequality $\varphi =\psi /\left\Vert \psi
\right\Vert _{L^{p}(Q\times \Omega ;A)}$ for any arbitrary $\psi \in 
\mathcal{C}_{0}^{\infty }(Q)\otimes \mathcal{C}^{\infty }(\Omega )\otimes A$
the same inequality holds for any arbitrary $\psi $ in place of $\varphi $,
which, together with the arbitrariness of the real number $t$ in $\left\vert
t\right\vert \leq 1$, gives $z=\left\vert v_{0}\right\vert ^{p-2}v_{0}$%
.\medskip

\noindent \textbf{Step 2}. Now, let us establish the convergence result $%
u_{\varepsilon }v_{\varepsilon }\rightarrow u_{0}v_{0}$ stoch. in $%
L^{r}(Q\times \Omega )$-weak $\Sigma $. First of all the sequence $%
(u_{\varepsilon }v_{\varepsilon })_{\varepsilon \in E}$ is bounded in $%
L^{r}(Q\times \Omega )$. Next, let $\varphi \in \mathcal{C}_{0}^{\infty
}(Q)\otimes \mathcal{C}^{\infty }(\Omega )\otimes A$ and set 
\begin{equation*}
\ell =\lim_{E\ni \varepsilon \rightarrow 0}\int_{Q\times \Omega
}u_{\varepsilon }v_{\varepsilon }\varphi ^{\varepsilon }dxd\mu \text{
(possibly up to a subsequence).}
\end{equation*}%
We need to show that $\ell =\iint_{Q\times \Omega \times \Delta (A)}\widehat{%
u}_{0}\widehat{v}_{0}\widehat{\varphi }dxd\mu d\beta $. First and foremost
we have $\varphi ^{\varepsilon }\in L^{r^{\prime }}(Q\times \Omega )$ and
so, $u_{\varepsilon }\varphi ^{\varepsilon }\in L^{p^{\prime }}(Q\times
\Omega )$ since $1/r^{\prime }+1/q=1/p^{\prime }$ and $u_{\varepsilon }\in
L^{q}(Q\times \Omega )$. Thus, once again by the second inequality in (\ref%
{3.13}) and keeping in mind the definition of $z_{\varepsilon }$ in Step 1,
one has 
\begin{eqnarray*}
\int_{Q\times \Omega }\left\vert z_{\varepsilon }+tu_{\varepsilon }\varphi
^{\varepsilon }\right\vert ^{p^{\prime }}dxd\mu &\leq &\int_{Q\times \Omega
}\left\vert z_{\varepsilon }\right\vert ^{p^{\prime }}dxd\mu +p^{\prime
}t\int_{Q\times \Omega }\left\vert z_{\varepsilon }\right\vert
^{p-2}z_{\varepsilon }u_{\varepsilon }\varphi ^{\varepsilon }dxd\mu \\
&&+c_{1}\left\vert t\right\vert ^{1+s} \\
&=&\int_{Q\times \Omega }\left\vert v_{\varepsilon }\right\vert ^{p}dxd\mu
+p^{\prime }t\int_{Q\times \Omega }v_{\varepsilon }u_{\varepsilon }\varphi
^{\varepsilon }dxd\mu +c_{1}\left\vert t\right\vert ^{1+s}
\end{eqnarray*}%
since $v_{\varepsilon }=\left\vert z_{\varepsilon }\right\vert ^{p^{\prime
}-2}z_{\varepsilon }$ and $\left\vert z_{\varepsilon }\right\vert
^{p^{\prime }}=\left\vert v_{\varepsilon }\right\vert ^{p}$. On the other
hand, one easily sees that the sequence $(u_{\varepsilon }\varphi
^{\varepsilon })_{\varepsilon \in E}$ is weakly stochastically $\Sigma $%
-convergent to $u_{0}\varrho (\varphi )$ in $L^{p^{\prime }}(Q\times \Omega
) $, so that, passing to the limit in the above inequality, using the lower
semicontinuity property (\ref{3.8}), we get 
\begin{eqnarray*}
\iint_{Q\times \Omega \times \Delta (A)}\left\vert \widehat{z}+t\widehat{u}%
_{0}\widehat{\varphi }\right\vert ^{p^{\prime }}dxd\mu d\beta &\leq
&\iint_{Q\times \Omega \times \Delta (A)}\left\vert \widehat{v}%
_{0}\right\vert ^{p}dxd\mu d\beta +p^{\prime }t\ell +c_{1}\left\vert
t\right\vert ^{1+s} \\
&=&\iint_{Q\times \Omega \times \Delta (A)}\left\vert \widehat{z}\right\vert
^{p^{\prime }}dxd\mu d\beta +p^{\prime }t\ell +c_{1}\left\vert t\right\vert
^{1+s},
\end{eqnarray*}%
since $z=\left\vert v_{0}\right\vert ^{p-2}v_{0}$ (as shown in Step 1), and
therefore, $\left\vert v_{0}\right\vert ^{p}=\left\vert z\right\vert
^{p^{\prime }}$. Besides, we have by the first inequality in (\ref{3.13})
that 
\begin{eqnarray*}
\iint_{Q\times \Omega \times \Delta (A)}\left\vert \widehat{z}+t\widehat{u}%
_{0}\widehat{\varphi }\right\vert ^{p^{\prime }}dxd\mu d\beta &\geq
&\iint_{Q\times \Omega \times \Delta (A)}\left\vert \widehat{z}\right\vert
^{p^{\prime }}dxd\mu d\beta \\
&&+p^{\prime }t\iint_{Q\times \Omega \times \Delta (A)}\left\vert \widehat{z}%
\right\vert ^{p^{\prime }-2}\widehat{z}\widehat{u}_{0}\widehat{\varphi }%
dxd\mu d\beta \\
&=&\iint_{Q\times \Omega \times \Delta (A)}\left\vert \widehat{z}\right\vert
^{p^{\prime }}dxd\mu d\beta \\
&&+p^{\prime }t\iint_{Q\times \Omega \times \Delta (A)}\widehat{v}_{0}%
\widehat{u}_{0}\widehat{\varphi }dxd\mu d\beta
\end{eqnarray*}%
since $v_{0}=\left\vert z\right\vert ^{p^{\prime }-2}z$. We are therefore
led to 
\begin{equation*}
p^{\prime }t\iint_{Q\times \Omega \times \Delta (A)}\widehat{v}_{0}\widehat{u%
}_{0}\widehat{\varphi }dxd\mu d\beta \leq p^{\prime }t\ell +c_{1}\left\vert
t\right\vert ^{1+s}\;\forall \left\vert t\right\vert \leq 1,
\end{equation*}%
hence $\ell =\iint_{Q\times \Omega \times \Delta (A)}\widehat{v}_{0}\widehat{%
u}_{0}\widehat{\varphi }dxd\mu d\beta $.
\end{proof}

The following result will be of great interest in practise. It is a mere
consequence of the preceding theorem.

\begin{corollary}
\label{c3.4}Let $(u_{\varepsilon })_{\varepsilon \in E}\subset L^{p}(Q\times
\Omega )$ and $(v_{\varepsilon })_{\varepsilon \in E}\subset L^{p^{\prime
}}(Q\times \Omega )\cap L^{\infty }(Q\times \Omega )$ ($1<p<\infty $ and $%
p^{\prime }=p/(p-1)$) be two sequences such that:

\begin{itemize}
\item[(i)] $u_{\varepsilon }\rightarrow u_{0}$ stoch. in $L^{p}(Q\times
\Omega )$-weak $\Sigma $;

\item[(ii)] $v_{\varepsilon }\rightarrow v_{0}$ stoch. in $L^{p^{\prime
}}(Q\times \Omega )$-strong $\Sigma $;

\item[(iii)] $(v_{\varepsilon })_{\varepsilon \in E}$ is bounded in $%
L^{\infty }(Q\times \Omega )$.
\end{itemize}

\noindent Then $u_{\varepsilon }v_{\varepsilon }\rightarrow u_{0}v_{0}$
stoch. in $L^{p}(Q\times \Omega )$-weak $\Sigma $.
\end{corollary}

\begin{proof}
By Theorem \ref{t3.2}, the sequence $(u_{\varepsilon }v_{\varepsilon
})_{\varepsilon \in E}$ weakly stochastically $\Sigma $-converges towards $%
u_{0}v_{0}$ in $L^{1}(Q\times \Omega )$. Besides the same sequence is
bounded in $L^{p}(Q\times \Omega )$ so that by the Theorem \ref{t3.1}, it
weakly stochastically $\Sigma $-converges in $L^{p}(Q\times \Omega )$
towards some $w_{0}\in L^{p}(Q\times \Omega ;\mathcal{B}_{A}^{p})$. This
gives as a result $w_{0}=u_{0}v_{0}$.
\end{proof}

The strong stochastic $\Sigma $-convergence is a generalization of the
strong convergence as one can easily see in the following result whose easy
proof is left to the reader.

\begin{proposition}
\label{p3.5}Let $(u_{\varepsilon })_{\varepsilon \in E}\subset L^{p}(Q\times
\Omega )$ $(1\leq p<\infty )$ be a strongly convergent sequence in $%
L^{p}(Q\times \Omega )$ to some $u_{0}\in L^{p}(Q\times \Omega )$. Then $%
(u_{\varepsilon })_{\varepsilon \in E}$ strongly stochastically $\Sigma $%
-converges in $L^{p}(Q\times \Omega )$ towards $u_{0}$.
\end{proposition}

In the first step of the proof of Theorem \ref{t3.2} we have proven the
following assertion: If $v_{\varepsilon }\rightarrow v_{0}$ stoch. in $%
L^{p}(Q\times \Omega )$-strong $\Sigma $ then $\left\vert v_{\varepsilon
}\right\vert ^{p-2}v_{\varepsilon }\rightarrow \left\vert v_{0}\right\vert
^{p-2}v_{0}$ stoch. in $L^{p^{\prime }}(Q\times \Omega )$-weak $\Sigma $.
One can weaken the above strong convergence condition and obtain, under an
additional weak convergence assumption, the following result: If $%
u_{\varepsilon }\rightarrow u_{0}$ stoch. in $L^{p}(Q\times \Omega )$-weak $%
\Sigma $ and $\left\vert u_{\varepsilon }\right\vert ^{p-2}u_{\varepsilon
}\rightarrow v_{0}$ stoch. in $L^{p^{\prime }}(Q\times \Omega )$-weak $%
\Sigma $, then 
\begin{equation*}
\iint_{Q\times \Omega \times \Delta (A)}\widehat{u}_{0}\widehat{v}_{0}dxd\mu
d\beta \leq \underset{\varepsilon \rightarrow 0}{\lim \inf }\int_{Q\times
\Omega }\left\vert u_{\varepsilon }\right\vert ^{p}dxd\mu \text{.}
\end{equation*}%
Moreover if the above inequality holds as an equality, then $%
v_{0}=\left\vert u_{0}\right\vert ^{p-2}u_{0}$.

The above result is a particular case of a general situation stated in the
following

\begin{theorem}
\label{t5}Let $(x,\omega ,y,\lambda )\mapsto a(x,\omega ,y,\lambda )$, from $%
\overline{Q}\times \Omega \times \mathbb{R}^{N}\times \mathbb{R}^{m}$ to $%
\mathbb{R}^{m}$ be a vector-valued function which is of Carath\'{e}odory's
type, i.e., \emph{(i)} and \emph{(ii)} below are satisfied:

\begin{itemize}
\item[(i)] $a(x,\cdot ,y,\lambda )$ is $d\mu $-measurable for any $%
(x,y,\lambda )\in \overline{Q}\times \mathbb{R}^{N}\times \mathbb{R}^{m}$

\item[(ii)] $a(\cdot ,\omega ,\cdot ,\cdot )$ is continuous for $d\mu $%
-almost all $\omega \in \Omega $,
\end{itemize}

\noindent and further satisfies the following conditions:

\begin{itemize}
\item[(iii)] $\left\vert a(x,\omega ,y,\lambda )\right\vert \leq
c(\left\vert \lambda \right\vert ^{p-1}+1)$

\item[(iv)] $\left( a(x,\omega ,y,\lambda )-a(x,\omega ,y,\lambda ^{\prime
})\right) \cdot (\lambda -\lambda ^{\prime })\geq 0$

\item[(v)] $a(x,\omega ,\cdot ,\lambda )\in (A)^{m}$
\end{itemize}

\noindent for all $(x,y)\in \overline{Q}\times \mathbb{R}^{N}$, all $\lambda
,\lambda ^{\prime }\in \mathbb{R}^{m}$ and for $d\mu $-almost all $\omega
\in \Omega $, where $c$ is a positive constant independent of $(x,\omega
,y,\lambda )$. Finally let $(v_{\varepsilon })_{\varepsilon \in E}\subset
L^{p}(Q\times \Omega )^{m}$ be a sequence which componentwise weakly
stochastically $\Sigma $-converges towards $v_{0}\in L^{p}(Q\times \Omega ;(%
\mathcal{B}_{A}^{p})^{m})$ as $E\ni \varepsilon \rightarrow 0$. Then the
sequence $(a^{\varepsilon }(\cdot ,v_{\varepsilon }))_{\varepsilon \in E}$
defined by $a^{\varepsilon }(\cdot ,v_{\varepsilon })(x,\omega
)=a(x,T(x/\varepsilon _{1})\omega ,x/\varepsilon _{2},v_{\varepsilon
}(x,\omega ))$ for $(x,\omega )\in Q\times \Omega $, is weakly
stochastically $\Sigma $-convergent in $L^{p^{\prime }}(Q\times \Omega )^{m}$
(up to a subsequence) to some $z_{0}\in L^{p^{\prime }}(Q\times \Omega ;(%
\mathcal{B}_{A}^{p^{\prime }})^{m})$ such that 
\begin{equation}
\iint_{Q\times \Omega \times \Delta (A)}\widehat{z}_{0}\widehat{v}_{0}dxd\mu
d\beta \leq \underset{E\ni \varepsilon \rightarrow 0}{\lim \inf }%
\int_{Q\times \Omega }a^{\varepsilon }(\cdot ,v_{\varepsilon })\cdot
v_{\varepsilon }dxd\mu .  \label{3.14}
\end{equation}%
Moreover if \emph{(\ref{3.14})} holds as an equality, then $z_{0}(x,\omega
,y)=a(x,\omega ,y,v_{0}(x,\omega ,y))$.
\end{theorem}

We will make use of the following lemma.

\begin{lemma}
\label{l6}Let $F_{1}$ and $F_{2}$ be two Banach spaces, $(Y,\mathfrak{M},\mu
)$ a measure space, $X$ a $\mu $-measurable subset of $Y$, and $g:X\times
F_{1}\rightarrow F_{2}$ a Carath\'{e}odory mapping. For each measurable
function $u:X\rightarrow F_{1}$, let $G(u)$ be the measurable function $%
x\mapsto g(x,u(x))$, from $X$ to $F_{2}$. If $G:u\mapsto G(u)$ maps $%
L^{p}(X;F_{1})$ into $L^{r}(X;F_{2})$ ($1\leq p,r<\infty $) then $G$ is
continuous in the norm topology.
\end{lemma}

\begin{proof}
A look at the proof of \cite[Chap. IV, Proposition 1.1]{18} shows that one
can replace in that proof, the Borel subset $\Omega $ of $\mathbb{R}^{n}$ by
the measurable subset $X$ of $Y$, $E$ by $F_{1}$, $F$ by $F_{2}$, and get
readily our result.
\end{proof}

\begin{proof}[\textit{Proof of Theorem} \ref{t5}]
By (iii) the sequence $(a^{\varepsilon }(\cdot ,v_{\varepsilon
}))_{\varepsilon \in E}$ is bounded in $L^{p^{\prime }}(Q\times \Omega )^{m}$%
, thus there exists a subsequence $E^{\prime }$ from $E$ and a function $%
z_{0}\in L^{p^{\prime }}(Q\times \Omega ;(\mathcal{B}_{A}^{p^{\prime
}})^{m}) $ such that $a^{\varepsilon }(\cdot ,v_{\varepsilon })\rightarrow
z_{0}$ stoch. in $L^{p^{\prime }}(Q\times \Omega )^{m}$-weak $\Sigma $ as $%
E\ni \varepsilon \rightarrow 0$. Let us show (\ref{3.14}). For that purpose,
let $\psi \in \lbrack \mathcal{C}_{0}^{\infty }(Q)\otimes \mathcal{C}%
^{\infty }(\Omega )\otimes A]^{m}$ (which is dense in $L^{p}(Q\times \Omega
;A)^{m}$); then the function $(x,\omega ,y)\mapsto a(x,\omega ,y,\psi
(x,\omega ,y))$ lies in $\mathcal{C}(\overline{Q};L^{\infty }(\Omega
;A))^{m} $. Indeed, as a result of (ii), the function $a(\cdot ,\omega
,y,\psi (\cdot ,\omega ,y))$ is continuous. Moreover for each fixed $x\in 
\overline{Q}$, $a(x,\cdot ,\cdot ,\psi (x,\cdot ,\cdot ))\in L^{\infty
}(\Omega ;A)^{m}$: in fact for any $y\in \mathbb{R}^{N}$ we have $\left\vert
a(x,\cdot ,y,\psi (x,\cdot ,y))\right\vert \leq c_{1}$ where $%
c_{1}=c(1+\left\Vert \psi \right\Vert _{\infty }^{p-1})$ and the function $%
a(x,\cdot ,y,\psi (x,\cdot ,y))$ is $\mu $-measurable; furthermore, for $\mu 
$-a.e. $\omega \in \Omega $, the function $a(x,\omega ,\cdot ,\psi (x,\omega
,\cdot ))$ belongs to $(A)^{m}$. In fact $\psi (x,\omega ,\cdot )\in (A)^{m}$%
, and it suffices to check that $a(x,\omega ,\cdot ,\phi )\in (A)^{m}$ for
any $\phi \in (A)^{m}$. But since the function $\phi $ is bounded, let $%
K\subset \mathbb{R}^{m}$ be a compact set such that $\phi (y)\in K$ for all $%
y\in \mathbb{R}^{N}$. Viewing $\lambda \mapsto a(x,\omega ,\cdot ,\lambda )$
as a function defined on $K$, we have that this function belongs to $%
\mathcal{C}(K;(A)^{m})$ (use also hypothesis (v)), so that by the classical
Stone-Weierstrass theorem one has $a(x,\omega ,\cdot ,\phi )\in (A)^{m}$;
see either \cite[Proposition 1]{41} or \cite[Proposition 3.1]{40} for the
justification. As a result, we end up with the fact that the function $%
(x,\omega ,y)\mapsto a(x,\omega ,y,\psi (x,\omega ,y))$ belongs to $\mathcal{%
C}(\overline{Q};L^{\infty }(\Omega ;A))^{m}$.

We now use (iv) to get 
\begin{equation*}
\int_{Q\times \Omega }\left( a^{\varepsilon }(\cdot ,v_{\varepsilon
})-a^{\varepsilon }(\cdot ,\psi ^{\varepsilon })\right) \cdot
(v_{\varepsilon }-\psi ^{\varepsilon })dxd\mu \geq 0
\end{equation*}%
or equivalently, 
\begin{eqnarray*}
\int_{Q\times \Omega }a^{\varepsilon }(\cdot ,v_{\varepsilon })\cdot
v_{\varepsilon }dxd\mu &\geq &\int_{Q\times \Omega }a^{\varepsilon }(\cdot
,v_{\varepsilon })\cdot \psi _{\varepsilon }dxd\mu +\int_{Q\times \Omega
}a^{\varepsilon }(\cdot ,\psi _{\varepsilon })\cdot v_{\varepsilon }dxd\mu \\
&&-\int_{Q\times \Omega }a^{\varepsilon }(\cdot ,\psi _{\varepsilon })\cdot
\psi _{\varepsilon }dxd\mu .
\end{eqnarray*}%
Taking the $\lim \inf_{E^{\prime }\ni \varepsilon \rightarrow 0}$ of both
sides of the above inequality we get 
\begin{eqnarray}
\underset{E^{\prime }\ni \varepsilon \rightarrow 0}{\lim \inf }\int_{Q\times
\Omega }a^{\varepsilon }(\cdot ,v_{\varepsilon })\cdot v_{\varepsilon
}dxd\mu &\geq &\iint_{Q\times \Omega \times \Delta (A)}\widehat{z}_{0}\cdot 
\widehat{\psi }dxd\mu d\beta  \label{3.15} \\
&&+\iint_{Q\times \Omega \times \Delta (A)}\widehat{a}(\cdot ,\widehat{\psi }%
)\cdot \widehat{v}_{0}dxd\mu d\beta  \notag \\
&&-\iint_{Q\times \Omega \times \Delta (A)}\widehat{a}(\cdot ,\widehat{\psi }%
)\cdot \widehat{\psi }dxd\mu d\beta  \notag
\end{eqnarray}%
where: for the first integral on the right-hand side of (\ref{3.15}), we
have used the definition of the weak stochastic $\Sigma $-convergence for
the sequence $(a^{\varepsilon }(\cdot ,v_{\varepsilon }))_{\varepsilon \in
E} $, for the second integral, we have used the definition of the weak
stochastic $\Sigma $-convergence of $(v_{\varepsilon })_{\varepsilon }$
associated with Proposition \ref{p3.4} by taking $a(\cdot ,\psi )$ as a test
function, and finally for the last integral, we use the same argument as for
the preceding integral. Therefore, subtracting $\iint_{Q\times \Omega \times
\Delta (A)}\widehat{z}_{0}\cdot \widehat{v}_{0}dxd\mu d\beta $ from each
member of (\ref{3.15}), we end up with 
\begin{equation}
\begin{array}{l}
\underset{E^{\prime }\ni \varepsilon \rightarrow 0}{\lim \inf }\int_{Q\times
\Omega }a^{\varepsilon }(\cdot ,v_{\varepsilon })\cdot v_{\varepsilon
}dxd\mu -\iint_{Q\times \Omega \times \Delta (A)}\widehat{z}_{0}\cdot 
\widehat{v}_{0}dxd\mu d\beta \\ 
\ \ \ \ \ \ \ \ \ \ \ \ \ \geq -\iint_{Q\times \Omega \times \Delta
(A)}\left( \widehat{z}_{0}-\widehat{a}(\cdot ,\widehat{\psi })\right) \cdot (%
\widehat{v}_{0}-\widehat{\psi })dxd\mu d\beta%
\end{array}
\label{3.16}
\end{equation}%
for any $\psi \in \lbrack \mathcal{C}_{0}^{\infty }(Q)\otimes \mathcal{C}%
^{\infty }(\Omega )\otimes A]^{m}$. The right-hand side of (\ref{3.16}) is
of the form $g(x,\omega ,s,\widehat{\psi }(x,\omega ,s))$ and, due to the
fact that $\widehat{z}_{0}\in L^{p^{\prime }}(Q\times \Omega \times \Delta
(A))^{m}$, one easily deduces from assumption (iii) (in Theorem \ref{t5})
that $g(x,\omega ,s,\widehat{\psi })\in L^{1}(Q\times \Omega \times \Delta
(A))$ for any $\widehat{\psi }\in L^{p}(Q\times \Omega \times \Delta
(A))^{m} $, so that the operator $G$ defined here as in Lemma \ref{l6} (by
taking there $X=Q\times \Omega \times \Delta (A)$, $F_{1}=L^{p}(Q\times
\Omega \times \Delta (A))^{m}$, $F_{2}=L^{1}(Q\times \Omega \times \Delta
(A))$) maps $L^{p}(X;F_{1})$ into $L^{1}(X;F_{2})$. In view of Lemma \ref{l6}%
, the map $G$ is continuous under the norm topology. As a result, the
inequality (\ref{3.16}) holds for any $\widehat{\psi }\in L^{p}(Q\times
\Omega \times \Delta (A))^{m}$ (that is for any $\psi \in L^{p}(Q\times
\Omega ;\mathcal{B}_{A}^{p})^{m}$). Hence taking in (\ref{3.16}) $\psi
=v_{0} $ we get readily (\ref{3.14}).

For the last part of the theorem, assuming that (\ref{3.14}) is actually an
equality, we return to (\ref{3.16}) and take there $\psi =v_{0}+tw$, $w\in
L^{p}(Q\times \Omega ;\mathcal{B}_{A}^{p})^{m}$ being arbitrarily fixed and $%
t>0$. Then, 
\begin{equation*}
\iint_{Q\times \Omega \times \Delta (A)}\left( \widehat{z}_{0}-\widehat{a}%
(\cdot ,\widehat{v}_{0}+t\widehat{w})\right) \cdot \widehat{w}dxd\mu d\beta
\leq 0\ \ \forall w\in L^{p}(Q\times \Omega ;\mathcal{B}_{A}^{p})^{m}\text{.}
\end{equation*}%
Letting $t\rightarrow 0$, and next changing $w$ for $-w$, we end up with 
\begin{equation*}
\iint_{Q\times \Omega \times \Delta (A)}\left( \widehat{z}_{0}-\widehat{a}%
(\cdot ,\widehat{v}_{0})\right) \cdot \widehat{w}dxd\mu d\beta =0\ \ \forall
w\in L^{p}(Q\times \Omega ;\mathcal{B}_{A}^{p})^{m}\text{,}
\end{equation*}%
which implies $z_{0}=a(\cdot ,v_{0})$.
\end{proof}

As was said before the statement of Theorem \ref{t5}, if we take $a(x,\omega
,y,\lambda )=\left| \lambda \right| ^{p-2}\lambda $ and $m=1$, then we
arrive at the claimed conclusion by the above theorem.

Now we assume in the sequel that the $H$-supralgebra $A$ is translation
invariant and moreover each of its elements is uniformly continuous, that
is, $A$ is an algebra with mean value. The next result requires some
preliminaries. Let $a\in \mathbb{\mathbb{R}}^{N}$. Since $A$ is translation
invariant, the translation operator $\tau _{a}:A\rightarrow A$ extends by
continuity to a unique translation operator still denoted by $\tau
_{a}:B_{A}^{p}\rightarrow B_{A}^{p}$ ($1\leq p<\infty $). Indeed $\tau _{a}$
is bijective and $\left\| \tau _{a}u\right\| _{p}=\left\| u\right\| _{p}$
since $M(\left| \tau _{a}u\right| ^{p})=M(\tau _{a}\left| u\right|
^{p})=M(\left| u\right| ^{p})$ for all $u\in A$. Besides, as each element of 
$A$ is uniformly continuous, the group of unitary operators $\{\tau
_{a}:a\in \mathbb{\mathbb{R}}^{N}\}$ thus defined is strongly continuous,
i.e. $\tau _{a}u\rightarrow u$ in $B_{A}^{p}$ as $\left| a\right|
\rightarrow 0$ for all $u\in B_{A}^{p}$. Moreover 
\begin{equation}
M(\tau _{a}u)=M(u)\text{ for all }u\in B_{A}^{p}\text{ and any }a\in \mathbb{%
\mathbb{R}}^{N}\text{.}  \label{5.2'}
\end{equation}%
Arguing as above we see that the group $\{\tau _{a}\}_{a\in \mathbb{\mathbb{R%
}}^{N}}$ yields a family of mappings still denoted by $\{\tau _{a}\}_{a\in 
\mathbb{\mathbb{R}}^{N}}$ (each of them sending $L^{p}(\Omega ;B_{A}^{p})$
into itself) verifying 
\begin{equation*}
\tau _{a}u\left( \omega ,y\right) =\tau _{a}u(\omega ,\cdot )(y)=u\left(
\omega ,y+a\right) \text{ for a.e. }(\omega ,y)\in \Omega \times \mathbb{R}%
^{N}\text{ and for }u\in L^{p}(\Omega ;B_{A}^{p}).
\end{equation*}%
With this in mind, we begin with the following preliminary result.

\begin{lemma}
\label{l5.1}Assume the $H$-supralgebra $A$ is an algebra with mean value on $%
\mathbb{R}_{y}^{N}$, i.e., it is translation invariant and each of its
elements is uniformly continuous. Let $(u_{\varepsilon })_{\varepsilon \in
E} $ be a sequence in $L^{p}(Q\times \Omega )$ $(1<p<\infty )$ which weakly
stochastically $\Sigma $-converges towards $u_{0}\in L^{p}(Q\times \Omega ;%
\mathcal{B}_{A}^{p})$. Let the sequence $(v_{\varepsilon })_{\varepsilon \in
E}$ be defined by 
\begin{equation*}
v_{\varepsilon }(x,\omega )=\int_{B_{r}}u_{\varepsilon }(x+\varepsilon
_{2}\rho ,\omega )d\rho \text{\ \ }((x,\omega )\in Q\times \Omega ).
\end{equation*}%
Then, as $E\ni \varepsilon \rightarrow 0$, 
\begin{equation}
v_{\varepsilon }\rightarrow v_{0}\text{ stoch. in }L^{p}(Q\times \Omega )%
\text{-weak }\Sigma  \label{5.3'}
\end{equation}%
where $v_{0}$ is defined by $v_{0}(x,\omega ,y)=\int_{B_{r}}u_{0}(x,\omega
,y+\rho )d\rho $\ for $(x,\omega ,y)\in Q\times \Omega \times \mathbb{%
\mathbb{R}}^{N}$.
\end{lemma}

\begin{remark}
\label{r5.2}\emph{Assume Lemma \ref{l5.1} holds. Then as }$E\ni \varepsilon
\rightarrow 0$\emph{, }%
\begin{equation}
\frac{1}{\left\vert B_{\varepsilon _{2}r}\right\vert }\int_{B_{\varepsilon
_{2}r}}u_{\varepsilon }(x+y,\omega )dy\rightarrow \frac{1}{\left\vert
B_{r}\right\vert }v_{0}\text{\ \emph{stoch.} \emph{in} }L^{p}(Q\times \Omega
)\text{\emph{-weak }}\Sigma \text{\emph{.}}  \label{5.5'}
\end{equation}%
\emph{The above convergence result will be of particular interest in the
next result.}
\end{remark}

\begin{proof}[\textit{Proof of Lemma} \ref{l5.1}]
Let $\varphi \in \mathcal{C}_{0}^{\infty }(Q)$, $f\in \mathcal{C}^{\infty
}(\Omega )$ and $g\in A$. One has 
\begin{eqnarray*}
&&\int_{Q\times \Omega }\left( \int_{B_{r}}u_{\varepsilon }(x+\varepsilon
_{2}\rho ,\omega )d\rho \right) \varphi (x)f\left( T\left( \frac{x}{%
\varepsilon _{1}}\right) \omega \right) g\left( \frac{x}{\varepsilon _{2}}%
\right) dxd\mu \\
&=&\int_{B_{r}}\left( \int_{Q\times \Omega }u_{\varepsilon }(x+\varepsilon
_{2}\rho ,\omega )\varphi (x)f\left( T\left( \frac{x}{\varepsilon _{1}}%
\right) \omega \right) g\left( \frac{x}{\varepsilon _{2}}\right) dxd\mu
\right) d\rho .
\end{eqnarray*}%
In view of the Lebesgue dominated convergence theorem, (\ref{5.3'}) will be
checked as soon as we show that for each fixed $\rho \in \mathbb{\mathbb{R}}%
^{N}$, 
\begin{eqnarray*}
&&\int_{Q\times \Omega }u_{\varepsilon }(x+\varepsilon _{2}\rho ,\omega
)\varphi (x)f\left( T\left( \frac{x}{\varepsilon _{1}}\right) \omega \right)
g\left( \frac{x}{\varepsilon _{2}}\right) dxd\mu \\
&\rightarrow &\int_{Q\times \Omega }\int_{\Delta (A)}\widehat{\tau _{-\rho
}u_{0}}(x,\omega ,s)\varphi (x)f(\omega )\widehat{g}(s)d\beta dxd\mu \text{
when }E\ni \varepsilon \rightarrow 0\text{.}
\end{eqnarray*}%
First of all, let us beginning by noticing that since $\mathcal{G}_{1}$ is a
bounded linear operator of $\mathcal{B}_{A}^{1}$ into $L^{1}(\Delta (A))$,
we have 
\begin{equation*}
\mathcal{G}_{1}\left( \int_{Br}u_{0}(x,\omega ,\cdot +\rho )d\rho \right)
=\int_{B_{r}}\mathcal{G}_{1}(u_{0}(x,\omega ,\cdot +\rho ))d\rho
\end{equation*}%
where $u_{0}$ is as above. So let $a\in \mathbb{\mathbb{R}}^{N}$ and let $%
\varphi $, $f$ and $g$ be as above. One has 
\begin{equation*}
\begin{array}{l}
\int_{Q\times \Omega }u_{\varepsilon }(x-\varepsilon _{2}a,\omega )\varphi
(x)f\left( T\left( \frac{x}{\varepsilon _{1}}\right) \omega \right) g\left( 
\frac{x}{\varepsilon _{2}}\right) dxd\mu = \\ 
\;\;=\int_{\left( Q-\varepsilon _{2}a\right) \times \Omega }u_{\varepsilon
}(x,\omega )\varphi (x+\varepsilon _{2}a)f\left( T\left( \frac{x}{%
\varepsilon _{1}}+\frac{\varepsilon _{2}}{\varepsilon _{1}}a\right) \omega
\right) g\left( \frac{x}{\varepsilon _{2}}+a\right) dxd\mu \\ 
\;\;=\int_{Q\times \Omega }u_{\varepsilon }(x,\omega )\varphi (x+\varepsilon
_{2}a)f\left( T\left( \frac{x}{\varepsilon _{1}}+\frac{\varepsilon _{2}}{%
\varepsilon _{1}}a\right) \omega \right) g\left( \frac{x}{\varepsilon _{2}}%
+a\right) dxd\mu \\ 
\;\;\;\;\;-\int_{\left( Q\backslash (Q-\varepsilon _{2}a)\right) \times
\Omega }u_{\varepsilon }(x,\omega )\varphi (x+\varepsilon _{2}a)f\left(
T\left( \frac{x}{\varepsilon _{1}}+\frac{\varepsilon _{2}}{\varepsilon _{1}}%
a\right) \omega \right) g\left( \frac{x}{\varepsilon _{2}}+a\right) dxd\mu
\\ 
\;\;\;\;\;\;\;+\int_{\left( (Q-\varepsilon _{2}a)\backslash Q\right) \times
\Omega }u_{\varepsilon }(x,\omega )\varphi (x+\varepsilon _{2}a)f\left(
T\left( \frac{x}{\varepsilon _{1}}+\frac{\varepsilon _{2}}{\varepsilon _{1}}%
a\right) \omega \right) g\left( \frac{x}{\varepsilon _{2}}+a\right) dxd\mu
\\ 
=(I)-(II)+(III).%
\end{array}%
\end{equation*}%
As for $(I)$ we have 
\begin{eqnarray*}
(I) &=&\int_{Q\times \Omega }u_{\varepsilon }(x,\omega )\varphi (x)f\left(
T\left( \frac{x}{\varepsilon _{1}}+\frac{\varepsilon _{2}}{\varepsilon _{1}}%
a\right) \omega \right) (\tau _{-a}g)\left( \frac{x}{\varepsilon _{2}}%
\right) dxd\mu \\
&&+\int_{Q\times \Omega }u_{\varepsilon }(x,\omega )[\varphi (x+\varepsilon
_{2}a)-\varphi (x)]f\left( T\left( \frac{x}{\varepsilon _{1}}+\frac{%
\varepsilon _{2}}{\varepsilon _{1}}a\right) \omega \right) (\tau
_{-a}g)\left( \frac{x}{\varepsilon _{2}}\right) dxd\mu \\
\;\;\;\;\;\; &=&(I_{1})+(I_{2}).
\end{eqnarray*}%
But 
\begin{eqnarray*}
(I_{1}) &=&\int_{Q\times \Omega }u_{\varepsilon }(x,\omega )\varphi
(x)f\left( T\left( \frac{x}{\varepsilon _{1}}\right) \omega \right) (\tau
_{-a}g)\left( \frac{x}{\varepsilon _{2}}\right) dxd\mu \\
&&+\int_{Q\times \Omega }u_{\varepsilon }(x,\omega )\varphi (x)(\tau
_{-a}g)\left( \frac{x}{\varepsilon _{2}}\right) \left[ f\left( T\left( \frac{%
x}{\varepsilon _{1}}+\frac{\varepsilon _{2}}{\varepsilon _{1}}a\right)
\omega \right) -f\left( T\left( \frac{x}{\varepsilon _{1}}\right) \omega
\right) \right] dxd\mu \\
&=&(I_{1}^{\prime })+(I_{2}^{\prime }).
\end{eqnarray*}%
The $H$-supralgebra $A$ being translation invariant, we have $\tau _{-a}g\in
A$ and so, 
\begin{equation*}
(I_{1}^{\prime })\rightarrow \int_{Q\times \Omega }\int_{\Delta (A)}\widehat{%
u}_{0}(x,\omega ,s)\varphi (x)f(\omega )\widehat{\tau _{-a}g}(s)d\beta
dxd\mu \text{\ as }E\ni \varepsilon \rightarrow 0.
\end{equation*}%
But 
\begin{eqnarray*}
\int_{\Delta (A)}\widehat{u}_{0}(x,\omega ,s)\widehat{\tau _{-a}g}(s)d\beta
&=&M(u_{0}(x,\omega ,\cdot )(\tau _{-a}g)) \\
&=&M(\tau _{-a}[\tau _{a}u_{0}(x,\omega ,\cdot )g]) \\
&=&M(\tau _{a}u_{0}(x,\omega ,\cdot )g)\text{\ (see (\ref{5.2'}))} \\
&=&\int_{\Delta (A)}\widehat{\tau _{a}u_{0}}(x,\omega ,s)\widehat{g}%
(s)d\beta .
\end{eqnarray*}%
Note that here we have identified $u_{0}(x,\omega ,\cdot )\in \mathcal{B}%
_{A}^{p}$ with its representative still denoted by $u_{0}(x,\omega ,\cdot
)\in B_{A}^{p}$ so that $M_{1}(u_{0}(x,\omega ,\cdot ))=M(u_{0}(x,\omega
,\cdot ))$, $u_{0}(x,\omega ,\cdot )$ on the left-hand side of the above
equality being an equivalence class whereas $u_{0}(x,\omega ,\cdot )$ on the
right-hand side is one of its representative. For $(I_{2}^{\prime })$, we
have 
\begin{eqnarray*}
&&\left\vert (I_{2}^{\prime })\right\vert \\
&\leq &\left\Vert u_{\varepsilon }\right\Vert _{L^{p}\left( Q\times \Omega
\right) }\left\Vert \varphi \right\Vert _{\infty }\left\Vert g\right\Vert
_{\infty }\left( \int_{Q\times \Omega }\left\vert f\left( T\left( \frac{x}{%
\varepsilon _{1}}+\frac{\varepsilon _{2}}{\varepsilon _{1}}a\right) \omega
\right) -f\left( T\left( \frac{x}{\varepsilon _{1}}\right) \omega \right)
\right\vert ^{p^{\prime }}dxd\mu \right) ^{1/p^{\prime }}.
\end{eqnarray*}%
But 
\begin{eqnarray*}
&&\int_{Q\times \Omega }\left\vert f\left( T\left( \frac{x}{\varepsilon _{1}}%
+\frac{\varepsilon _{2}}{\varepsilon _{1}}a\right) \omega \right) -f\left(
T\left( \frac{x}{\varepsilon _{1}}\right) \omega \right) \right\vert
^{p^{\prime }}dxd\mu \\
&=&\int_{Q\times \Omega }\left\vert \left( U\left( \frac{x}{\varepsilon _{1}}%
+\frac{\varepsilon _{2}}{\varepsilon _{1}}a\right) f\right) (\omega )-\left(
U\left( \frac{x}{\varepsilon _{1}}\right) f\right) (\omega )\right\vert
^{p^{\prime }}dxd\mu \text{.}
\end{eqnarray*}%
Since the group $U(x)$ is strongly continuous in $L^{p^{\prime }}\left(
\Omega \right) $ (see Subsection 2.1) we get immediately (using the Lebesgue
dominated convergence theorem) that 
\begin{equation*}
\int_{Q\times \Omega }\left\vert f\left( T\left( \frac{x}{\varepsilon _{1}}+%
\frac{\varepsilon _{2}}{\varepsilon _{1}}a\right) \omega \right) -f\left(
T\left( \frac{x}{\varepsilon _{1}}\right) \omega \right) \right\vert
^{p^{\prime }}dxd\mu \rightarrow 0\text{ as }\varepsilon \rightarrow 0.
\end{equation*}%
Thus $(I_{2}^{\prime })\rightarrow 0$ as $E\ni \varepsilon \rightarrow 0$.
Finally since the sequence $(u_{\varepsilon })_{\varepsilon \in E}$ is
bounded in $L^{p}(Q\times \Omega )$ and as $p>1$, this sequence is uniformly
integrable in $L^{1}(Q\times \Omega )$, so that from the inequality 
\begin{eqnarray*}
&&\int_{\left( (Q-\varepsilon _{2}a)\Delta Q\right) \times \Omega
}\left\vert u_{\varepsilon }(x,\omega )\right\vert \left\vert \varphi
(x+\varepsilon _{2}a)\right\vert \left\vert f\left( T\left( \frac{x}{%
\varepsilon _{1}}+\frac{\varepsilon _{2}}{\varepsilon _{1}}a\right) \omega
\right) \right\vert \left\vert g\left( \frac{x}{\varepsilon _{2}}+a\right)
\right\vert dxd\mu \\
&\leq &\left\Vert \varphi \right\Vert _{\infty }\left\Vert f\right\Vert
_{\infty }\left\Vert g\right\Vert _{\infty }\int_{\left( (Q-\varepsilon
_{2}a)\Delta Q\right) \times \Omega }\left\vert u_{\varepsilon }(x,\omega
)\right\vert dxd\mu ,
\end{eqnarray*}%
we see that $(II)$ and $(III)$ go towards $0$ as $E\ni \varepsilon
\rightarrow 0$; here the symbol $\Delta $ between the sets $(Q-\varepsilon
_{2}a)$ and $Q$ denotes the symmetric difference between these two sets.
Hence the lemma.
\end{proof}

We are now able to state and prove the most important compactness result of
the paper. It will be of capital interest in the next sections.

\begin{theorem}
\label{t3.3}Let $1<p<\infty $. Let $X$ be a normed closed convex subset of $%
W^{1,p}\left( Q\right) $, $Q$ being an open subset of $\mathbb{\mathbb{R}}%
^{N}$. Let $A$ be an ergodic supralgebra on $\mathbb{\mathbb{R}}_{y}^{N}$.
Assume $(u_{\varepsilon })_{\varepsilon \in E}$ is a sequence in $%
L^{p}(Q\times \Omega )$ such that:

\begin{itemize}
\item[(i)] $u_{\varepsilon }(\cdot ,\omega )\in X$ for all $\varepsilon \in
E $ and for $\mu $-a.e. $\omega \in \Omega $;

\item[(ii)] $(u_{\varepsilon })_{\varepsilon \in E}$ is bounded in $%
L^{p}(\Omega ;W^{1,p}(Q))$.
\end{itemize}

Then there exist $u_{0}\in W^{1,p}(Q;I_{nv}^{p}(\Omega ))$, $u_{1}\in
L^{p}(Q;\mathcal{W}^{1,p}(\Omega ))$, $u_{2}\in L^{p}(Q\times \Omega ;%
\mathcal{B}_{\#A}^{1,p})$ and a subsequence $E^{\prime }$ from $E$ such that

\begin{itemize}
\item[(iii)] $u_{0}(\cdot ,\omega )\in X$ for $\mu $-a.e. $\omega \in \Omega 
$ and $Du_{0}(x,\cdot )\in (I_{nv}^{p}(\Omega ))^{N}$ for a.e. $x\in Q$ and,
as $E^{\prime }\ni \varepsilon \rightarrow 0$,

\item[(iv)] $u_{\varepsilon }\rightarrow u_{0}$ stoch. in $L^{p}(Q\times
\Omega )$-weak;

\item[(v)] $Du_{\varepsilon }\rightarrow Du_{0}+\overline{D}_{\omega }u_{1}+%
\overline{D}_{y}u_{2}$ stoch. in $L^{p}(Q\times \Omega )^{N}$-weak $\Sigma $.
\end{itemize}
\end{theorem}

\begin{proof}
By Theorem \ref{t3.1}, there exist a subsequence $E^{\prime }$ from $E$, a
function $u_{0}\in L^{p}(Q\times \Omega ;\mathcal{B}_{A}^{p})$ and a vector
function $\mathbf{v}=(v_{i})_{1\leq i\leq N}\in L^{p}(Q\times \Omega ;%
\mathcal{B}_{A}^{p})^{N}$ such that, as $E^{\prime }\ni \varepsilon
\rightarrow 0$, we have $u_{\varepsilon }\rightarrow u_{0}$ stoch. in $%
L^{p}(Q\times \Omega )$-weak $\Sigma $ and $Du_{\varepsilon }\rightarrow 
\mathbf{v}$ stoch. in $L^{p}(Q\times \Omega )^{N}$-weak $\Sigma $.

\noindent At this level, the proof consists of three parts. We must check
that:\ Part (\textbf{I}) ($a$) $u_{0}$ does not depend upon $y$, that is $%
\overline{D}_{y}u_{0}=0$, ($b$) $u_{0}(x,\cdot )\in I_{nv}^{p}(\Omega )$,
that is $D_{\omega }u_{0}(x,\cdot )=0$ or equivalently $\int_{\Omega
}u_{0}(x,\cdot )D_{i,p}\varphi d\mu =0\;\forall \varphi \in \mathcal{C}%
^{\infty }(\Omega )$ and ($c$) $u_{0}\in W^{1,p}(Q;I_{nv}^{p}(\Omega ))$;
Part (\textbf{II}) $u_{0}(\cdot ,\omega )\in X$ for $\mu $-a.e. $\omega \in
\Omega $ and $Du_{0}(x,\cdot )\in I_{nv}^{p}(\Omega )^{N}$ for a.e. $x\in Q$%
; Part (\textbf{III}) There exist two functions $u_{1}\in L^{p}(Q;\mathcal{W}%
^{1,p}(\Omega ))$ and $u_{2}\in L^{p}(Q\times \Omega ;\mathcal{B}%
_{\#A}^{1,p})$ such that $\mathbf{v}=Du_{0}+\overline{D}_{\omega }u_{1}+%
\overline{D}_{y}u_{2}$.

Let us first prove (\textbf{I}). ($a$) Let $\Phi _{\varepsilon }(x,\omega
)=\varepsilon _{2}\varphi (x)f(T(x/\varepsilon _{1})\omega )g(x/\varepsilon
_{2})$ for $(x,\omega )\in Q\times \Omega $, where $\varphi \in \mathcal{C}%
_{0}^{\infty }(Q)$, $f\in \mathcal{C}^{\infty }(\Omega )$ and $g\in
A^{\infty }$. Then 
\begin{eqnarray*}
\int_{Q\times \Omega }\frac{\partial u_{\varepsilon }}{\partial x_{i}}\Phi
_{\varepsilon }dxd\mu &=&-\int_{Q\times \Omega }\varepsilon
_{2}u_{\varepsilon }f^{\varepsilon }g^{\varepsilon }\frac{\partial \varphi }{%
\partial x_{i}}dxd\mu -\int_{Q\times \Omega }u_{\varepsilon }\varphi
f^{\varepsilon }\left( D_{y_{i}}g\right) ^{\varepsilon }dxd\mu \\
&&-\int_{Q\times \Omega }\frac{\varepsilon _{2}}{\varepsilon _{1}}%
u_{\varepsilon }\varphi g^{\varepsilon }(D_{i,\omega }f)^{\varepsilon }dxd\mu
\end{eqnarray*}%
where $D_{y_{i}}g=\partial g/\partial y_{i}$. Letting $E^{\prime }\ni
\varepsilon \rightarrow 0$ we get 
\begin{equation*}
\iint_{Q\times \Omega \times \Delta (A)}\widehat{u}_{0}\varphi \widehat{%
D_{y_{i}}g}fdxd\mu d\beta =0\text{,}
\end{equation*}%
hence $\int_{\Delta (A)}\widehat{u}_{0}(x,\omega ,\cdot )\widehat{D_{y_{i}}g}%
d\beta =0$\ for all $g\in A^{\infty }$ and all $1\leq i\leq N$, which means
that $u_{0}$ does not depend on $y$ since the $H$-supralgebra $A$ is ergodic.

($b$) Let $\Phi _{\varepsilon }(x,\omega )=\varepsilon _{1}\varphi
(x)f(T(x/\varepsilon _{1})\omega )$ for $(x,\omega )\in Q\times \Omega $
where $\varphi \in \mathcal{C}_{0}^{\infty }(Q)$ and $f\in \mathcal{C}%
^{\infty }(\Omega )$. Then proceeding as above we get $\int_{\Omega
}u_{0}(x,\cdot )D_{i,\omega }fd\mu =0$ for all $1\leq i\leq N$ and $f\in 
\mathcal{C}^{\infty }(\Omega )$, which is equivalent to say that $%
u_{0}(x,\cdot )\in I_{nv}^{p}(\Omega )$ for a.e. $x\in Q$ (use property (P2)
in Subsection 2.1).

($c$) Hypothesis (ii) implies that the sequence $(u_{\varepsilon
})_{\varepsilon \in E^{\prime }}$ is bounded in $W^{1,p}(Q;L^{p}(\Omega ))$,
which yields the existence of a subsequence of $E^{\prime }$ not relabeled
and of a function $u\in W^{1,p}(Q;L^{p}(\Omega ))$ such that $u_{\varepsilon
}\rightarrow u$ in $W^{1,p}(Q;L^{p}(\Omega ))$-weak as $E^{\prime }\ni
\varepsilon \rightarrow 0$. In particular $\int_{\Omega }u_{\varepsilon
}(\cdot ,\omega )\psi (\omega )d\mu \rightarrow \int_{\Omega }u(\cdot
,\omega )\psi (\omega )d\mu $ in $L^{1}(Q)$-weak for all $\psi \in
I_{nv}^{p^{\prime }}(\Omega )$. Therefore using \cite{10} (see in particular
Lemma 3.6 therein) we get at once $u_{0}\in W^{1,p}(Q;L^{p}(\Omega ))$, so
that $u_{0}\in W^{1,p}(Q;I_{nv}^{p}(\Omega ))$.

As for (\textbf{II}), repeating the proof of [parts (iii) and (vi) of ] %
\cite[Theorem 3.7 (b)]{10} we are immediately led to (\textbf{II}). It
remains to check (\textbf{III}) here above. We begin by deriving the
existence of $u_{2}\in L^{p}(Q\times \Omega ;\mathcal{B}_{\#A}^{1,p})$. For
that purpose, let $r>0$ be freely fixed. Let $B_{\varepsilon _{2}r}$ denote
the open ball in $\mathbb{\mathbb{R}}^{N}$ centered at the origin and of
radius $\varepsilon _{2}r$. By the equalities 
\begin{eqnarray*}
&&\frac{1}{\varepsilon _{2}}\left( u_{\varepsilon }(x,\omega )-\frac{1}{%
\left\vert B_{\varepsilon _{2}r}\right\vert }\int_{B_{\varepsilon
_{2}r}}u_{\varepsilon }(x+\rho ,\omega )d\rho \right) \\
&=&\frac{1}{\varepsilon _{2}}\frac{1}{\left\vert B_{\varepsilon
_{2}r}\right\vert }\int_{B_{\varepsilon _{2}r}}\left( u_{\varepsilon
}(x,\omega )-u_{\varepsilon }(x+\rho ,\omega )\right) d\rho \\
&=&\frac{1}{\varepsilon _{2}}\frac{1}{\left\vert B_{r}\right\vert }%
\int_{B_{r}}\left( u_{\varepsilon }(x,\omega )-u_{\varepsilon
}(x+\varepsilon _{2}\rho ,\omega )\right) d\rho \\
&=&-\frac{1}{\left\vert B_{r}\right\vert }\int_{B_{r}}d\rho
\int_{0}^{1}Du_{\varepsilon }(x+t\varepsilon _{2}\rho ,\omega )\cdot \rho dt
\end{eqnarray*}%
where the dot denotes the usual Euclidean inner product in $\mathbb{\mathbb{R%
}}^{N}$, we deduce from the boundedness of $(u_{\varepsilon })_{\varepsilon
\in E^{\prime }}$ in $L^{p}(\Omega ;W^{1,p}(Q))$ that the sequence $%
(z_{\varepsilon }^{r})_{\varepsilon \in E^{\prime }}$ defined by 
\begin{equation*}
z_{\varepsilon }^{r}(x,\omega )=\frac{1}{\varepsilon _{2}}\left(
u_{\varepsilon }(x,\omega )-\frac{1}{\left\vert B_{\varepsilon
_{2}r}\right\vert }\int_{B_{\varepsilon _{2}r}}u_{\varepsilon }(x+\rho
,\omega )d\rho \right) \;((x,\omega )\in Q\times \Omega ,\varepsilon \in
E^{\prime })
\end{equation*}%
is bounded in $L^{p}(Q\times \Omega )$. It is important to note that in
general the function $z_{\varepsilon }^{r}$ is well defined since $%
u_{\varepsilon }$ and $Du_{\varepsilon }$ can be naturally extended off $Q$
as elements of $L^{p}(\Omega ;L_{\text{loc}}^{p}(\mathbb{R}^{N}))$ and $%
L^{p}(\Omega ;L_{\text{loc}}^{p}(\mathbb{R}^{N})^{N})$, respectively. Once
more, by virtue of Theorem \ref{t3.1} we find that there exist a subsequence
from $E^{\prime }$ (not relabeled) and a function $z_{r}$ in $L^{p}(Q\times
\Omega ;\mathcal{B}_{A}^{p})$ such that, as $E^{\prime }\ni \varepsilon
\rightarrow 0$ 
\begin{equation}
z_{\varepsilon }^{r}\rightarrow z_{r}\text{\ stoch. in }L^{p}(Q\times \Omega
)\text{-weak }\Sigma \text{.}  \label{5.8'}
\end{equation}%
As $(z_{\varepsilon }^{r})_{\varepsilon \in E^{\prime }}$ is bounded in $%
L^{p}(Q\times \Omega )$ we have (since $\varepsilon _{2}$, $\varepsilon
_{2}/\varepsilon _{1}\rightarrow 0$ as $E^{\prime }\ni \varepsilon
\rightarrow 0$) that 
\begin{equation}
\varepsilon _{2}z_{\varepsilon }^{r}\rightarrow 0\text{\ in }L^{p}(Q\times
\Omega )\text{ and }\frac{\varepsilon _{2}}{\varepsilon _{1}}z_{\varepsilon
}^{r}\rightarrow 0\text{ in }L^{p}(Q\times \Omega )\text{ as }E^{\prime }\ni
\varepsilon \rightarrow 0.  \label{5.9'}
\end{equation}%
Now, for $\varphi \in \mathcal{C}_{0}^{\infty }(Q)$, $f\in \mathcal{C}%
^{\infty }(\Omega )$ and $g\in A^{\infty }$ we have 
\begin{equation}
\begin{array}{l}
\int_{Q\times \Omega }\left( \frac{\partial u_{\varepsilon }}{\partial x_{i}}%
(x,\omega )-\frac{1}{\left\vert B_{\varepsilon _{2}r}\right\vert }%
\int_{B_{\varepsilon _{2}r}}\frac{\partial u_{\varepsilon }}{\partial x_{i}}%
(x+\rho ,\omega )d\rho \right) \varphi (x)f\left( T\left( \frac{x}{%
\varepsilon _{1}}\right) \omega \right) g(\frac{x}{\varepsilon _{2}})dxd\mu
\\ 
\;\;=-\int_{Q\times \Omega }\varepsilon _{2}z_{\varepsilon }^{r}(x,\omega
)f\left( T\left( \frac{x}{\varepsilon _{1}}\right) \omega \right) g(\frac{x}{%
\varepsilon _{2}})\frac{\partial \varphi }{\partial x_{i}}(x)dxd\mu \\ 
-\int_{Q\times \Omega }\frac{\varepsilon _{2}}{\varepsilon _{1}}%
z_{\varepsilon }^{r}(x,\omega )\varphi (x)g(\frac{x}{\varepsilon _{2}}%
)\left( D_{i,\omega }f\right) (T\left( \frac{x}{\varepsilon _{1}}\right)
\omega )dxd\mu \\ 
\;\;\;\;\;-\int_{Q\times \Omega }z_{\varepsilon }^{r}(x,\omega )\varphi
(x)f\left( T\left( \frac{x}{\varepsilon _{1}}\right) \omega \right) \frac{%
\partial g}{\partial y_{i}}\left( \frac{x}{\varepsilon _{2}}\right) dxd\mu .%
\end{array}
\label{5.10'}
\end{equation}%
Passing to the limit in (\ref{5.10'}) (as $E^{\prime }\ni \varepsilon
\rightarrow 0$) using conjointly (\ref{5.8'}), (\ref{5.9'}) and Remark \ref%
{r5.2} (see (\ref{5.5'}) therein) one gets 
\begin{equation*}
\begin{array}{l}
\iint_{Q\times \Omega \times \Delta (A)}\mathcal{G}_{1}\left( v_{i}(x,\omega
,\cdot )-\frac{1}{\left\vert B_{r}\right\vert }\int_{B_{r}}v_{i}(x,\omega
,\cdot +\rho )d\rho \right) (s)\varphi (x)f(\omega )\widehat{g}(s)dxd\mu
d\beta \\ 
\;\;\;=-\iint_{Q\times \Omega \times \Delta (A)}\widehat{z}_{r}(x,\omega
,s)\varphi (x)f(\omega )\partial _{i}\widehat{g}(s)dxd\mu d\beta ,%
\end{array}%
\end{equation*}%
the derivative $\partial _{i}$ in front of $\widehat{g}$ being the partial
derivative of index $i$ with respect to $\Delta (A)$ defined in the
preceding section as $\partial _{i}\widehat{g}=\mathcal{G}(\partial
g/\partial y_{i})$ (see also (\ref{2.7'}) therein). Therefore, because of
the arbitrariness of $\varphi $, $f$ and $g$, we are led to 
\begin{equation*}
\begin{array}{l}
\partial _{i}\widehat{z}_{r}(x,\omega ,\cdot )=\mathcal{G}_{1}\left(
v_{i}(x,\omega ,\cdot )-\frac{1}{\left\vert B_{r}\right\vert }%
\int_{B_{r}}v_{i}(x,\omega ,\cdot +\rho )d\rho \right) \text{\ a.e. in }%
\Delta (A)\text{ } \\ 
\text{for }(x,\omega )\in Q\times \Omega .%
\end{array}%
\end{equation*}%
But $\partial _{i}\widehat{z}_{r}(x,\omega ,\cdot )=\partial _{i}\mathcal{G}%
_{1}(z_{r}(x,\omega ,\cdot ))=\mathcal{G}_{1}\left( \overline{\partial }%
z_{r}/\partial y_{i}(x,\omega ,\cdot )\right) $, hence, for $1\leq i\leq N$, 
\begin{equation*}
\frac{\overline{\partial }z_{r}}{\partial y_{i}}(x,\omega ,\cdot
)=v_{i}(x,\omega ,\cdot )-\frac{1}{\left\vert B_{r}\right\vert }%
\int_{B_{r}}v_{i}(x,\omega ,\cdot +\rho )d\rho \text{\ a.e. in }\mathbb{R}%
_{y}^{N}\text{\ for }(x,\omega )\in Q\times \Omega
\end{equation*}%
(recall that $\mathcal{G}_{1}$ is an isomorphism of $\mathcal{B}_{A}^{1}$
onto $L^{1}(\Delta (A))$ which carries over $\mathcal{B}_{A}^{p}$ onto $%
L^{p}(\Delta (A))$ isomorphically and isometrically). Set $f_{r}(x,\omega
,y)=z_{r}(x,\omega ,y)-M_{y}(z_{r}(x,\omega ,\cdot ))$ where here, $%
z_{r}(x,\omega ,\cdot )\in \mathcal{B}_{A}^{p}$ is viewed as its
representative in $B_{A}^{p}$ and $M_{y}=M$ standing here for the mean value
with respect to $y$ defined as in the preceding section (see in particular
property (\textbf{2}) and equality (\ref{0.1}) in Subsection 2.3). Then $%
M_{y}(f_{r})=0$ and moreover $\overline{D}_{y}f_{r}=\overline{D}_{y}z_{r}$
so that $f_{r}\in L^{p}(Q\times \Omega ;\mathcal{B}_{A}^{p})$ with $%
\overline{\partial }f_{r}/\partial y_{i}\in L^{p}(Q\times \Omega ;\mathcal{B}%
_{A}^{p})$, that is, 
\begin{equation*}
f_{r}\in L^{p}(Q\times \Omega ;\mathcal{B}_{A}^{1,p}/\mathbb{C})\text{.}%
\;\;\;\;\;\;\;\;\;\;\;\;\;\;\;\;\;\;\;\;
\end{equation*}%
So let $g_{r}=J_{1}\circ f_{r}$, where $J_{1}$ denotes the canonical mapping
of $\mathcal{B}_{A}^{1,p}/\mathbb{C}$ into its separated completion $%
\mathcal{B}_{\#A}^{1,p}$. Then $g_{r}\in L^{p}(Q\times \Omega ;\mathcal{B}%
_{\#A}^{1,p})$ and moreover 
\begin{equation*}
\frac{\overline{\partial }g_{r}}{\partial y_{i}}(x,\omega ,\cdot
)=v_{i}(x,\omega ,\cdot )-\frac{1}{\left\vert B_{r}\right\vert }%
\int_{B_{r}}v_{i}(x,\omega ,\cdot +\rho )d\rho \;\ \ (1\leq i\leq N)
\end{equation*}%
since $\frac{\overline{\partial }g_{r}}{\partial y_{i}}(x,\omega ,\cdot )=%
\frac{\overline{\partial }f_{r}}{\partial y_{i}}(x,\omega ,\cdot )=\frac{%
\overline{\partial }z_{r}}{\partial y_{i}}(x,\omega ,\cdot )$. Now, we also
view $v_{i}(x,\omega ,\cdot )$ as its representative in $B_{A}^{p}$. Taking
this into account, we have 
\begin{equation}
\begin{array}{l}
\left\Vert g_{r}(x,\omega ,\cdot )-g_{r^{\prime }}(x,\omega ,\cdot
)\right\Vert _{\mathcal{B}_{\#A}^{1,p}} \\ 
\leq \left\Vert \overline{D}_{y}g_{r}(x,\omega ,\cdot )-\mathbf{v}(x,\omega
,\cdot )+M_{y}(\mathbf{v}(x,\omega ,\cdot ))\right\Vert _{p} \\ 
+\left\Vert \overline{D}_{y}g_{r^{\prime }}(x,\omega ,\cdot )-\mathbf{v}%
(x,\omega ,\cdot )+M_{y}(\mathbf{v}(x,\omega ,\cdot ))\right\Vert _{p}.%
\end{array}
\label{5.11'}
\end{equation}%
But 
\begin{equation*}
\begin{array}{l}
\left\Vert \overline{D}_{y}g_{r}(x,\omega ,\cdot )-\mathbf{v}(x,\omega
,\cdot )+M_{y}(\mathbf{v}(x,\omega ,\cdot ))\right\Vert _{p} \\ 
\;\;\;=\left\Vert \frac{1}{\left\vert B_{r}\right\vert }\int_{B_{r}}\mathbf{v%
}(x,\omega ,\cdot +\rho )d\rho -M_{y}(\mathbf{v}(x,\omega ,\cdot
))\right\Vert _{p}.%
\end{array}%
\end{equation*}%
Therefore, since the algebra $A$ is ergodic, the right-hand side (and hence
the left-hand side) of (\ref{5.11'}) goes to zero when $r,r^{\prime
}\rightarrow +\infty $. Thus, the sequence $(g_{r}(x,\omega ,\cdot ))_{r>0}$
is a Cauchy sequence in the Banach space $\mathcal{B}_{\#A}^{1,p}$, whence
the existence of a unique $u_{2}(x,\omega ,\cdot )\in \mathcal{B}%
_{\#A}^{1,p} $ such that 
\begin{equation*}
g_{r}(x,\omega ,\cdot )\rightarrow u_{2}(x,\omega ,\cdot )\text{\ in }%
\mathcal{B}_{\#A}^{1,p}\text{\ as }r\rightarrow +\infty ,
\end{equation*}%
that is 
\begin{equation*}
\overline{D}_{y}g_{r}(x,\omega ,\cdot )\rightarrow \overline{D}%
_{y}u_{2}(x,\omega ,\cdot )\text{\ in }(\mathcal{B}_{A}^{p})^{N}\text{\ as }%
r\rightarrow +\infty .
\end{equation*}%
Once again the ergodicity of $A$ and the uniqueness of the limit leads at
once to 
\begin{equation*}
\overline{D}_{y}u_{2}(x,\omega ,\cdot )=\mathbf{v}(x,\omega ,\cdot )-M_{y}(%
\mathbf{v}(x,\omega ,\cdot ))\text{\ a.e. in }\mathbb{R}^{N}\text{\ and for
a.e. }(x,\omega )\in Q\times \Omega \text{.}
\end{equation*}%
We deduce the existence of a function $u_{2}:Q\times \Omega \rightarrow 
\mathcal{B}_{A}^{p}$, $(x,\omega )\mapsto u_{2}(x,\omega ,\cdot )$, lying in 
$L^{p}(Q\times \Omega ;\mathcal{B}_{A}^{p})$ such that 
\begin{equation}
\mathbf{v}-M(\mathbf{v})=\overline{D}_{y}u_{2}.\;\;\;\;  \label{Equ1}
\end{equation}%
Let us finally derive the existence of $u_{1}$. Let $\Phi _{\varepsilon
}(x,\omega )=\varphi (x)\Psi (T(x/\varepsilon _{1})\omega )$ ($(x,\omega
)\in Q\times \Omega )$) with $\varphi \in \mathcal{C}_{0}^{\infty }(Q)$ and $%
\Psi =(\psi _{j})_{1\leq j\leq N}\in \mathcal{V}_{\func{div}}$ (i.e. $\func{%
div}_{\omega ,p^{\prime }}\Psi =0$). Clearly 
\begin{equation*}
\sum_{j=1}^{N}\int_{Q\times \Omega }\frac{\partial u_{\varepsilon }}{%
\partial x_{j}}\varphi \psi _{j}^{\varepsilon }dxd\mu
=-\sum_{j=1}^{N}\int_{Q\times \Omega }u_{\varepsilon }\psi _{j}^{\varepsilon
}\frac{\partial \varphi }{\partial x_{j}}dxd\mu
\end{equation*}%
where $\psi _{j}^{\varepsilon }(x,\omega )=\psi _{j}(T(x/\varepsilon
_{1})\omega )$. Passing to the limit when $E^{\prime }\ni \varepsilon
\rightarrow 0$ yields 
\begin{equation*}
\sum_{j=1}^{N}\iint_{Q\times \Omega \times \Delta (A)}\widehat{v}_{j}\varphi
\psi _{j}dxd\mu d\beta =-\sum_{j=1}^{N}\iint_{Q\times \Omega \times \Delta
(A)}u_{0}\psi _{j}\frac{\partial \varphi }{\partial x_{j}}dxd\mu d\beta 
\text{,}
\end{equation*}%
or equivalently, 
\begin{equation*}
\iint_{Q\times \Omega \times \Delta (A)}\left( \widehat{\mathbf{v}}%
-Du_{0}\right) \cdot \Psi \varphi dxd\mu d\beta =0\text{,}
\end{equation*}%
and so, as $\varphi $ is arbitrarily fixed in $\mathcal{C}_{0}^{\infty }(Q)$%
, 
\begin{equation*}
\iint_{\Omega \times \Delta (A)}\left( \widehat{\mathbf{v}}(x,\omega
,s)-Du_{0}(x,\omega )\right) \cdot \Psi (\omega )d\mu d\beta =0\;\;\forall
\Psi \in \mathcal{V}_{\func{div}}.
\end{equation*}%
This is also equivalent to 
\begin{equation*}
\int_{\Omega }\left( M(\mathbf{v})-Du_{0}\right) \cdot \Psi d\mu =0\;\;\text{%
for all }\Psi \in \mathcal{V}_{\func{div}}\text{.}
\end{equation*}%
Therefore, the Proposition \ref{p2.1} provides us with a function $%
u_{1}(x,\cdot )\in \mathcal{W}^{1,p}(\Omega )$ such that 
\begin{equation}
M(\mathbf{v})-Du_{0}=\overline{D}_{\omega }u_{1}(x,\cdot ).  \label{Equ2}
\end{equation}%
Putting (\ref{Equ1}) and (\ref{Equ2}) together leads at once at $\mathbf{v}%
=Du_{0}+\overline{D}_{\omega }u_{1}+\overline{D}_{y}u_{2}$, where the
function $u_{1}:x\mapsto u_{1}(x,\cdot )$ lies in $L^{p}(Q;\mathcal{W}%
^{1,p}(\Omega ))$. This completes the proof.
\end{proof}

\begin{remark}
\label{r3.2}\emph{The preceding theorem generalizes its homologue (see
Theorem 3.5 in \cite{CMP}) as follows: In Theorem \ref{t3.3} above, take }$%
\Omega =\Delta (A_{z})$\emph{\ where }$A_{z}$\emph{\ is any }$H$\emph{%
-supralgebra on }$\mathbb{R}^{N}$\emph{\ which is translation invariant and
whose elements are uniformly continuous. Then thanks to Theorem \ref{t2.2},
a dynamical system can be constructed on }$\Delta (A_{z})$\emph{\ such that
the corresponding invariant probability measure is the }$M$\emph{-measure }$%
\beta _{z}$\emph{\ associated to }$A_{z}$\emph{. Therefore by the equality (%
\ref{2.12}), our claim is justified since in \cite[Theorem 3.5]{CMP}, both
the algebras }$A$\emph{\ and }$A_{z}$\emph{\ are assumed to be ergodic while
here, the algebra }$A_{z}$\emph{\ is not assumed to be ergodic. We will see
in the next section how the above result is used, and how it generalizes the
one in \cite{CMP} as claimed.}
\end{remark}

\section{Application to the homogenization of a linear partial differential
equation}

We need to show how the preceding result arises in the homogenization of
partial differential equations. To illustrate this we begin by focusing our
attention on the rather simple case of an elliptic linear differential
operator of order two, in divergence form, namely, we consider the following
boundary-value problem 
\begin{equation}
\begin{array}{l}
-\overset{N}{\underset{i,j=1}{\sum }}\frac{\partial }{\partial x_{i}}\left(
a_{ij}(x,T(x/\varepsilon _{1})\omega ,x/\varepsilon _{2})\frac{\partial
u_{\varepsilon }}{\partial x_{j}}\right) =f\text{ in }Q \\ 
\ \ \ \ \ \ \ \ \ \ \ \ \ \ \ \ \ \ \ \ \ \ \ \ \ \ \ \ \ \ \ \ \ \ \ \ \ \
\ \ \ \ \ \ \ \ \ u_{\varepsilon }=0\text{ on }\partial Q%
\end{array}
\label{4.1}
\end{equation}%
where $Q$ is a bounded open subset in $\mathbb{R}^{N}$, $f\in L^{\infty
}(\Omega ;H^{-1}(Q))=L^{\infty }(\Omega ;W^{-1,2}(Q))$, $a_{ij}\in \mathcal{C%
}(\overline{Q};L^{\infty }(\Omega ;\mathcal{B}(\mathbb{R}_{y}^{N})))$, $%
a_{ij}=\overline{a}_{ji}$ (the complex conjugate of $a_{ji}$), and $%
(a_{ij})_{1\leq i,j\leq N}$ satisfies the following ellipticity condition:
there exists a constant $\alpha >0$ such that $\sum_{i,j=1}^{N}a_{ij}(x,%
\omega ,y)\lambda _{i}\overline{\lambda }_{j}\geq \alpha \left\vert \lambda
\right\vert ^{2}$ for all $(x,y)\in \overline{Q}\times \mathbb{R}^{N}$, for $%
d\mu $-almost all $\omega \in \Omega $ and for all $\lambda \in \mathbb{C}%
^{N}$. It is a well-known fact that for each $\varepsilon >0$ (\ref{4.1})
uniquely determines $u_{\varepsilon }=u_{\varepsilon }(x,\omega )\in
H_{0}^{1}(Q;L^{2}(\Omega ))$ in such a way that we have in hands a
generalized sequence $(u_{\varepsilon })_{\varepsilon >0}$. The fundamental
problem in homogenization theory is the study of the asymptotic behavior of
such a sequence under a suitable assumption made on the coefficients $a_{ij}$
of the operator in (\ref{4.1}). Here, as we will see in the sequel, it will
be sufficient to make this assumption with respect to the variable $y\in 
\mathbb{R}^{N}$. Prior to this, it is worth to recall the following facts:
firstly, in the case when the functions $a_{ij}$ do not depend on the
variable $y$, the homogenization of (\ref{4.1}) has been conducted in \cite%
{10}; secondly, in the case when the coefficients $a_{ij}$ depend only on
the variables $x,y$ (i.e. the functions $a_{ij}(x,\cdot ,y)$ are constants),
it is commonly known that under the periodicity assumption on the functions $%
a_{ij}$ (with respect to $y$), the homogenization problem for (\ref{4.1})
has already been solved by many authors and the results are available in the
literature. In the same case, it is also known that in the general framework
of deterministic homogenization theory the same results are available in the
ergodic environment; see e.g. \cite{26}. However, in contrast with the
ergodic setting, no result is available in the non-ergodic framework so far.
The following theorem provides us with a general homogenization result in
all settings: the stochastic one, the coupled stochastic-deterministic one
and the deterministic one as well.

\begin{theorem}
\label{t4.1}Assume the following assumption holds: 
\begin{equation}
a_{ij}(x,\omega ,\cdot )\in A\text{\ for all }x\in \overline{Q}\text{ and }%
\mu \text{-a.e. }\omega \in \Omega \text{ }(1\leq i,j\leq N)\text{.}
\label{4.2}
\end{equation}%
For each fixed $\varepsilon >0$ let $u_{\varepsilon }$ be the unique
solution to \emph{(\ref{4.1})}. Then, as $\varepsilon \rightarrow 0$, 
\begin{equation}
u_{\varepsilon }\rightarrow u_{0}\text{\ stoch. in }L^{2}(Q\times \Omega )%
\text{-weak}  \label{4.3}
\end{equation}%
and%
\begin{equation}
\frac{\partial u_{\varepsilon }}{\partial x_{j}}\rightarrow \frac{\partial
u_{0}}{\partial x_{j}}+\overline{D}_{j,\omega }u_{1}+\frac{\overline{%
\partial }u_{2}}{\partial y_{j}}\text{ stoch. in }L^{2}(Q\times \Omega )%
\text{-weak }\Sigma \text{ }(1\leq j\leq N)  \label{4.4}
\end{equation}%
where the triple $(u_{0},u_{1},u_{2})\in \mathbb{F}%
^{1}=H_{0}^{1}(Q;I_{nv}^{2}(\Omega ))\times L^{2}(Q;\mathcal{W}^{1,2}(\Omega
))\times L^{2}(Q\times \Omega ;\mathcal{B}_{\#A}^{1,2})$ is the unique
solution to the following variational problem 
\begin{equation}
\left\{ 
\begin{array}{l}
\mathbf{u}=(u_{0},u_{1},u_{2})\in \mathbb{F}^{1}: \\ 
\underset{i,j=1}{\overset{N}{\sum }}\iint_{Q\times \Omega \times \Delta (A)}%
\widehat{a}_{ij}\mathbb{D}_{j}\mathbf{u}\overline{\mathbb{D}_{i}\mathbf{v}}%
dxd\mu d\beta =\left\langle f,\overline{v}_{0}\right\rangle \text{ for all }%
(v_{0},v_{1},v_{2})\in \mathbb{F}^{1}%
\end{array}%
\right.  \label{4.5}
\end{equation}%
with $\mathbb{D}_{j}\mathbf{u}=\partial u_{0}/\partial x_{j}+\overline{D}%
_{j,\omega }u_{1}+\mathcal{G}_{1}(\overline{\partial }u_{2}/\partial y_{j})$
(same definition for $\mathbb{D}_{i}\mathbf{v}$) and $\left\langle f,%
\overline{v}_{0}\right\rangle =\int_{\Omega }\left( f(\omega ),\overline{v}%
_{0}(\omega )\right) d\mu $, $\left( \cdot ,\cdot \right) $ denoting the
duality pairing between $H^{-1}(Q)$ and $H_{0}^{1}(Q)$.
\end{theorem}

\begin{proof}
We have 
\begin{equation}
\underset{i,j=1}{\overset{N}{\sum }}\int_{Q}a_{ij}^{\varepsilon }(\cdot
,\omega )\frac{\partial u_{\varepsilon }(\cdot ,\omega )}{\partial x_{j}}%
\overline{\frac{\partial v}{\partial x_{i}}}dx=\left( f(\omega ),\overline{v}%
\right)  \label{4.6}
\end{equation}%
for all $v\in H_{0}^{1}(Q)$, where $a_{ij}^{\varepsilon }(x,\omega
)=a_{ij}(x,T(x/\varepsilon _{1})\omega ,x/\varepsilon _{2})$ for $(x,\omega
)\in Q\times \Omega $. By taking the particular $v=u_{\varepsilon }(\cdot
,\omega )$ and and making use of the properties of the matrix $%
(a_{ij})_{1\leq i,j\leq N}$ and of the function $f$ we get the existence of
an absolute constant $c>0$ such that $\sup_{\varepsilon >0}\left\Vert
u_{\varepsilon }(\cdot ,\omega )\right\Vert _{H_{0}^{1}(Q)}\leq c$ for $\mu $%
-a.e. $\omega \in \Omega $. Hence, by Theorem \ref{t3.3} (where we take
there $X=H_{0}^{1}(Q)$), given any fundamental sequence $E$, there exists a
subsequence $E^{\prime }$ of $E$ and a triple $\mathbf{u}%
=(u_{0},u_{1},u_{2})\in \mathbb{F}^{1}$ such that, as $E^{\prime }\ni
\varepsilon \rightarrow 0$ we have (\ref{4.3})-(\ref{4.4}). Thus the theorem
will be proven as soon as we check that $\mathbf{u}$ verifies the
variational equation (\ref{4.5}). In fact it is easy to see that equation (%
\ref{4.5}) has at most one solution, so that checking that $\mathbf{u}$
verifies (\ref{4.5}) will prove that $\mathbf{u}$ does not depend on the
subsequence $E^{\prime }$, but on the whole sequence $\varepsilon >0$ which
will therefore establish Theorem \ref{t4.1}. Before we can do this, let us,
however notice that the space $\mathcal{F}_{0}^{\infty }=[\mathcal{D}%
(Q)\otimes (I_{nv}^{2}(\Omega ))]\times \lbrack \mathcal{D}(Q)\otimes I_{2}(%
\mathcal{C}^{\infty }(\Omega ))]\times \lbrack \mathcal{D}(Q)\otimes 
\mathcal{C}^{\infty }(\Omega )\otimes (J_{1}\circ \varrho )(A^{\infty }/%
\mathbb{C})]$ is dense in $\mathbb{F}^{1}$. Indeed, this comes from the fact
that $I_{2}(\mathcal{C}^{\infty }(\Omega ))$ (resp. $(J_{1}\circ \varrho
)(A^{\infty }/\mathbb{C})$) is dense in $\mathcal{W}^{1,2}(\Omega )$ (resp. $%
\mathcal{B}_{\#A}^{1,2}$), where $J_{1}$ (resp. $\varrho $, $I_{2}$) denotes
the canonical mapping of $\mathcal{B}_{A}^{1,2}/\mathbb{C}$ (resp. $%
B_{A}^{2} $, $\mathcal{C}^{\infty }(\Omega )$) into its separated completion 
$\mathcal{B}_{\#A}^{1,2}$ (resp. $\mathcal{B}_{A}^{2}$, $\mathcal{W}%
^{1,2}(\Omega )$); see Section 2 (particularly the Subsections 2.1 and 2.3
therein).

With this in mind, let $\Phi =(\psi _{0},I_{2}(\psi _{1}),(J_{1}\circ
\varrho )(\psi _{2}))\in \mathcal{F}_{0}^{\infty }$ and define 
\begin{equation*}
\Phi _{\varepsilon }(x,\omega )=\psi _{0}(x,\omega )+\varepsilon _{1}\psi
_{1}(x,T(x/\varepsilon _{1})\omega )+\varepsilon _{2}\psi
_{2}(x,T(x/\varepsilon _{1})\omega ,x/\varepsilon _{2})\text{ }((x,\omega
)\in Q\times \Omega )
\end{equation*}%
where $\psi _{0}\in \mathcal{C}_{0}^{\infty }(Q)\otimes I_{nv}^{2}(\Omega )$%
, $\psi _{1}\in \mathcal{C}_{0}^{\infty }(Q)\otimes \mathcal{C}^{\infty
}(\Omega )$ and $\psi _{2}\in \mathcal{C}_{0}^{\infty }(Q)\otimes \mathcal{C}%
^{\infty }(\Omega )\otimes (A^{\infty }/\mathbb{C})$; clearly $\Phi
_{\varepsilon }(\cdot ,\omega )\in \mathcal{C}_{0}^{\infty }(Q)$. Taking in (%
\ref{4.6}) the particular $v=\Phi _{\varepsilon }(\cdot ,\omega )$ and
integrating the resulting equality over $\Omega $ with respect to $\mu $, we
get 
\begin{equation}
\underset{i,j=1}{\overset{N}{\sum }}\int_{Q\times \Omega
}a_{ij}^{\varepsilon }\frac{\partial u_{\varepsilon }}{\partial x_{j}}%
\overline{\frac{\partial \Phi _{\varepsilon }}{\partial x_{i}}}dxd\mu
=\left\langle f,\overline{\Phi }_{\varepsilon }\right\rangle .  \label{4.7}
\end{equation}%
One easily show that as $\varepsilon \rightarrow 0$, $\left\langle f,%
\overline{\Phi }_{\varepsilon }\right\rangle \rightarrow \left\langle f,%
\overline{\psi }_{0}\right\rangle $ and $\partial \Phi _{\varepsilon
}/\partial x_{i}\rightarrow \partial \psi _{0}/\partial x_{i}+\overline{D}%
_{i,\omega }I_{2}(\psi _{1})+\overline{\partial }(J_{1}\circ \varrho )(\psi
_{2})/\partial y_{i}=\partial \psi _{0}/\partial x_{i}+D_{i,\omega }\psi
_{1}+\partial \psi _{2}/\partial y_{i}$ stoch. in $L^{2}(Q\times \Omega )$%
-strong $\Sigma $ ($1\leq i\leq N$). Putting together this convergence
result with (\ref{4.4}), we get by Theorem \ref{t3.2} that, as $E^{\prime
}\ni \varepsilon \rightarrow 0$, 
\begin{equation*}
\frac{\partial u_{\varepsilon }}{\partial x_{j}}\overline{\frac{\partial
\Phi _{\varepsilon }}{\partial x_{i}}}\rightarrow \mathfrak{D}_{j}\mathbf{u}%
\overline{\mathfrak{D}_{i}\Phi }\text{ stoch. in }L^{2}(Q\times \Omega )%
\text{-weak }\Sigma
\end{equation*}%
where $\mathfrak{D}_{j}\mathbf{u}=\partial u_{0}/\partial x_{j}+\overline{D}%
_{j,\omega }u_{1}+\overline{\partial }u_{2}/\partial x_{j}$ and $\mathfrak{D}%
_{i}\Phi =\partial \psi _{0}/\partial x_{i}+D_{i,\omega }\psi _{1}+\partial
\psi _{2}/\partial y_{i}$. Note that $\mathcal{G}_{1}(\mathfrak{D}_{j}%
\mathbf{u})=\mathbb{D}_{j}\mathbf{u}$ and $\mathcal{G}(\mathfrak{D}_{i}\Phi
)=\mathcal{G}_{1}(\varrho (\mathfrak{D}_{i}\Phi ))=\mathbb{D}_{i}\Phi $.
Hence a passage to the limit in (\ref{4.7}) using Proposition \ref{p3.4}
(recall that the function $a_{ij}\in \mathcal{C}(\overline{Q};L^{\infty
}(\Omega ;A))$ and is therefore an admissible function in the sense of
Definition \ref{d3.1'}) yields 
\begin{equation*}
\underset{i,j=1}{\overset{N}{\sum }}\iint_{Q\times \Omega \times \Delta (A)}%
\widehat{a}_{ij}\mathbb{D}_{j}\mathbf{u}\overline{\mathbb{D}_{i}\Phi }dxd\mu
d\beta =\left\langle f,\overline{\psi }_{0}\right\rangle \text{ for all }%
\Phi \in \mathcal{F}_{0}^{\infty }\text{.}
\end{equation*}%
We therefore get (\ref{4.5}) by density of $\mathcal{F}_{0}^{\infty }$ in $%
\mathbb{F}^{1}$. This completes the proof of the theorem.
\end{proof}

Now, we need to show that Theorem \ref{t4.1} generalizes all the existing
results in the framework of homogenization of linear elliptic equations. To
that end, we will distinguish three cases: (1) the case when the functions $%
a_{ij}$ and $f$ do not depend on the random variable $\omega $; (2) the case
when the functions $a_{ij}$ do not depend on the deterministic variable $y$
and, (3) the case when the functions $a_{ij}$ depend upon both $\omega $ and 
$y$ and the function $f$ does not depend on $\omega $, but with $\Omega
=\Delta (A_{z})$ where $A_{z}$ is some algebra with mean value on $\mathbb{R}%
_{z}^{N}$.\medskip

For the first case we have $a_{ij}(x,\omega ,y)=a_{ij}(x,y)$ and $f(x,\omega
)=f(x)$ for $(x,\omega ,y)\in Q\times \Omega \times \mathbb{R}_{y}^{N}$. In
that case, the problem (\ref{4.1}) is a deterministic one, and its solution $%
u_{\varepsilon }$ does not depend on $\omega $. A rapid survey of the proof
of Theorem \ref{t3.3} gives, by the Remark \ref{r3.0}, that the functions $%
\mathbf{v}$ and $u_{0}$ therein can be chosen in $L^{p}(Q;\mathcal{B}%
_{A}^{p})^{N}$ and in $L^{p}(Q;\mathcal{B}_{A}^{p})$, respectively. This
yields immediately the fact that the function $u_{1}$ there is equal to
zero, so that $\mathbf{v}=Du_{0}+\overline{D}_{y}u_{2}$ with $u_{2}\in
L^{p}(Q;\mathcal{B}_{\#A}^{1,p})$. The continuity assumption on $%
a_{ij}(x,\cdot )$ can therefore be replaced by a measurability assumption $%
a_{ij}(x,\cdot )\in L^{\infty }(\mathbb{R}_{y}^{N})$, so that the
homogenization result for (\ref{4.1}) therefore reads as

\begin{theorem}
\label{t4.2}Assume the following assumption holds: 
\begin{equation*}
a_{ij}(x,\cdot )\in B_{A}^{2}\text{\ for all }x\in \overline{Q}\text{ }%
(1\leq i,j\leq N)\text{.}
\end{equation*}%
For each fixed $\varepsilon >0$ let $u_{\varepsilon }$ be the unique
solution to \emph{(\ref{4.1})}. Then, as $\varepsilon \rightarrow 0$, 
\begin{equation*}
u_{\varepsilon }\rightarrow u_{0}\text{ in }H_{0}^{1}(Q)\text{-weak}
\end{equation*}%
and%
\begin{equation*}
\frac{\partial u_{\varepsilon }}{\partial x_{j}}\rightarrow \frac{\partial
u_{0}}{\partial x_{j}}+\frac{\overline{\partial }u_{2}}{\partial y_{j}}\text{
in }L^{2}(Q)\text{-weak }\Sigma \text{ }(1\leq j\leq N)
\end{equation*}%
where the couple $(u_{0},u_{2})\in \mathbb{F}^{1}=H_{0}^{1}(Q)\times L^{2}(Q;%
\mathcal{B}_{\#A}^{1,2})$ is the unique solution to the following
variational problem 
\begin{equation}
\left\{ 
\begin{array}{l}
\mathbf{u}=(u_{0},u_{2})\in \mathbb{F}^{1}: \\ 
\underset{i,j=1}{\overset{N}{\sum }}\iint_{Q\times \Delta (A)}\widehat{a}%
_{ij}\left( \frac{\partial u_{0}}{\partial x_{j}}+\mathcal{G}_{1}(\frac{%
\overline{\partial }u_{2}}{\partial y_{j}})\right) \overline{\left( \frac{%
\partial v_{0}}{\partial x_{i}}+\mathcal{G}_{1}(\frac{\overline{\partial }%
v_{2}}{\partial y_{i}})\right) }dxd\beta =\left( f,\overline{v}_{0}\right)
\\ 
\text{for all }(v_{0},v_{2})\in \mathbb{F}^{1}.%
\end{array}%
\right.  \label{4.8}
\end{equation}
\end{theorem}

Equation (\ref{4.8}) is well-known in the literature of deterministic
homogenization; see in particular \cite{26}.\medskip

For the case of a random operator ($a_{ij}(x,\omega ,y)=a_{ij}(x,\omega )$),
a similar type of reasoning as the one in the previous case (using once
again the Remark \ref{r3.0}) yields the following result.

\begin{theorem}
\label{t4.3}For each fixed $\varepsilon >0$ let $u_{\varepsilon }$ be the
unique solution to \emph{(\ref{4.1})}. Then, as $\varepsilon \rightarrow 0$, 
\begin{equation*}
u_{\varepsilon }\rightarrow u_{0}\text{\ stoch. in }L^{2}(Q\times \Omega )%
\text{-weak}
\end{equation*}%
and%
\begin{equation*}
\frac{\partial u_{\varepsilon }}{\partial x_{j}}\rightarrow \frac{\partial
u_{0}}{\partial x_{j}}+\overline{D}_{j,\omega }u_{1}\text{ stoch. in }%
L^{2}(Q\times \Omega )\text{-weak }(1\leq j\leq N)
\end{equation*}%
where the couple $(u_{0},u_{1})$ is the unique solution of the following
variational problem 
\begin{equation}
\left\{ 
\begin{array}{l}
(u_{0},u_{1})\in \mathbb{F}^{1}=H_{0}^{1}(Q;I_{nv}^{2}(\Omega ))\times
L^{2}(Q;\mathcal{W}^{1,2}(\Omega )) \\ 
\underset{i,j=1}{\overset{N}{\sum }}\iint_{Q\times \Omega }a_{ij}\left( 
\frac{\partial u_{0}}{\partial x_{j}}+\overline{D}_{j,\omega }u_{1}\right) 
\overline{\left( \frac{\partial v_{0}}{\partial x_{i}}+\overline{D}%
_{i,\omega }v_{1}\right) }dxd\mu =\left\langle f,\overline{v}%
_{0}\right\rangle \\ 
\text{for all }(v_{0},v_{1})\in \mathbb{F}^{1}.%
\end{array}%
\right.  \label{4.9}
\end{equation}
\end{theorem}

Equation (\ref{4.9}) is also well-known in the literature; see \cite{10,
Kozlov1, Kozlov2}.\medskip

For the last case, we assume that $\Omega =\Delta (A_{z})$, $A_{z}$ being
some algebra with mean value on $\mathbb{R}_{z}^{N}$ for the action $%
\mathcal{H}=(H_{\varepsilon })_{\varepsilon >0}$ defined by $H_{\varepsilon
}(x)=x/\varepsilon _{1}$ ($x\in \mathbb{R}^{N}$). The ergodic algebra $A$
(in Theorem \ref{t4.1}) is denoted here by $A_{y}$ and its associated $M$%
-measure by $\beta _{y}$. We use the same letter $\mathcal{G}$ to denote the
Gelfand transformation on $A_{y}$ and on $A_{z}$ as well. Points in $\Delta
(A_{y})$ (resp. $\Delta (A_{z})$) are denoted by $s$ (resp. $\omega $). The
compact space $\Delta (A_{z})$ is equipped with the $M$-measure $\beta _{z}$
for $A_{z}$. We know by Theorem \ref{t2.2} that one can define a dynamical
system on the spectrum $\Delta (A_{z})$ of $A_{z}$ so that the invariant
probability measure is precisely the $M$-measure $\beta _{z}$ for $A_{z}$.

With the above preliminaries, our concern here is not to reformulate the
statement of Theorem \ref{t4.1}, but to show how it includes the general
setting of reiterated deterministic homogenization. For that purpose, let $%
b_{ij}\in \mathcal{C}(\overline{Q};L^{\infty }(\mathbb{R}_{z}^{N};\mathcal{B}%
(\mathbb{R}_{y}^{N})))$ with $b_{ij}=\overline{b}_{ji}$ and $(b_{ij})_{1\leq
i,j\leq N}$ satisfying the following ellipticity condition: there exists a
constant $\alpha >0$ such that $\sum_{i,j=1}^{N}b_{ij}(x,z,y)\lambda _{i}%
\overline{\lambda }_{j}\geq \alpha \left\vert \lambda \right\vert ^{2}$ for
all $(x,y)\in \overline{Q}\times \mathbb{R}^{N}$, for almost all $z\in 
\mathbb{R}^{N}$ and for all $\lambda \in \mathbb{C}^{N}$. Assume moreover
that the following hypothesis is satisfied: 
\begin{equation}
b_{ij}(x,\cdot ,y)\in B_{A_{z}}^{2}\text{ and }b_{ij}(x,z,\cdot )\in A_{y}%
\text{ a.e. }(z,y)\in \mathbb{R}_{z}^{N}\times \mathbb{R}_{y}^{N}\text{ }%
(1\leq i,j\leq N).  \label{4.9'}
\end{equation}%
Since by the construction of the dynamical system $T(z)$ on $\Delta (A_{z})$
(see the proof of Theorem \ref{t2.2}) we have 
\begin{equation*}
\mathcal{G}(b_{ij}(x,\cdot +z,y))(\omega )=\mathcal{G}(b_{ij}(x,\cdot
,y))(T(z)\omega )
\end{equation*}%
(for all $z\in \mathbb{R}^{N}$, for almost $\omega \in \Delta (A_{z})$ and
for any fixed $(x,y)\in Q\times \mathbb{R}_{y}^{N}$), we set in a natural
way 
\begin{equation*}
a_{ij}(x,\omega ,y)=\mathcal{G}(b_{ij}(x,\cdot ,y))(\omega )\ ((x,\omega
,y)\in Q\times \Delta (A_{z})\times \mathbb{R}^{N}),
\end{equation*}%
$\mathcal{G}$ being here the extension of the Gelfand transformation to the
Besicovitch space $B_{A_{z}}^{2}$; see property (\textbf{2}) in Subsection
2.3. Then, in view of the properties of $\mathcal{G}$, the functions $a_{ij}$
thus defined satisfy all the requirements of Theorem \ref{t4.1} so that we
get a homogenization result for the following problem: 
\begin{equation}
-\overset{N}{\underset{i,j=1}{\sum }}\frac{\partial }{\partial x_{i}}\left(
b_{ij}\left( x,\frac{x}{\varepsilon _{1}},\frac{x}{\varepsilon _{2}}\right) 
\frac{\partial u_{\varepsilon }}{\partial x_{j}}\right) =f\text{ in }Q\text{%
, }u_{\varepsilon }\in H_{0}^{1}(Q).  \label{4.10}
\end{equation}%
It is important to note that we do not say that problem (\ref{4.1}) is
equivalent to problem (\ref{4.10}). In fact let us assume that the functions 
$b_{ij}$ are such that $b_{ij}(x,\cdot ,y)\in A_{z}$. Then 
\begin{eqnarray*}
\frac{\partial }{\partial x_{i}}\left( a_{ij}\left( x,T\left( \frac{x}{%
\varepsilon _{1}}\right) \omega ,\frac{x}{\varepsilon _{2}}\right) \frac{%
\partial u_{\varepsilon }}{\partial x_{j}}\right) &=&\frac{\partial }{%
\partial x_{i}}\left( \mathcal{G}\left( b_{ij}\left( x,\cdot ,\frac{x}{%
\varepsilon _{2}}\right) \right) \left( T\left( \frac{x}{\varepsilon _{1}}%
\right) \omega \right) \frac{\partial u_{\varepsilon }}{\partial x_{j}}%
\right) \\
&=&\frac{\partial }{\partial x_{i}}\left( \mathcal{G}\left( b_{ij}\left(
x,\cdot +\frac{x}{\varepsilon _{1}},\frac{x}{\varepsilon _{2}}\right)
\right) (\omega )\frac{\partial u_{\varepsilon }}{\partial x_{j}}\right) \\
&=&\mathcal{G}\left[ \frac{\partial }{\partial x_{i}}\left( b_{ij}\left(
x,\cdot +\frac{x}{\varepsilon _{1}},\frac{x}{\varepsilon _{2}}\right) \frac{%
\partial u_{\varepsilon }}{\partial x_{j}}\right) \right] (\omega ),\text{ }
\end{eqnarray*}%
the last equality being due to the fact that $\mathcal{G}$ is linear
continuous. If in particular $\omega =\delta _{z}$ ($z\in \mathbb{R}^{N}$),
the Dirac mass at $z$, then 
\begin{eqnarray*}
\frac{\partial }{\partial x_{i}}\left( a_{ij}\left( x,T\left( \frac{x}{%
\varepsilon _{1}}\right) \omega ,\frac{x}{\varepsilon _{2}}\right) \frac{%
\partial u_{\varepsilon }}{\partial x_{j}}\right) &=&\left\langle \delta
_{z},\frac{\partial }{\partial x_{i}}\left( b_{ij}\left( x,\cdot +\frac{x}{%
\varepsilon _{1}},\frac{x}{\varepsilon _{2}}\right) \frac{\partial
u_{\varepsilon }}{\partial x_{j}}\right) \right\rangle \\
&=&\frac{\partial }{\partial x_{i}}\left( b_{ij}\left( x,z+\frac{x}{%
\varepsilon _{1}},\frac{x}{\varepsilon _{2}}\right) \frac{\partial
u_{\varepsilon }}{\partial x_{j}}\right) ,
\end{eqnarray*}%
$\left\langle ,\right\rangle $ denoting the duality pairing between $%
A_{z}^{\prime }$ and $A_{z}$ (see Subsection 2.2). Thus, in this particular
case, Equation (\ref{4.1}) becomes 
\begin{equation}
-\overset{N}{\underset{i,j=1}{\sum }}\frac{\partial }{\partial x_{i}}\left(
b_{ij}\left( x,z+\frac{x}{\varepsilon _{1}},\frac{x}{\varepsilon _{2}}%
\right) \frac{\partial u_{\varepsilon }}{\partial x_{j}}\right) =f\text{ in }%
Q\text{, }u_{\varepsilon }\in H_{0}^{1}(Q).  \label{4.11}
\end{equation}%
One therefore sees that (\ref{4.10}) comes from (\ref{4.11}) by taking there 
$z=0$. It is also to be noted that if the algebra $A_{z}$ is ergodic, then
the dynamical system $T(z)$ is ergodic, in such a way that (\ref{4.11}) is
equivalent to (\ref{4.10}). However, still assuming $A_{z}$ to be ergodic
and taking the functions $b_{ij}$ in $\mathcal{C}(\overline{Q};L^{\infty }(%
\mathbb{R}_{z}^{N};\mathcal{B}(\mathbb{R}_{y}^{N})))$ with (\ref{4.9'}), and
finally arguing as in \cite[Section 4]{Efendiev}, one also obtains the
equivalence of (\ref{4.10}) and (\ref{4.11}). We also note here that the
algebra $A_{z}$ is not assumed to be ergodic in general, so that we have a
great flexibility in Theorem \ref{t4.1} in the particular case where $\Omega
=\Delta (A_{z})$. Indeed, Theorem \ref{t4.1} works in all the environments:
the ergodic one and the non ergodic one for $A_{z}$. In the particular case
when the algebra $A_{z}$ is ergodic, $I_{nv}^{2}(\Delta (A_{z}))$ consists
of constants, so that the function $u_{0}$ lies in $H_{0}^{1}(Q)$. We thus
recover the well-known results in that environment. If we assume that $A_{z}$
is not ergodic, our result is then new. For the sake of completeness, let us
give some concrete situations in which our result applies.

\begin{example}[Homogenization in ergodic algebras]
\label{e1}\emph{One can solve the homogenization problem for (\ref{4.10})
under each of the following hypotheses:}

\begin{itemize}
\item[(H)$_{1}$] \emph{The function }$b_{ij}(x,\cdot ,\cdot )$\emph{\ is
periodic in }$y$\emph{\ and in }$z$\emph{;}

\item[(H)$_{2}$] \emph{The function }$b_{ij}(x,\cdot ,\cdot )$\emph{\ is
almost periodic in }$y$\emph{\ and in }$z$\emph{;}

\item[(H)$_{3}$] \emph{The function }$b_{ij}(x,\cdot ,y)$\emph{\ is almost
periodic and the function }$b_{ij}(x,z,\cdot )$\emph{\ is weakly almost
periodic \cite{17};}

\item[(H)$_{4}$] \emph{The functions }$b_{ij}(x,\cdot ,y)$\emph{\ and }$%
b_{ij}(x,z,\cdot )$\emph{\ are both weakly almost periodic.}
\end{itemize}
\end{example}

\begin{example}[Homogenization in non ergodic algebras]
\label{e2}\emph{For the sake of simplicity, we assume here that }$N=1$\emph{%
. Let }$A_{z}$\emph{\ be the algebra generated by the function }$f(z)=\cos 
\sqrt[3]{z}$\emph{\ (}$z\in \mathbb{R}$\emph{) and all its translates }$%
f(\cdot +a)$\emph{, }$a\in \mathbb{R}$\emph{. It is known that }$A$\emph{\
is an algebra with mean value which is not ergodic; see \cite{20} for
details. However, as said above, one can solve the homogenization problem
for (\ref{4.10}) under the following hypothesis: }$b_{ij}(x,\cdot ,y)\in
B_{A_{z}}^{2}$\emph{\ and }$b_{ij}(x,z,\cdot )\in A_{y}$\emph{, where }$%
A_{y} $\emph{\ is any ergodic algebra with mean value on }$\mathbb{R}$\emph{%
. The homogenization problem solved here is new. One can also consider other
homogenization problems in the present setting of non ergodic algebras.}
\end{example}

\section{Application to the homogenization of nonlinear Reynolds-type
equations}

In this section we study the homogenization problem for nonlinear
Reynolds-type equations. More precisely, let $1<p<\infty $ be fixed and let
the function $(x,\omega ,y,\lambda )\mapsto a(x,\omega ,y,\lambda )$ from $%
\overline{Q}\times \Omega \times \mathbb{R}^{N}\times \mathbb{R}^{N}$ to $%
\mathbb{R}^{N}$ satisfy the following conditions: 
\begin{equation}
\begin{array}{l}
\text{For all }(x,y,\lambda )\in \mathbb{R}^{N}\times \mathbb{R}^{N}\times 
\mathbb{R}^{nN}\text{ and for almost all }\omega \in \Omega \text{, }%
a(x,\cdot ,y,\lambda ) \\ 
\text{is measurable and }a(\cdot ,\omega ,\cdot ,\cdot )\text{\ is
continuous;}%
\end{array}
\label{5.1}
\end{equation}%
\begin{equation}
\begin{array}{l}
\text{There are four constants }c_{0},\,c_{1},c_{2}>0\text{, }0<\alpha \leq
\min (1,p-1) \\ 
\text{and a continuity modulus }\nu \text{ (i.e., a nondecreasing continuous}
\\ 
\text{function on }[0,+\infty )\text{ such that }\nu (0)=0,\nu (r)>0\text{\
if }r>0\text{, and} \\ 
\nu (r)=1\text{ if }r>1\text{) such that for a.e. }y\in \mathbb{R}^{N}\text{
and for }\mu \text{-a.e }\omega \in \Omega \text{,} \\ 
\text{(i) }a(x,\omega ,y,\lambda )\cdot \lambda \geq c_{0}\left\vert \lambda
\right\vert ^{p} \\ 
\text{(ii) }\left\vert a(x,\omega ,y,\lambda )\right\vert \leq c_{1}\left(
1+\left\vert \lambda \right\vert ^{p-1}\right) \\ 
\text{(iii) }\left\vert a(x,\omega ,y,\lambda )-a(x^{\prime },\omega
,y,\lambda ^{\prime })\right\vert \leq \nu (\left\vert x-x^{\prime
}\right\vert )(1+\left\vert \lambda \right\vert ^{p-1}+\left\vert \lambda
^{\prime }\right\vert ^{p-1}) \\ 
\;\;\;\;\;\;\;\;+c_{2}\left( 1+\left\vert \lambda \right\vert +\left\vert
\lambda ^{\prime }\right\vert \right) ^{p-1-\alpha }\left\vert \lambda
-\lambda ^{\prime }\right\vert ^{\alpha } \\ 
\text{(iv) }\left( a(x,\omega ,y,\lambda )-a(x,\omega ,y,\lambda ^{\prime
})\right) \cdot \left( \lambda -\lambda ^{\prime }\right) >0\text{\ if }%
\lambda \neq \lambda ^{\prime } \\ 
\text{for all }x,x^{\prime }\in \overline{Q}\text{, all }\lambda ,\lambda
^{\prime }\in \mathbb{R}^{N}\text{, where the dot denotes the usual} \\ 
\text{Euclidean inner product in }\mathbb{R}^{N}\text{ and }\left\vert
\,\cdot \,\right\vert \text{ the associated norm,} \\ 
Q\text{ being a bounded open set in }\mathbb{R}^{N}\text{.}%
\end{array}
\label{5.3}
\end{equation}

We consider the boundary value problem 
\begin{equation}
\begin{array}{l}
u_{\varepsilon }(\cdot ,\omega )\in W_{0}^{1,p}(\Omega ): \\ 
\func{div}a\left( x,T\left( \frac{x}{\varepsilon _{1}}\right) \omega ,\frac{x%
}{\varepsilon _{2}},Du_{\varepsilon }\right) =\func{div}b\left( x,T\left( 
\frac{x}{\varepsilon _{1}}\right) \omega ,\frac{x}{\varepsilon _{2}}\right) 
\text{\ in }Q%
\end{array}
\label{5.4}
\end{equation}%
with $b\in \mathcal{C}(\overline{Q};L^{\infty }(\Omega ;\mathcal{B}(\mathbb{R%
}_{y}^{N})))$, where we assume that the scales $\varepsilon _{1}$ and $%
\varepsilon _{2}$ are well-separated as in Section 3. It is easily seen that
the realization $a(x,T(z)\omega ,y,\lambda )$ is well-defined for almost all 
$\omega \in \Omega $, such that the functions $x\mapsto a(x,T(x/\varepsilon
_{1})\omega ,x/\varepsilon _{2},\mathbf{v}(x))$ of $Q$ into $\mathbb{R}^{N}$
(defined as element of $L^{p^{\prime }}(Q)^{N}$ for $\mathbf{v}\in
L^{p}(Q)^{N}$) and $(x,\omega )\mapsto a(x,T(x/\varepsilon _{1})\omega
,x/\varepsilon _{2},\Psi (x,T(x/\varepsilon _{1})\omega ,x/\varepsilon
_{2})) $ (for $\Psi \in \mathcal{C}(\overline{Q};L^{\infty }(\Omega ;%
\mathcal{B}(\mathbb{R}_{y}^{N})^{N}))$) of $Q\times \Omega $ into $\mathbb{R}%
^{N}$ denoted by $a^{\varepsilon }(\cdot ,\Psi ^{\varepsilon })$ (as element
of $L^{\infty }(Q\times \Omega )^{N}$) are well-defined. With this in mind,
we see that the problem (\ref{5.4}) admits a unique solution $u_{\varepsilon
}\in W_{0}^{1,p}(Q;L^{p}(\Omega ))$ (for each fixed $\varepsilon >0$); see
e.g., \cite[Chap. 2]{Lions}.

The main advantage of considering this problem lies in its application in
hydrodynamics. One of the difficulties encountered in homogenizing the above
problem is that the right-hand side of (\ref{5.4}) depends upon $\varepsilon 
$ and rather weakly converges in $L^{p^{\prime }}(\Omega ;W^{-1,p^{\prime
}}(Q))$, contrary to what is usually considered in the literature.

Throughout the rest of this section, all the vector spaces are assumed to be
real vector spaces, and the scalar functions are assumed real valued.
Obviously, this entails that the results of Section 3 are still valid, the
only difference being that all the function spaces are real. Now, let $A$ be
an ergodic $H$-supralgebra on $\mathbb{R}_{y}^{N}$. Our goal here is to
investigate the limiting behavior of $(u_{\varepsilon })_{\varepsilon >0}$
(the sequence of solutions to (\ref{5.4})) under the assumptions 
\begin{equation}
b\in \mathcal{C}(\overline{Q};L^{\infty }(\Omega ;A))\ \ \ \ \ \ \ \ \ \ \ \
\ \ \ \ \ \ \ \ \ \ \ \ \ \ \ \ \ \ \ \
\;\;\;\;\;\;\;\;\;\;\;\;\;\;\;\;\;\;\;\;\;\;\;\;  \label{5.7}
\end{equation}%
\begin{equation}
a_{i}(x,\omega ,\cdot ,\lambda )\in A\text{\ for all }(x,\omega ,\lambda
)\in \overline{Q}\times \Omega \times \mathbb{R}^{N}\;(1\leq i\leq N)
\label{5.8}
\end{equation}%
where $a_{i}$ denotes the $i$th component of the function $a$. Assuming (\ref%
{5.8}), it follows as in the proof of Theorem \ref{t5} that, for any $\Psi
\in \lbrack \mathcal{C}_{0}^{\infty }(Q)\otimes \mathcal{C}^{\infty }(\Omega
)\otimes A]^{N}$, the function $a(\cdot ,\Psi ):(x,\omega ,y)\mapsto
a(x,\omega ,y,\Psi (x,\omega ,y))$ lies in $\mathcal{C}(\overline{Q}%
;L^{\infty }(\Omega ;A)^{N})$ so that by Proposition \ref{p3.4} we have the
following convergence result: 
\begin{equation}
a_{i}^{\varepsilon }(\cdot ,\Psi ^{\varepsilon })\rightarrow \varrho \circ
a_{i}(\cdot ,\Psi )\text{\ stoch. in }L^{p^{\prime }}(Q\times \Omega )\text{%
-weak }\Sigma \text{ as }\varepsilon \rightarrow 0\text{ }(1\leq i\leq N),
\label{5.9}
\end{equation}%
$\varrho $ being the canonical mapping from $B_{A}^{p^{\prime }}$ into $%
\mathcal{B}_{A}^{p^{\prime }}$ and $(\varrho \circ a_{i}(\cdot ,\Psi
))(x,\omega ,y)=\varrho (a_{i}(x,\omega ,\cdot ,\Psi (x,\omega ,\cdot )))(y)$
for $(x,\omega ,y)\in Q\times \Omega \times \mathbb{R}^{N}$. Moreover the
following result holds.

\begin{proposition}
\label{p5.2}The mapping $\Psi \mapsto a(\cdot ,\Psi )$ from $[\mathcal{C}%
_{0}^{\infty }(Q)\otimes \mathcal{C}^{\infty }(\Omega )\otimes A]^{N}$ into $%
L^{p^{\prime }}(Q\times \Omega ;B_{A}^{p^{\prime }})^{N}$, extends by
continuity to a unique mapping still denoted by $a$, of $L^{p}(Q\times
\Omega ;(B_{A}^{p})^{N})$ into $L^{p^{\prime }}(Q\times \Omega
;B_{A}^{p^{\prime }})^{N}$ such that 
\begin{equation*}
(a(\cdot ,\mathbf{v})-a(\cdot ,\mathbf{w}))\cdot (\mathbf{v}-\mathbf{w})\geq
0\text{\ a.e. in }Q\times \Omega \times \mathbb{R}_{y}^{N}
\end{equation*}%
\begin{equation*}
\left\| a(\cdot ,\mathbf{v})\right\| _{L^{p^{\prime }}(Q\times \Omega
;B_{A}^{p^{\prime }})^{N}}\leq c_{1}^{\prime }(1+\left\| \mathbf{v}\right\|
_{L^{p}(Q\times \Omega ;(B_{A}^{p})^{N})}^{p-1})
\end{equation*}%
\begin{equation*}
\left\| a(\cdot ,\mathbf{v})-a(\cdot ,\mathbf{w})\right\| _{L^{p^{\prime
}}(Q\times \Omega ;B_{A}^{p^{\prime }})^{N}}\leq c_{2}\left\| 1+\left| 
\mathbf{v}\right| +\left| \mathbf{w}\right| \right\| _{L^{p}(Q\times \Omega
;B_{A}^{p})}^{p-1-\alpha }\left\| \mathbf{v}-\mathbf{w}\right\|
_{L^{p}(Q\times \Omega ;B_{A}^{p})^{N}}^{\alpha }
\end{equation*}%
\begin{equation*}
\left| a(x,\omega ,y,\mathbf{w})-a(x^{\prime },\omega ,y,\mathbf{w})\right|
\leq \nu (\left| x-x^{\prime }\right| )(1+\left| \mathbf{w}\right| ^{p-1})%
\text{ a.e. in }\Omega \times \mathbb{R}_{y}^{N}
\end{equation*}%
for all $\mathbf{v},\mathbf{w}\in L^{p}(Q\times \Omega ;(B_{A}^{p})^{N})$
and all $x,x^{\prime }\in Q$, where the constant $c_{1}^{\prime }$ depends
only on $c_{1}$ and on $Q$.
\end{proposition}

\begin{proof}
It is immediate that for $\Psi \in \lbrack \mathcal{C}_{0}^{\infty
}(Q)\otimes \mathcal{C}^{\infty }(\Omega )\otimes A]^{N}$, the function $%
a(\cdot ,\Psi )$ verifies properties of the same type as in (\ref{5.3}) (see
in particular inequality (iii) therein), so that arguing as in the proof of %
\cite[Proposition 3.1]{36} we get the result.
\end{proof}

As a consequence of the convergence result (\ref{5.9}) we have the following
result whose proof is quite similar to that of \cite[Corollary 3.9]{40} (see
also \cite[Corollary 3.3]{CPAA}).

\begin{proposition}
\label{p5.3}For $\psi _{0}\in \mathcal{C}_{0}^{\infty }(Q)\otimes
I_{nv}^{p}(\Omega )$, $\psi _{1}\in \mathcal{C}_{0}^{\infty }(Q)\otimes 
\mathcal{C}^{\infty }(\Omega )$ and $\psi _{2}\in \mathcal{C}_{0}^{\infty
}(Q)\otimes \mathcal{C}^{\infty }(\Omega )\otimes A^{\infty }$, define the
function $\Phi _{\varepsilon }$ \emph{(}$\varepsilon >0$\emph{)} by 
\begin{equation}
\Phi _{\varepsilon }=\psi _{0}+\varepsilon _{1}\psi _{1}^{\varepsilon
}+\varepsilon _{2}\psi _{2}^{\varepsilon },\;\;\;\;\;\;\;\;\;\;\;\;\;\;
\label{5.10}
\end{equation}%
i.e., $\Phi _{\varepsilon }(x,\omega )=\psi _{0}(x,\omega )+\varepsilon
_{1}\psi _{1}(x,T(x/\varepsilon _{1})\omega )+\varepsilon _{2}\psi
_{2}(x,T(x/\varepsilon _{1})\omega ,x/\varepsilon _{2})$ $((x,\omega )\in
Q\times \Omega )$. Let $(v_{\varepsilon })_{\varepsilon \in E}$ be a
sequence in $L^{p}(Q\times \Omega )^{N}$ such that $v_{\varepsilon
}\rightarrow v_{0}$ stoch. in $L^{p}(Q\times \Omega )^{N}$-weak $\Sigma $ as 
$E\ni \varepsilon \rightarrow 0$ where $v_{0}\in L^{p}(Q\times \Omega ;%
\mathcal{B}_{A}^{p})^{N}$. Then, as $E\ni \varepsilon \rightarrow 0$, 
\begin{equation*}
\int_{Q\times \Omega }a^{\varepsilon }(\cdot ,D\Phi _{\varepsilon })\cdot
v_{\varepsilon }dxd\mu \rightarrow \iint_{Q\times \Omega \times \Delta (A)}%
\widehat{a}(\cdot ,D\psi _{0}+D_{\omega }\psi _{1}+\partial \widehat{\psi }%
_{2})\cdot \widehat{v}_{0}dxd\mu d\beta
\end{equation*}%
where $\partial \widehat{\psi }_{2}=(\partial _{j}\widehat{\psi }%
_{2})_{1\leq j\leq N}$ with $\partial _{j}\widehat{\psi }_{2}=\mathcal{G}%
(\partial \psi _{2}/\partial y_{j})$.
\end{proposition}

We recall that the algebra $A$ is as stated earlier in this section. For $%
1<p<\infty $ we put $\mathbb{F}_{0}^{1,p}=W_{0}^{1,p}(Q;I_{nv}^{p}(\Omega
))\times L^{p}(Q;\mathcal{W}^{1,p}(\Omega ))\times L^{p}(Q\times \Omega ;%
\mathcal{B}_{\#A}^{1,p})$. We endow $\mathbb{F}_{0}^{1,p}$ with the norm 
\begin{eqnarray*}
\left\Vert \mathbf{u}\right\Vert _{\mathbb{F}_{0}^{1,p}} &=&\sum_{i=1}^{N}%
\left[ \left\Vert D_{x_{i}}u_{0}\right\Vert _{L^{p}(Q\times \Omega
)}+\left\Vert \overline{D}_{i,\omega }u_{1}\right\Vert _{L^{p}(Q\times
\Omega )}+\left\Vert \overline{D}_{y_{i}}u_{2}\right\Vert _{L^{p}(Q\times
\Omega ;\mathcal{B}_{A}^{p})}\right] \\
\mathbf{u} &=&(u_{0},u_{1},u_{2})\in \mathbb{F}_{0}^{1,p}.
\end{eqnarray*}%
In this norm, $\mathbb{F}_{0}^{1,p}$ is a Banach space admitting $\mathcal{F}%
_{0}^{\infty }=[\mathcal{C}_{0}^{\infty }(Q)\otimes (I_{nv}^{p}(\Omega
))]\times \lbrack \mathcal{C}_{0}^{\infty }(Q)\otimes I_{p}(\mathcal{C}%
^{\infty }(\Omega ))]\times \lbrack \mathcal{C}_{0}^{\infty }(Q)\otimes 
\mathcal{C}^{\infty }(\Omega )\otimes (J_{1}\circ \varrho )(A^{\infty }/%
\mathbb{C})]$ as a dense subspace where, $J_{1}$ (resp. $\varrho $, $I_{p}$)
denotes the canonical mapping of $\mathcal{B}_{A}^{1,p}/\mathbb{C}$ (resp. $%
B_{A}^{p}$, $\mathcal{C}^{\infty }(\Omega )$) into its separated completion $%
\mathcal{B}_{\#A}^{1,p}$ (resp. $\mathcal{B}_{A}^{p}$, $\mathcal{W}%
^{1,p}(\Omega )$). With this in mind, we have the following homogenization
result.

\begin{theorem}
\label{t5.1}Let $1<p<\infty $. Assume \emph{(\ref{5.7})} and \emph{(\ref{5.8}%
)} hold with $A$ an ergodic $H$-supralgebra on $\mathbb{R}^{N}$ which is
moreover an algebra with mean value. For each real $\varepsilon >0$, let $%
u_{\varepsilon }$ be the unique solution of \emph{(\ref{5.4})}. Then, as $%
\varepsilon \rightarrow 0$, 
\begin{equation}
u_{\varepsilon }\rightarrow u_{0}\text{\ stoch. in }L^{p}(Q\times \Omega )%
\text{-weak}  \label{5.11}
\end{equation}%
and%
\begin{equation}
\frac{\partial u_{\varepsilon }}{\partial x_{j}}\rightarrow \frac{\partial
u_{0}}{\partial x_{j}}+\overline{D}_{j,\omega }u_{1}+\frac{\overline{%
\partial }u_{2}}{\partial y_{j}}\text{ stoch. in }L^{p}(Q\times \Omega )%
\text{-weak }\Sigma \text{ }(1\leq j\leq N)  \label{5.12}
\end{equation}%
where $\mathbf{u}=(u_{0},u_{1},u_{2})\in \mathbb{F}_{0}^{1,p}$ is the unique
solution of the variational equation 
\begin{equation}
\begin{array}{l}
\iint_{Q\times \Omega \times \Delta (A)}\widehat{a}(\cdot ,\mathbb{D}\mathbf{%
u})\cdot \mathbb{D}\mathbf{v}dxd\mu d\beta =\int_{Q\times \Omega \times
\Delta (A)}\widehat{b}(x,\omega ,s)\widehat{\func{div}}\mathbf{v}dxd\mu
d\beta \\ 
\text{for all }\mathbf{v}=(v_{0},v_{1},v_{2})\in \mathbb{F}_{0}^{1,p}%
\end{array}
\label{5.13}
\end{equation}%
with $\mathbb{D}\mathbf{w}=Dw_{0}+\overline{D}_{\omega }w_{1}+\mathcal{G}%
_{1}^{N}(\overline{D}_{y}w_{2})$ for $\mathbf{w}=(w_{0},w_{1},w_{2})\in 
\mathbb{F}_{0}^{1,p}$ where: $\mathcal{G}_{1}^{N}(\overline{D}_{y}w_{2})=(%
\mathcal{G}_{1}(\overline{\partial }w_{2}/\partial y_{i}))_{1\leq i\leq N}$, 
$\widehat{\func{div}}\mathbf{w}=\func{div}w_{0}+\func{div}_{\omega }w_{1}+%
\widehat{\func{div}}_{y}w_{2}$ with $\func{div}_{\omega }w_{1}=\sum_{i=1}^{N}%
\overline{D}_{i,\omega }w_{1}$, $\widehat{\func{div}}_{y}w_{2}=\mathcal{G}%
_{1}(\overline{\func{div}}_{y}w_{2})$, $\overline{\func{div}}%
_{y}w_{2}=\sum_{i=1}^{N}\overline{\partial }w_{2}/\partial y_{i}$, $\mathcal{%
G}_{1}$ being the isometric isomorphism of $\mathcal{B}_{A}^{p}$ onto $%
L^{p}(\Delta (A))$.
\end{theorem}

\begin{proof}
First of all, it is evident that due to the properties of the mapping $a$ we
have 
\begin{equation}
c_{0}\left\| u_{\varepsilon }(\cdot ,\omega )\right\| _{W_{0}^{1,p}(\Omega
)}^{p-1}\leq \left\| \func{div}b^{\varepsilon }(\cdot ,\omega )\right\|
_{W^{-1,p^{\prime }}(\Omega )}.\;\;  \label{5.14}
\end{equation}%
Thus, rising the two members of the inequality (\ref{5.14}) to the power $%
p^{\prime }$ and integrating the resulting inequality over $\Omega $ we get 
\begin{equation*}
c_{0}^{p^{\prime }}\left\| u_{\varepsilon }\right\|
_{W_{0}^{1,p}(Q;L^{p}(\Omega ))}^{p}\leq \left\| \func{div}b^{\varepsilon
}\right\| _{L^{p^{\prime }}(\Omega ;W^{-1,p^{\prime }}(Q))}^{p^{\prime }}.
\end{equation*}%
Hence. as the sequence $(\func{div}b^{\varepsilon })$ is bounded in $%
L^{p^{\prime }}(\Omega ;W^{-1,p^{\prime }}(Q))$, it results that the
sequence $(u_{\varepsilon })_{\varepsilon >0}$ is bounded in $%
W_{0}^{1,p}(Q;L^{p}(\Omega ))$. Therefore the sequence $(a^{\varepsilon
}(\cdot ,Du_{\varepsilon }))_{\varepsilon >0}$ is bounded in $L^{p^{\prime
}}(Q\times \Omega )^{N}$. Thus, given a fundamental sequence $E$, Theorem %
\ref{t3.3} guarantees the existence of a subsequence $E^{\prime }$ extracted
from $E$ and a triplet $\mathbf{u}=(u_{0},u_{1},u_{2})\in \mathbb{F}%
_{0}^{1,p}$ such that (\ref{5.11})-(\ref{5.12}) hold when $E^{\prime }\ni
\varepsilon \rightarrow 0$. The next part of the proof is to show that $%
\mathbf{u}$ solves equation (\ref{5.13}). For that purpose, let $\Phi =(\psi
_{0},I_{p}(\psi _{1}),(J_{1}\circ \varrho )(\psi _{2}))\in \mathcal{F}%
_{0}^{\infty }$ and define $\Phi _{\varepsilon }$ ($\varepsilon >0$) as in (%
\ref{5.10}). Then $\Phi _{\varepsilon }(\cdot ,\omega )\in \mathcal{C}%
_{0}^{\infty }(Q)$ and further, in view of [part (iv) of] (\ref{5.3}) we
have 
\begin{equation*}
0\leq \int_{Q\times \Omega }(a^{\varepsilon }(\cdot ,Du_{\varepsilon
})-a^{\varepsilon }(\cdot ,D\Phi _{\varepsilon }))\cdot (Du_{\varepsilon
}-D\Phi _{\varepsilon })dxd\mu ,
\end{equation*}%
or, 
\begin{equation}
0\leq -\int_{Q\times \Omega }a^{\varepsilon }(\cdot ,D\Phi _{\varepsilon
})\cdot (Du_{\varepsilon }-D\Phi _{\varepsilon })dxd\mu +\int_{Q\times
\Omega }b^{\varepsilon }\func{div}(u_{\varepsilon }-\Phi _{\varepsilon
})dxd\mu .  \label{5.15}
\end{equation}%
But, 
\begin{eqnarray*}
\int_{Q\times \Omega }b^{\varepsilon }\func{div}(\text{\thinspace }%
u_{\varepsilon }-\Phi _{\varepsilon })dxd\mu &=&\int_{Q\times \Omega
}b^{\varepsilon }\func{div}u_{\varepsilon }dxd\mu -\int_{Q\times \Omega
}b^{\varepsilon }\func{div}\psi _{0}dxd\mu \\
&&+\int_{Q\times \Omega }b^{\varepsilon }[(\text{div}_{\omega }\text{%
\thinspace }\psi _{1})^{\varepsilon }+(\text{div}_{y}\text{\thinspace }\psi
_{2})^{\varepsilon }]dxd\mu \\
&&-\int_{Q\times \Omega }b^{\varepsilon }[\varepsilon _{1}(\func{div}\psi
_{1})^{\varepsilon }+\varepsilon _{2}(\func{div}\psi _{2})^{\varepsilon
})]dxd\mu \\
&&\;\;\;\;\;+\int_{Q\times \Omega }\frac{\varepsilon _{2}}{\varepsilon _{1}}(%
\text{div}_{\omega }\text{\thinspace }\psi _{2})^{\varepsilon }dxd\mu ,
\end{eqnarray*}%
and, as $E^{\prime }\ni \varepsilon \rightarrow 0$ we have by using (\ref%
{5.12}), 
\begin{equation}
\int_{Q\times \Omega }b^{\varepsilon }\func{div}(u_{\varepsilon }-\Phi
_{\varepsilon })dxd\mu \rightarrow \iint_{Q\times \Omega \times \Delta (A)}%
\widehat{b}\,\widehat{\func{div}}(\mathbf{u}-\Phi )dxd\mu d\beta .
\label{5.17}
\end{equation}%
Therefore passing to the limit in (\ref{5.15}) using (\ref{5.11})-(\ref{5.12}%
), the above convergence result (\ref{5.17}) together with Proposition \ref%
{p5.3}, we get 
\begin{equation}
0\leq -\iint_{Q\times \Omega \times \Delta (A)}\widehat{a}(\cdot ,\mathbb{D}%
\Phi )\cdot \mathbb{D}(\mathbf{u}-\Phi )dxd\mu d\beta +\iint_{Q\times \Omega
\times \Delta (A)}\widehat{b}\,\widehat{\func{div}}(\mathbf{u}-\Phi )dxd\mu
d\beta ,  \label{5.18}
\end{equation}%
$\Phi $ ranging over $\mathcal{F}_{0}^{\infty }$, and hence over $\mathbb{F}%
_{0}^{1,p}$ too (by a density argument). Taking in (\ref{5.18}) the
particular functions $\Phi =\mathbf{u}-t\mathbf{v}$ with $t>0$ and $\mathbf{v%
}=(v_{0},v_{1},v_{2})\in \mathbb{F}_{0}^{1,p}$, then dividing both sides of
the resulting inequality by $t$, and finally letting $t\rightarrow 0$, we
get (\ref{5.13}).

The next point to check is to show that $\mathbf{u}$ is unique. We begin by
showing that $u_{2}$ is unique. For that, we take in (\ref{5.13}) the
function $\mathbf{v}=(v_{0},v_{1},v_{2})$ with $v_{0}=0$ and $v_{1}=0$; then
for each fixed $(x,\omega )\in Q\times \Omega $, $u_{2}(x,\omega ,\cdot )$
is solution to 
\begin{equation}
\begin{array}{l}
\int_{\Delta (A)}\widehat{a}(x,\omega ,s,Du_{0}(x,\omega )+\overline{D}%
_{\omega }u_{1}(x,\omega )+\partial \widehat{u}_{2}(x,\omega ,s))\cdot
\partial \widehat{w}d\beta \\ 
\;\;\;=\int_{\Delta (A)}\widehat{b}(x,\omega ,s)\,\widehat{\overline{\text{%
div}}_{y}\,w}d\beta \text{\ for all }w\in \mathcal{B}_{\#A}^{1,p}.%
\end{array}
\label{5.20}
\end{equation}%
So, let $(x,\omega ,\xi )\in Q\times \Omega \times \mathbb{R}^{N}$ be freely
fixed, and let $\zeta =\zeta (x,\omega ,\xi )\in \mathcal{B}_{\#A}^{1,p}$ be
defined by the cell equation 
\begin{equation}
\int_{\Delta (A)}\widehat{a}(x,\omega ,\cdot ,\xi +\partial \widehat{\zeta }%
)\cdot \partial \widehat{w}d\beta =\int_{\Delta (A)}\widehat{b}(x,\omega
,\cdot )\widehat{\overline{\text{div}}_{y}\,w}d\beta \text{\ for all }w\in 
\mathcal{B}_{\#A}^{1,p}\text{.}  \label{5.21}
\end{equation}%
Since the linear functional $w\mapsto \int_{\Delta (A)}\widehat{b}(x,\omega
,\cdot )\widehat{\overline{\func{div}}_{y}\,w}d\beta $ is continuous on $%
\mathcal{B}_{\#A}^{1,p}$, it follows by \cite[Chap. 2]{Lions} that equation (%
\ref{5.21}) admits at least a solution. But this solution is unique; indeed
if $\zeta _{1}$ and $\zeta _{2}$ are solutions to (\ref{5.21}), then we have 
\begin{equation*}
\begin{array}{l}
\int_{\Delta (A)}(\widehat{a}(x,\omega ,\cdot ,\xi +\partial \widehat{\zeta }%
_{1})-\widehat{a}(x,\omega ,\cdot ,\xi +\partial \widehat{\zeta }_{2}))\cdot
\partial \widehat{w}d\beta =0 \\ 
\text{for all }w\in \mathcal{B}_{\#A}^{1,p}.%
\end{array}%
\end{equation*}%
Taking in particular $w=\zeta _{1}-\zeta _{2}$ in the above equation it
follows by [part (iv) of] (\ref{5.3}) that $\partial \widehat{\zeta }%
_{1}=\partial \widehat{\zeta }_{2}$ and hence $\overline{D}_{y}\zeta _{1}=%
\overline{D}_{y}\zeta _{2}$. We deduce that $\zeta _{1}=\zeta _{2}$ since
they belong to $\mathcal{B}_{\#A}^{1,p}$. Now, taking in (\ref{5.21}) $\xi
=Du_{0}(x)+\overline{D}_{\omega }u_{1}(x,\omega )$, and comparing the
resulting equation with (\ref{5.20}), we get by the uniqueness of the
solution of (\ref{5.21}) that $u_{2}(x,\omega ,\cdot )=\zeta (x,\omega
,Du_{0}(x,\omega )+\overline{D}_{\omega }u_{1}(x,\omega ))$ for a.e. $%
(x,\omega )\in Q\times \Omega $. This shows the uniqueness of $u_{2}$. The
same process shows the uniqueness of $u_{1}$ and of $u_{0}$. We conclude
that $\mathbf{u}$ is unique so that the convergence results (\ref{5.11}) and
(\ref{5.12}) hold for the whole sequence $\varepsilon $ as expected. This
completes the proof of the theorem.
\end{proof}

One can work out some homogenization problems related to problem (\ref{5.4}%
), (\ref{5.7}) and (\ref{5.8}). In particular one can solve:\bigskip

\begin{itemize}
\item[\textbf{(P)}$_{1}$] \textbf{The coupled stochastic-periodic
homogenization problem} stated as follows: For each fixed $(x,\omega
,\lambda )\in Q\times \Omega \times \mathbb{R}^{N}$, the functions $y\mapsto
a(x,\omega ,y,\lambda )$ and $y\mapsto b(x,\omega ,y)$, are $Y$-periodic
where $Y=(0,1)^{N}$. Here we get the homogenization of (\ref{5.4}) with $A=%
\mathcal{C}_{\text{per}}(Y)$.\bigskip

\item[\textbf{(P)}$_{2}$] \textbf{The coupled stochastic-almost periodic
homogenization problem} stated as follows: 
\begin{equation*}
\begin{array}{l}
a(x,\omega ,\cdot ,\lambda )\in (AP(\mathbb{R}^{N}))^{N}\text{ for any }%
(x,\omega ,\lambda )\in \overline{Q}\times \Omega \times \mathbb{R}^{N}; \\ 
b(x,\omega ,\cdot )\in AP(\mathbb{R}^{N})\text{ for a.e. }(x,\omega )\in
Q\times \Omega%
\end{array}%
\end{equation*}%
where here, $AP(\mathbb{R}^{N})$ is the algebra of all Bohr almost periodic
complex functions \cite{8} defined as the algebra of functions on $\mathbb{R}%
^{N}$ that are uniformly approximated by finite linear combinations of
functions in the set $\{\gamma _{k}:k\in \mathbb{R}^{N}\}$ with $\gamma
_{k}(y)=\exp (2i\pi k\cdot y)$ ($y\in \mathbb{R}^{N}$). It is known that $AP(%
\mathbb{R}^{N})$ satisfies assumptions of Theorem \ref{t3.3} (see \cite{CMP}%
). We are led to the homogenization of (\ref{5.4}) with $A=AP(\mathbb{R}%
^{N}) $.\bigskip

\item[\textbf{(P)}$_{3}$] \textbf{The coupled stochastic-weakly almost
periodic homogenization problem I}: 
\begin{equation*}
\begin{array}{l}
a(x,\omega ,\cdot ,\lambda )\in (WAP(\mathbb{R}^{N}))^{N}\text{ for any }%
(x,\omega ,\lambda )\in \overline{Q}\times \Omega \times \mathbb{R}^{N}; \\ 
b(x,\omega ,\cdot )\in AP(\mathbb{R}^{N})\text{ for a.e. }(x,\omega )\in
Q\times \Omega%
\end{array}%
\end{equation*}%
where $WAP(\mathbb{R}^{N})$ is the algebra of weakly almost periodic
functions on $\mathbb{R}^{N}$ \cite{17, CMP}. It is known \cite{CMP} that $%
WAP(\mathbb{R}^{N})$ satisfies hypotheses of Theorem \ref{t3.3} with
moreover, $AP(\mathbb{R}^{N})\subset WAP(\mathbb{R}^{N})$. This leads to the
homogenization of (\ref{5.4}) with $A=WAP(\mathbb{R}^{N})$.\bigskip

\item[\textbf{(P)}$_{4}$] \textbf{The coupled stochastic-weakly almost
periodic homogenization problem II}: 
\begin{equation*}
\begin{array}{l}
a(x,\omega ,\cdot ,\lambda )\in (AP(\mathbb{R}^{N}))^{N}\text{ for any }%
(x,\omega ,\lambda )\in \overline{Q}\times \Omega \times \mathbb{R}^{N}; \\ 
b(x,\omega ,\cdot )\in WAP(\mathbb{R}^{N})\text{ for a.e. }(x,\omega )\in
Q\times \Omega%
\end{array}%
\end{equation*}%
which yields the homogenization of (\ref{5.4}) with $A=WAP(\mathbb{R}^{N})$%
.\bigskip

\item[\textbf{(P)}$_{5}$] \textbf{The fully coupled stochastic-weakly almost
periodic homogenization problem III}: 
\begin{equation*}
\begin{array}{l}
a(x,\omega ,\cdot ,\lambda )\in (WAP(\mathbb{R}^{N}))^{N}\text{ for any }%
(x,\omega ,\lambda )\in \overline{Q}\times \Omega \times \mathbb{R}^{N}; \\ 
b(x,\omega ,\cdot )\in WAP(\mathbb{R}^{N})\text{ for a.e. }(x,\omega )\in
Q\times \Omega .%
\end{array}%
\end{equation*}%
Here the suitable $H$-supralgebra is $A=WAP(\mathbb{R}^{N})$.\bigskip
\end{itemize}

The same remark as the one made at the end of the preceding section is also
valid in this case, namely, the results of this section apply in all the
environments: the deterministic one and the stochastic one as well. This is
also true for the reiterated deterministic framework as seen in the
preceding section. This therefore extends all the results of the paper \cite%
{CMP} since in the case when $\Omega =\Delta (A_{z})$ we do not need to make
any ergodicity assumption on the algebra $A_{z}$ as it was the case in \cite%
{CMP}.

\section{Application to the homogenization of a model of rotating fluids}

Throughout this section, all the vector spaces are assumed to be real vector
spaces, and the scalar functions are assumed real valued.

\subsection{Introduction and preliminary results}

It is well known that the flows of\ commonly encountered Newtonian fluids
are modeled by the Navier-Stokes equations. These flows are sometimes
laminar, sometimes turbulent. Unfortunately, in reality, the flows of fluids
are almost always turbulent. Thus starting from two identical situations,
the flow may evolve very differently. This explains its dual nature of being
both deterministic and unpredictable (random).

In this section, our goal is not to establish the conditions for the
prediction of \ the turbulence, but to describe the asymptotic behavior of a
model of turbulence. More precisely we study the asymptotic behavior, as $%
0<\varepsilon \rightarrow 0$, of the following three dimensional Stokes
equation 
\begin{equation}
\begin{array}{l}
P^{\varepsilon }\mathbf{u}_{\varepsilon }+\mathbf{h}^{\varepsilon }\times 
\mathbf{u}_{\varepsilon }+\func{grad}p_{\varepsilon }=\mathbf{f}\text{\ in }Q
\\ 
\;\;\;\;\;\;\;\;\;\;\;\;\;\;\;\;\;\;\;\;\ \ \ \ \ \ \ \;\ \func{div}\mathbf{u%
}_{\varepsilon }=0\text{\ in }Q \\ 
\;\;\;\;\;\;\;\;\;\;\;\;\;\;\;\;\;\;\;\;\;\ \ \ \ \ \ \ \ \ \ \;\;\;\mathbf{u%
}_{\varepsilon }=0\text{\ on }\partial Q.%
\end{array}
\label{6.1}
\end{equation}%
Let us make precise the data in (\ref{6.1}). Let $Q$ be a smooth bounded
open set in $\mathbb{R}_{x}^{N}$ ($N=3$); in $Q$ we consider the partial
differential operator 
\begin{equation*}
P^{\varepsilon }=-\overset{N}{\underset{i,j=1}{\sum }}\frac{\partial }{%
\partial x_{i}}\left( a_{ij}\left( x,T(x/\varepsilon _{1})\omega
,x/\varepsilon _{2}\right) \frac{\partial }{\partial x_{j}}\right)
\end{equation*}%
where $T$ is an $N$-dimensional dynamical system acting on the probability
space $(\Omega ,\mathcal{M},\mu )$, the functions $a_{ij}\in \mathcal{C}(%
\overline{Q};L^{\infty }(\Omega ;\mathcal{B}(\mathbb{R}_{y}^{N})))$ ($1\leq
i,j\leq N$) satisfy the following assumptions: 
\begin{equation}
a_{ij}=a_{ji}  \label{6.2}
\end{equation}%
and there exists a constant $\alpha >0$ such that 
\begin{equation}
\begin{array}{l}
\overset{N}{\underset{i,j=1}{\sum }}a_{ij}(x,\omega ,y)\lambda _{i}\lambda
_{j}\geq \alpha \left\vert \lambda \right\vert ^{2}\text{ for all }\lambda
=(\lambda _{i})\in \mathbb{R}^{N}\text{, all }x\in \overline{Q} \\ 
\text{and for almost all }(\omega ,y)\in \Omega \times \mathbb{R}^{N}\text{.}%
\end{array}
\label{6.3}
\end{equation}%
The operator $P^{\varepsilon }$ defined above is assumed to act on vector
functions as follows: for $\mathbf{u}=(u^{i})_{1\leq i\leq N}\in
H^{1}(Q)^{N}=(W^{1,2}(Q))^{N}$ we have $P^{\varepsilon }\mathbf{u}%
=(P^{\varepsilon }u^{i})_{1\leq i\leq N}$. The function $\mathbf{h}%
^{\varepsilon }$ is defined by $\mathbf{h}^{\varepsilon }(x,\omega )=\mathbf{%
h}(T(x/\varepsilon _{1})\omega ,x/\varepsilon _{2})$ for $(x,\omega )\in
Q\times \Omega $, where $\mathbf{h}=(h_{i})\in L^{\infty }(\Omega ;\mathcal{B%
}(\mathbb{R}_{y}^{N}))^{N}$. Likewise, for two vector functions $\mathbf{u}%
=(u^{i})$ and $\mathbf{v}=(v^{i})$ both in $L^{2}(Q)^{N}$, $\mathbf{u}\times 
\mathbf{v}$ denotes the exterior product of $\mathbf{u}$ and $\mathbf{v}$
defined to be the vector $\mathbf{w}=(w^{i})$ with 
\begin{equation}
w^{i}=\sum_{j,k=1}^{N=3}\varepsilon _{ijk}u^{j}v^{k}\;\;\;\;(1\leq i\leq N=3)
\label{6.1'}
\end{equation}%
where $\varepsilon _{ijk}$ is the totally antisymmetric tensor defined as
follows: $\varepsilon _{iii}=0$ for $1\leq i\leq 3$, and $\varepsilon
_{123}=\varepsilon _{231}=\varepsilon _{312}=1$ and $\varepsilon
_{321}=\varepsilon _{213}=\varepsilon _{132}=-1$; $\partial \mathbf{u}%
/\partial x_{j}$ stands for the vector $(\partial u^{1}/\partial
x_{j},\ldots ,\partial u^{N}/\partial x_{j})$. Finally, the function $%
\mathbf{f}$ lies in $L^{\infty }(\Omega ;H^{-1}(Q)^{N})=L^{\infty }(\Omega
;(W^{-1,2}(Q))^{N})$ and $\func{grad}p$ (for $p\in L^{2}(Q)$) designates the
gradient of $p$, sometimes denoted by $Dp$.

It is known that the problem (\ref{6.1}) (for each fixed $\varepsilon >0$
and for $\mu $-almost all $\omega \in \Omega $) uniquely determines a couple 
$(\mathbf{u}_{\varepsilon }(\cdot ,\omega ),p_{\varepsilon }(\cdot ,\omega
))\in H_{0}^{1}(Q)^{N}\times (L^{2}(Q)/\mathbb{R})$, which therefore yields
a unique couple $(\mathbf{u}_{\varepsilon },p_{\varepsilon })\in
H_{0}^{1}(Q;L^{2}(\Omega ))^{N}\times L^{2}(\Omega ;L^{2}(Q)/\mathbb{R})$.
Thus we have in hand a sequence $((\mathbf{u}_{\varepsilon },p_{\varepsilon
}))_{\varepsilon >0}$ and we aim at investigating its asymptotic behavior,
as $\varepsilon \rightarrow 0$, under suitable assumptions on $a_{ij}$ ($%
1\leq i,j\leq N$) and on $\mathbf{h}$. It is worth noting that there exist
many references on the homogenization of Stokes equations in the periodic
setting as well as in the stochastic setting.

We assume throughout this section that $A$ is an algebra with mean value on $%
\mathbb{R}_{y}^{N}$. In the study of the problem (\ref{6.1}) the following
issues arise:

\begin{enumerate}
\item To establish the conditions under which the solutions of (\ref{6.1})
converge as $\varepsilon \rightarrow 0$;

\item To determine the boundary value problem for the limit function.
\end{enumerate}

These issues will be addressed in the next subsection. Prior to this, we
introduce the following space: 
\begin{equation*}
\mathbb{H}_{0}^{1}(Q)=\{\mathbf{u}\in H_{0}^{1}(Q)^{N}:\func{div}\mathbf{u}%
=0\}.
\end{equation*}%
This is a Hilbert space under the Hilbertian norm of $H_{0}^{1}(Q)^{N}$
defined by 
\begin{equation*}
\left\| \mathbf{v}\right\| _{H_{0}^{1}(Q)^{N}}=\left(
\sum_{k=1}^{N}\int_{Q}\left| \nabla v^{k}\right| ^{2}dx\right) ^{1/2},\ \ 
\mathbf{v}=(v^{k})\in H_{0}^{1}(Q)^{N}
\end{equation*}%
with $\nabla v^{k}$ denoting the gradient of $v^{k}\in H_{0}^{1}(Q)$. Next
we introduce the bilinear form $a^{\varepsilon }(\omega ;\cdot ,\cdot )$
defined as follows:%
\begin{equation*}
a^{\varepsilon }(\omega ;\mathbf{u},\mathbf{v})=\sum_{i,j,k=1}^{N}%
\int_{Q}a_{ij}^{\varepsilon }(x,\omega )\frac{\partial u^{k}}{\partial x_{j}}%
\frac{\partial v^{k}}{\partial x_{i}}dx\text{\ \ }(\mathbf{u}=(u^{k}),\ 
\mathbf{v}=(v^{k})\in H_{0}^{1}(Q)^{N}).
\end{equation*}%
One easily sees by (\ref{6.2})-(\ref{6.3}) that $a^{\varepsilon }(\omega
;\cdot ,\cdot )$ is symmetric and satisfies the coercivity assumption 
\begin{equation}
a^{\varepsilon }(\omega ;\mathbf{v},\mathbf{v})\geq \alpha \left\| \mathbf{v}%
\right\| _{H_{0}^{1}(Q)^{N}}^{2}\ \ (\mathbf{v}\in
H_{0}^{1}(Q)^{N},0<\varepsilon <1).  \label{6.4}
\end{equation}%
Moreover we have $\left| a^{\varepsilon }(\omega ;\mathbf{u},\mathbf{v}%
)\right| \leq c\left\| \mathbf{u}\right\| _{H_{0}^{1}(Q)^{N}}\left\| \mathbf{%
v}\right\| _{H_{0}^{1}(Q)^{N}}$ for every $\mathbf{u},\mathbf{v}\in
H_{0}^{1}(Q)^{N}$ and $0<\varepsilon <1$, where $c$ is a positive constant
independent of $\varepsilon $ and of $\omega \in \Omega $.

In the sequel we will use the following notation: the stochastic divergence
operator $\func{div}_{\omega ,2}$ will be merely denoted by $\func{div}%
_{\omega }$. With this in mind, we will make use of the following spaces: 
\begin{equation*}
W_{\func{div}_{\omega }}^{1,2}(\Omega )=\{\mathbf{u}\in \mathcal{W}%
^{1,2}(\Omega )^{N}:\text{div}_{\omega }\mathbf{u}=0\}
\end{equation*}%
and 
\begin{equation*}
\mathcal{B}_{\func{div}_{y}}^{1,2}=\{\mathbf{u}\in (\mathcal{B}%
_{\#A}^{1,2})^{N}:\overline{\func{div}}_{y}\mathbf{u}=0\}
\end{equation*}%
where $\overline{\func{div}}_{y}\mathbf{u}=\sum_{i=1}^{N}\overline{\partial }%
u^{i}/\partial y_{i}$ and div$_{\omega }\mathbf{u}=\sum_{i=1}^{N}\overline{D}%
_{i,\omega }u^{i}$, and of their smooth counterparts 
\begin{equation*}
\mathcal{W}_{\func{div}_{\omega }}^{\infty }(\Omega )=\{\mathbf{u}\in 
\mathcal{C}^{\infty }(\Omega )^{N}:\text{div}_{\omega }\mathbf{u}=0\}
\end{equation*}%
and 
\begin{equation*}
A_{\func{div}_{y}}^{\infty }/\mathbb{R}=\{\mathbf{u}\in (A^{\infty })^{N}:M(%
\mathbf{u})=0\text{ and div}_{y}\mathbf{u}=0\}.
\end{equation*}%
The following result holds.

\begin{lemma}
\label{l6.1}The space $I_{2}^{N}(\mathcal{W}_{\func{div}_{\omega }}^{\infty
}(\Omega ))$ (resp. $(J_{1}\circ \varrho )^{N}(A_{\func{div}_{y}}^{\infty }/%
\mathbb{R})$) is dense in $W_{\func{div}_{\omega }}^{1,2}(\Omega )$ (resp. $%
\mathcal{B}_{\func{div}_{y}}^{1,2}$) where, for $\mathbf{u}=(u^{i})_{i}\in
A_{\func{div}_{y}}^{\infty }/\mathbb{R}$, $(J_{1}\circ \varrho )^{N}(\mathbf{%
u})=(J_{1}(\varrho (u^{i})))_{i}$, and, for $\mathbf{u}=(u^{i})_{i}\in 
\mathcal{W}_{\func{div}_{\omega }}^{\infty }(\Omega )$, $I_{2}^{N}(\mathbf{u}%
)=(I_{2}(u^{i}))_{i}$.
\end{lemma}

\begin{proof}
As regard the denseness of $\mathcal{W}_{\func{div}_{\omega }}^{\infty
}(\Omega )$ in $W_{\func{div}_{\omega }}^{1,2}(\Omega )$, this follows in
the same way as the proof of \cite[Lemma 2.3]{Wright1}. Concerning the next
part, as it was seen in Section 2, when $\Omega =\Delta (A)$, the space $%
\mathcal{C}^{\infty }(\Omega )$ is just replaced by the space $\mathcal{G}%
_{1}(\varrho (A^{\infty }))=\mathcal{G}(A^{\infty })=\mathcal{D}(\Delta (A))$%
. Moreover the algebra $A$ being ergodic, the invariant functions (for the
dynamical system induced by the translations on $\Delta (A)$) consist of
constants. Therefore we have $\mathcal{W}^{1,2}(\Delta (A))=\overline{%
\mathcal{G}}_{1}(\mathcal{B}_{\#A}^{1,2})=W_{\#}^{1,2}(\Delta (A))$ (see (%
\ref{2.5'}) and (\ref{2.6'}) for the properties of $\overline{\mathcal{G}}%
_{1}$). Let us recall that $\mathcal{W}^{1,2}(\Delta (A))$ is the completion
of $\mathcal{C}^{\infty }(\Delta (A))\equiv \mathcal{D}(\Delta (A))$ with
respect to the seminorm (\ref{2.0}) (see Section 2), which is also the
completion with respect to the same seminorm of $\mathcal{C}^{\infty
}(\Delta (A))/\mathbb{R}=\{u\in \mathcal{D}(\Delta (A)):\int_{\Delta
(A)}ud\beta =0\}=\mathcal{D}(\Delta (A))/\mathbb{R}$ (that is, $%
W_{\#}^{1,2}(\Delta (A))$) since any $u\in \mathcal{C}^{\infty }(\Delta (A))$
invariant (for the dynamical system induced on $\Delta (A)$ by the
translations on $\mathbb{R}^{N}$) is constant. Hence, using once again %
\cite[Lemma 2.3]{Wright1} we get the last part, and the lemma is proved.
\end{proof}

Now, let $\mathcal{D}_{\func{div}}(Q)=\{\mathbf{\varphi }\in \mathcal{C}%
_{0}^{\infty }(Q)^{N}:\func{div}\mathbf{\varphi }=0\}$. It is known \cite%
{Lions, Temam} that $\mathcal{D}_{\func{div}}(Q)$ is dense in $\mathbb{H}%
_{0}^{1}(Q)$. Next set 
\begin{equation*}
\mathbb{F}_{0}^{1}=\mathbb{H}_{0}^{1}(Q;I_{nv}^{2}(\Omega ))\times
L^{2}(Q;W_{\func{div}_{\omega }}^{1,2}(\Omega ))\times L^{2}(Q\times \Omega ;%
\mathcal{B}_{\func{div}_{y}}^{1,2})
\end{equation*}%
and 
\begin{eqnarray*}
\mathcal{F}_{0}^{\infty } &=&[\mathcal{D}_{\func{div}}(Q;I_{nv}^{2}(\Omega
))]\times \lbrack \mathcal{C}_{0}^{\infty }(Q)\otimes I_{2}^{N}(\mathcal{W}_{%
\func{div}_{\omega }}^{\infty }(\Omega ))]\times \lbrack \mathcal{C}%
_{0}^{\infty }(Q)\otimes \mathcal{C}^{\infty }(\Omega )\otimes \\
&&(J_{1}\circ \varrho )^{N}(A_{\func{div}_{y}}^{\infty }/\mathbb{R})]
\end{eqnarray*}%
where: $\mathcal{D}_{\func{div}}(Q;I_{nv}^{2}(\Omega ))$ is defined to be
the space of those $\mathbf{u}\in (\mathcal{C}_{0}^{\infty }(Q)\otimes
I_{nv}^{2}(\Omega ))^{N}$ such that $\func{div}\mathbf{u}=0$ (and the other
members of the Cartesian product in $\mathcal{F}_{0}^{\infty }$ are defined
as usual) and $\mathbb{H}_{0}^{1}(Q;I_{nv}^{2}(\Omega ))=\{\mathbf{u}\in
H_{0}^{1}(Q;I_{nv}^{2}(\Omega ))^{N}:\func{div}\mathbf{u}=0\}$. Then thanks
to Lemma \ref{l6.1}, the space $\mathcal{F}_{0}^{\infty }$ is dense in $%
\mathbb{F}_{0}^{1}$.

\subsection{Homogenization result}

Our goal in this subsection is the study of the asymptotic behavior of $(%
\mathbf{u}_{\varepsilon })_{\varepsilon >0}$ (the sequence of solutions to (%
\ref{6.1})) under the following assumptions: 
\begin{equation}
a_{ij}(x,\omega ,\cdot )\in A\text{ for all }(x,\omega )\in Q\times \Omega 
\text{, }1\leq i,j\leq N;  \label{6.8}
\end{equation}%
\begin{equation}
\mathbf{h}\in L^{\infty }(\Omega ;A)^{N}.\text{\ \ \ \ \ \ \ \ \ \ \ \ \ \ \
\ \ \ \ \ \ \ \ \ \ \ \ \ \ \ \ \ \ \ \ \ \ \ \ \ \ \ \ \ \ \ \ \ \ \ \ \ \
\ \ \ \ \ \ \ }  \label{6.9}
\end{equation}

We are now able to state and prove the homogenization result of this section.

\begin{theorem}
\label{t6.1}Assume \emph{(\ref{6.8})-(\ref{6.9})} hold. For each $%
0<\varepsilon <1$ and for a.e. $\omega \in \Omega $ let $\mathbf{u}%
_{\varepsilon }(\cdot ,\omega )=(u_{\varepsilon }^{k}(\cdot ,\omega ))\in 
\mathbb{H}_{0}^{1}(Q)$ be defined by \emph{(\ref{6.1})}. Then as $%
\varepsilon \rightarrow 0$, 
\begin{equation}
\mathbf{u}_{\varepsilon }\rightarrow \mathbf{u}_{0}\text{ stoch. in }%
L^{2}(Q\times \Omega )^{N}\text{-weak}  \label{6.13}
\end{equation}%
and 
\begin{equation}
\frac{\partial u_{\varepsilon }^{k}}{\partial x_{j}}\rightarrow \frac{%
\partial u_{0}^{k}}{\partial x_{j}}+\overline{D}_{j,\omega }u_{1}^{k}+\frac{%
\overline{\partial }u_{2}^{k}}{\partial y_{j}}\text{ stoch. in }%
L^{2}(Q\times \Omega )\text{-weak }\Sigma \text{ }(1\leq j,k\leq N)
\label{6.14}
\end{equation}%
where $\mathbf{u}=(\mathbf{u}_{0},\mathbf{u}_{1},\mathbf{u}_{2})\in \mathbb{F%
}_{0}^{1}$ (with $\mathbf{u}_{i}=(u_{i}^{k})_{1\leq k\leq N}$, $0\leq i\leq
2 $) is the unique solution to the following variational problem: 
\begin{equation}
\left\{ 
\begin{array}{l}
a(\mathbf{u},\mathbf{v})+\iint_{Q\times \Omega }(\widetilde{\mathbf{h}}%
\times \mathbf{u}_{0})\cdot \mathbf{v}_{0}dxd\mu =\left\langle \mathbf{f},%
\mathbf{v}_{0}\right\rangle \\ 
\text{for all }\mathbf{v}=(\mathbf{v}_{0},\mathbf{v}_{1},\mathbf{v}_{2})\in 
\mathbb{F}_{0}^{1}%
\end{array}%
\right.  \label{6.15}
\end{equation}%
with: 
\begin{eqnarray*}
a(\mathbf{u},\mathbf{v}) &=&\sum_{i,j,k=1}^{N}\iint_{Q\times \Omega \times
\Delta (A)}\widehat{a}_{ij}(x,\omega ,s)\left( \frac{\partial u_{0}^{k}}{%
\partial x_{j}}+\overline{D}_{j,\omega }u_{1}^{k}+\partial _{j}\widehat{%
u_{2}^{k}}\right) \\
&&\;\;\;\;\ \ \ \ \ \ \ \ \times \left( \frac{\partial v_{0}^{k}}{\partial
x_{i}}+\overline{D}_{i,\omega }v_{1}^{k}+\partial _{i}\widehat{v_{2}^{k}}%
\right) dxd\mu d\beta ;
\end{eqnarray*}%
\begin{equation*}
\widetilde{\mathbf{h}}(\omega )=\int_{\Delta (A)}\widehat{\mathbf{h}}(\omega
,s)d\beta ;\ \ \ \ \ \ \ \ \ \ \ \ \ \ \ \ \ \ \ \ \ \ \ \ \ \ \ \ \ \ \ \ \ 
\end{equation*}%
\begin{equation*}
\left\langle \mathbf{f},\mathbf{v}_{0}\right\rangle =\int_{\Omega }\left( 
\mathbf{f}(\cdot ,\omega ),\mathbf{v}_{0}(\cdot ,\omega )\right)
_{H^{-1}(Q)^{N},H_{0}^{1}(Q)^{N}}d\mu
\end{equation*}%
and%
\begin{equation*}
\partial _{j}\widehat{u_{2}^{k}}=\mathcal{G}_{1}(\overline{\partial }%
u_{2}^{k}/\partial y_{j})\text{ (and a same definition for }\partial _{i}%
\widehat{v_{2}^{k}}\text{).}
\end{equation*}
\end{theorem}

\begin{proof}
We have (for each $0<\varepsilon <1$ and for $\mu $-a.e. $\omega \in \Omega $%
) 
\begin{equation}
\begin{array}{l}
a^{\varepsilon }(\omega ;\mathbf{u}_{\varepsilon }(\cdot ,\omega ),\mathbf{v}%
)+\int_{Q}(\mathbf{h}^{\varepsilon }(\cdot ,\omega )\times \mathbf{u}%
_{\varepsilon }(\cdot ,\omega ))\cdot \mathbf{v}dx-\int_{Q}p_{\varepsilon
}(\cdot ,\omega )\func{div}\mathbf{v}dx \\ 
=\left( \mathbf{f}(\cdot ,\omega ),\mathbf{v}\right) \text{ for all }\mathbf{%
v}\in H_{0}^{1}(Q)^{N}%
\end{array}
\label{6.16}
\end{equation}%
where $(\mathbf{f}(\cdot ,\omega ),\mathbf{v})\equiv (\mathbf{f}(\cdot
,\omega ),\mathbf{v})_{H^{-1}(Q)^{N},H_{0}^{1}(Q)^{N}}$. Taking the
particular $\mathbf{v}=\mathbf{u}_{\varepsilon }(\cdot ,\omega )$ in (\ref%
{6.16}) and using the fact that $(\mathbf{h}^{\varepsilon }(\cdot ,\omega
)\times \mathbf{u}_{\varepsilon }(\cdot ,\omega ))\cdot \mathbf{u}%
_{\varepsilon }(\cdot ,\omega )=0$, we obtain immediately that the sequence $%
(\mathbf{u}_{\varepsilon }(\cdot ,\omega ))_{\varepsilon >0}$ is bounded in $%
H_{0}^{1}(Q)^{N}$ uniformly in $\omega $. On the other hand, it is an easy
task (using the above boundedness condition) to see that 
\begin{equation*}
\left\vert \left( \func{grad}p_{\varepsilon }(\cdot ,\omega ),\mathbf{v}%
\right) \right\vert \leq c\left\Vert \mathbf{v}\right\Vert
_{H_{0}^{1}(Q)^{N}}\text{ for all }\mathbf{v}\in H_{0}^{1}(Q)^{N}\text{,}
\end{equation*}%
$c$ being a positive constant independent of $\omega $ and of $\mathbf{v}\in
H_{0}^{1}(Q)^{N}$. Hence, the sequence $(\func{grad}p_{\varepsilon }(\cdot
,\omega ))_{0<\varepsilon <1}$ is bounded in $H^{-1}(Q)^{N}$ independently
of $\omega \in \Omega $. Therefore, using a well-known argument (see e.g. %
\cite[p. 15]{Temam}) we deduce that the sequence $(p_{\varepsilon }(\cdot
,\omega ))_{0<\varepsilon <1}$ is bounded in $L^{2}(Q)^{N}$ independently of 
$\omega $. So, given an arbitrary ordinary sequence $E$, Theorems \ref{t3.1}
and \ref{t3.3} give rise to a subsequence $E^{\prime }$ from $E$ and
functions $\mathbf{u}_{0}=(u_{0}^{k})\in H_{0}^{1}(Q;I_{nv}^{2}(\Omega
))^{N} $, $\mathbf{u}_{1}=(u_{1}^{k})\in L^{2}(Q;\mathcal{W}^{1,2}(\Omega
))^{N}$, $\mathbf{u}_{2}=(u_{2}^{k})\in L^{2}(Q\times \Omega ;\mathcal{B}%
_{\#A}^{1,2})^{N}$, $p_{0}\in L^{2}(Q\times \Omega ;\mathcal{B}_{A}^{2})$
such that, as $E^{\prime }\ni \varepsilon \rightarrow 0$ we have (\ref{6.13}%
)-(\ref{6.14}) and 
\begin{equation}
p_{\varepsilon }\rightarrow p_{0}\text{ stoch. in }L^{2}(Q\times \Omega )%
\text{-weak }\Sigma \text{.}  \label{6.17}
\end{equation}%
It is easy to see that, due to the equality $\func{div}\mathbf{u}%
_{\varepsilon }=0$, we have $\func{div}\mathbf{u}_{0}=0$, $\func{div}%
_{\omega }\mathbf{u}_{1}=0$ and $\overline{\func{div}}_{y}\mathbf{u}_{2}=0$.
Therefore $\mathbf{u}=(\mathbf{u}_{0},\mathbf{u}_{1},\mathbf{u}_{2})\in 
\mathbb{F}_{0}^{1}$. The next step is to show that $\mathbf{u}$ solves
equation (\ref{6.15}). To this end, let $\mathbf{\Phi }=(\Psi
_{0},I_{2}^{N}(\Psi _{1}),(J_{1}\circ \varrho )^{N}(\Psi _{2}))\in \mathcal{F%
}_{0}^{\infty }$; define $\mathbf{\Phi }_{\varepsilon }:=\Psi
_{0}+\varepsilon _{1}\Psi _{1}^{\varepsilon }+\varepsilon _{2}\Psi
_{2}^{\varepsilon }$, that is, $\Phi _{\varepsilon }(x,\omega )=\Psi
_{0}(x,\omega )+\varepsilon _{1}\Psi _{1}(x,T(x/\varepsilon _{1})\omega
)+\varepsilon _{2}\Psi _{2}(x,T(x/\varepsilon _{1})\omega ,x/\varepsilon
_{2})$ for $(x,\omega )\in Q\times \Omega $. We have, in view of (\ref{6.16}%
), 
\begin{equation}
\begin{array}{l}
\int_{\Omega }a^{\varepsilon }(\omega ;\mathbf{u}_{\varepsilon },\mathbf{%
\Phi }_{\varepsilon })d\mu +\int_{Q\times \Omega }(\mathbf{h}^{\varepsilon
}\times \mathbf{u}_{\varepsilon })\cdot \mathbf{\Phi }_{\varepsilon }dxd\mu
-\int_{Q\times \Omega }p_{\varepsilon }\func{div}\mathbf{\Phi }_{\varepsilon
}dxd\mu \\ 
\ \ \ \ =\int_{\Omega }(\mathbf{f}(\cdot ,\omega ),\mathbf{\Phi }%
_{\varepsilon }(\cdot ,\omega ))d\mu .%
\end{array}
\label{6.18}
\end{equation}%
We need to pass to the limit in (\ref{6.18}). Starting from the term $%
\int_{\Omega }a^{\varepsilon }(\omega ;\mathbf{u}_{\varepsilon },\mathbf{%
\Phi }_{\varepsilon })d\mu $, proceeding as in the proof of Theorem \ref%
{t4.1} we get 
\begin{equation*}
\int_{\Omega }a^{\varepsilon }(\omega ;\mathbf{u}_{\varepsilon },\mathbf{%
\Phi }_{\varepsilon })d\mu \rightarrow a(\mathbf{u},\mathbf{\Phi })\text{ as 
}E^{\prime }\ni \varepsilon \rightarrow 0\text{.}
\end{equation*}%
From the definition of $\mathbf{h}^{\varepsilon }\times \mathbf{u}%
_{\varepsilon }$, it readily follows from Proposition \ref{p3.4} (taking
there $h_{i}\psi _{0,j}\in \mathcal{K}(Q;L^{\infty }(\Omega ,A))$ as a test
function, where $\Psi _{0}=(\psi _{0,i})_{1\leq i\leq N}$) that, as $%
E^{\prime }\ni \varepsilon \rightarrow 0$, 
\begin{equation*}
\int_{Q\times \Omega }(\mathbf{h}^{\varepsilon }\times \mathbf{u}%
_{\varepsilon })\cdot \mathbf{\Phi }_{\varepsilon }dxd\mu \rightarrow
\int_{Q\times \Omega }(\widetilde{\mathbf{h}}\times \mathbf{u}_{0})\cdot
\Psi _{0}dxd\mu ,
\end{equation*}%
and due to (\ref{6.17}), as $E^{\prime }\ni \varepsilon \rightarrow 0$, 
\begin{equation*}
\int_{Q\times \Omega }p_{\varepsilon }\func{div}\mathbf{\Phi }_{\varepsilon
}dxd\mu \rightarrow \iint_{Q\times \Omega \times \Delta (A)}\widehat{p}_{0}(%
\func{div}\Psi _{0}+\text{div}_{\omega }\Psi _{1}+\widehat{\text{div}%
_{y}\Psi _{2}})dxd\mu d\beta ,
\end{equation*}%
which, with the fact that $\mathbf{\Phi }\in \mathcal{F}_{0}^{\infty }$
(which yields $\func{div}\Psi _{0}=0$, $\func{div}_{\omega }\Psi _{1}=0$ and 
$\func{div}_{y}\Psi _{2}=0$) gives 
\begin{equation*}
\int_{Q\times \Omega }p_{\varepsilon }\func{div}\mathbf{\Phi }_{\varepsilon
}dxd\mu \rightarrow 0\text{ when }E^{\prime }\ni \varepsilon \rightarrow 0.
\end{equation*}%
Moreover one obviously has $\int_{\Omega }(\mathbf{f}(\cdot ,\omega ),%
\mathbf{\Phi }_{\varepsilon }(\cdot ,\omega ))d\mu \rightarrow \int_{\Omega
}(\mathbf{f}(\cdot ,\omega ),\Psi _{0}(\cdot ,\omega ))d\mu $ when $%
E^{\prime }\ni \varepsilon \rightarrow 0$.

Finally, taking into account all the above facts, a passage to the limit in (%
\ref{6.18}) when $E^{\prime }\ni \varepsilon \rightarrow 0$ yields 
\begin{equation*}
a(\mathbf{u},\mathbf{\Phi })+\iint_{Q\times \Omega }(\widetilde{\mathbf{h}}%
\times \mathbf{u}_{0})\cdot \Psi _{0}dxd\mu =\left\langle \mathbf{f},\Psi
_{0}\right\rangle
\end{equation*}%
for all $\mathbf{\Phi }\in \mathcal{F}_{0}^{\infty }$. Using the continuity
of the linear form $\mathbf{v}_{0}\mapsto \iint_{Q\times \Omega }(\widetilde{%
\mathbf{h}}\times \mathbf{u}_{0})\cdot \mathbf{v}_{0}dxd\mu $ on $\mathbb{H}%
_{0}^{1}(Q;I_{nv}^{2}(\Omega ))$ (recall that $\widetilde{\mathbf{h}}\in
L^{\infty }(\Omega )^{N}$) associated to the density of $\mathcal{F}%
_{0}^{\infty }$ in $\mathbb{F}_{0}^{1}$, we are led at once to (\ref{6.15}).
Finally, from the equality $\iint_{Q\times \Omega }(\widetilde{\mathbf{h}}%
\times \mathbf{u}_{0})\cdot \mathbf{u}_{0}dxd\mu =0$ it classically follows
that the solution of (\ref{6.15}) is unique. Therefore (\ref{6.13})-(\ref%
{6.14}) hold for the whole sequence $\varepsilon >0$ as claimed.
\end{proof}

\subsection{Some applications of Theorem \ref{t6.1}}

We give in this subsection some concrete situations in which Theorem \ref%
{t6.1} is applicable. First of all, we recall that we will only need to
satisfy assumptions (\ref{6.8})-(\ref{6.9}). With this in mind, we see that
one can solve the following homogenization problems:

\begin{itemize}
\item[(P)$_{1}$] \textbf{The coupled stochastic-periodic homogenization
problem} stated as follows: For each fixed $1\leq i,j\leq N$ and for $\mu $%
-a.e. $\omega \in \Omega $ and a.e. $x\in Q$, the functions $y\mapsto
a_{ij}(x,\omega ,y)$ and $y\mapsto \mathbf{h}(\omega ,y)$ are $Y$-periodic,
where $Y=(0,1)^{N}$. Thus, we are led to the homogenization of (\ref{6.1})
under the above assumptions, but with $A=\mathcal{C}_{\text{per}}(Y)$%
.\bigskip

\item[(P)$_{2}$] \textbf{The coupled stochastic-almost periodic
homogenization problem} stated as follows: 
\begin{equation*}
a_{ij}(x,\omega ,\cdot )\in AP(\mathbb{R}^{N})\text{ and }\mathbf{h}(\omega
,\cdot )\in (AP(\mathbb{R}^{N}))^{N}.
\end{equation*}%
The homogenization of (\ref{6.1}) follows with $A=AP(\mathbb{R}^{N})$%
.\bigskip

\item[(P)$_{3}$] \textbf{The coupled stochastic-perturbed almost periodic
homogenization problem}: 
\begin{equation*}
a_{ij}(x,\omega ,\cdot )\in AP(\mathbb{R}^{N})+\mathcal{C}_{0}(\mathbb{R}%
^{N})\text{ and }\mathbf{h}(\omega ,\cdot )\in (\mathcal{C}_{\text{per}%
}(Y))^{N}
\end{equation*}%
where $\mathcal{C}_{0}(\mathbb{R}^{N})$ is the space of functions on $%
\mathbb{R}^{N}$ that vanish at infinity. It is a fact that $A=AP(\mathbb{R}%
^{N})+\mathcal{C}_{0}(\mathbb{R}^{N})$ is an ergodic $H$-supralgebra (called
the algebra of perturbed almost periodic functions) satisfying the
assumptions of Theorem \ref{t3.3}; see \cite{CMP}. Thus we get the
homogenization of (\ref{6.1}) with the above $A$.\bigskip

\item[(P)$_{4}$] \textbf{The coupled stochastic-weakly almost periodic
homogenization} \textbf{problem} stated either as 
\begin{equation*}
a_{ij}(x,\omega ,\cdot )\in WAP(\mathbb{R}^{N})\text{ and }\mathbf{h}(\omega
,\cdot )\in (AP(\mathbb{R}^{N}))^{N}
\end{equation*}%
or 
\begin{equation*}
a_{ij}(x,\omega ,\cdot )\in WAP(\mathbb{R}^{N})\text{ and }\mathbf{h}(\omega
,\cdot )\in (WAP(\mathbb{R}^{N}))^{N}.
\end{equation*}%
In each of the above cases we are led to the homogenization of (\ref{6.1})
with $A=WAP(\mathbb{R}^{N})$.\bigskip

\item[(P)$_{5}$] \textbf{The coupled stochastic-deterministic homogenization
problem in the Fourier-Stieltjes algebra}. We first need to define the
Fourier-Stieltjes algebra $FS(\mathbb{R}^{N})$: The Fourier-Stieltjes algebra%
\emph{\ }on $\mathbb{R}^{N}$\ is defined as the closure in $\mathcal{B}(%
\mathbb{R}^{N})$\ of the space 
\begin{equation*}
FS_{\ast }(\mathbb{R}^{N})=\left\{ f:\mathbb{R}^{N}\rightarrow \mathbb{R}%
,\;f(x)=\int_{\mathbb{R}^{N}}\exp (ix\cdot y)d\nu (y)\text{\ for some }\nu
\in \mathcal{M}_{\ast }(\mathbb{R}^{N})\right\}
\end{equation*}%
where $\mathcal{M}_{\ast }(\mathbb{R}^{N})$\ denotes the space of complex
valued measures $\nu $\ with finite total variation: $\left\vert \nu
\right\vert (\mathbb{R}^{N})<\infty $. We denote it by $FS(\mathbb{R}^{N})$.
Since by \cite{17} any function in $FS_{\ast }(\mathbb{R}^{N})$ is a weakly
almost periodic continuous function, we have that $FS(\mathbb{R}^{N})\subset
WAP(\mathbb{R}^{N})$. Moreover thanks to \cite[Theorem 4.5]{Chou} $FS(%
\mathbb{R}^{N})$ is a proper subalgebra of $WAP(\mathbb{R}^{N})$.

As $FS(\mathbb{R}^{N})$ is an ergodic algebra which is translation invariant
(this is easily seen: indeed $FS_{\ast }(\mathbb{R}^{N})$ is translation
invariant) we see that the hypotheses of Theorem \ref{t3.3} are satisfied
with algebra $A=FS(\mathbb{R}_{y}^{N})$.

With all the above in mind, we see that one can solve the homogenization
problem for (\ref{6.1}) under the assumption: 
\begin{equation*}
a_{ij}(x,\omega ,\cdot )\in FS(\mathbb{R}^{N})\text{ and }\mathbf{h}(\omega
,\cdot )\in (AP(\mathbb{R}^{N}))^{N}.
\end{equation*}
\end{itemize}

\begin{acknowledgement}
\emph{The work of the second author is supported by the University of
Pretoria through a postdoctoral fellowship. Both authors acknowledge the
support of the National Research Foundation of South Africa through a "focus
area" grant.}
\end{acknowledgement}

\end{document}